\documentclass[12pt,a4paper,reqno, final]{amsart}

\usepackage[utf8]{inputenc}
\usepackage{amsmath, amsfonts, amsthm, amssymb, amscd}
\usepackage[foot]{amsaddr}
\usepackage{graphicx}
\usepackage{lmodern} 
\usepackage{esint}
\PassOptionsToPackage{dvipsnames}{xcolor}
\usepackage{tikz}
\usepackage{epsfig}
\usepackage{floatflt}
\usepackage[english]{babel}
\usepackage{showkeys}

\usepackage{comment}
\usepackage[T1]{fontenc}       
\usepackage{url}
\usepackage{color}
\usepackage{setspace}
\usepackage{nicefrac}
\usepackage[colorlinks, linkcolor = black, citecolor = black, filecolor = black, urlcolor = blue]{hyperref}
\usepackage{esint} 
\usepackage{fancyref}
\usepackage{ifthen}
\usepackage{xcolor}
\usepackage{todonotes}
\usepackage{commath}
\usepackage{bookmark}
\allowdisplaybreaks[2]

\begin{document}
	\numberwithin{equation}{section}
	\newcommand{\proofend}{\begin{flushright}$\Box$\end{flushright}}
	\renewcommand{\proof}{\textbf{Proof}: }
	\newcommand{\R}{\mathbb{R}}
	\newcommand{\C}{\mathbb{C}}
	\newcommand{\N}{\mathbb{N}}
	\newcommand{\Z}{\mathbb{Z}}
	\newcommand{\Suml}{\sum\limits}
	\newcommand{\Cupl}{\bigcup\limits}
	\newcommand{\Capl}{\bigcap\limits}
	\newcommand{\Intl}{\int\limits}
	\newcommand{\Liml}{\lim\limits}
	\newcommand{\supl}{\sup\limits}
	\newcommand{\infl}{\inf\limits}
	\newcommand{\rL}{\mathrm{L}}
	\newcommand{\cL}{\mathcal{L}}
	\newcommand{\cW}{\mathcal{W}}
	\newcommand{\cD}{\mathcal{D}}
	\newcommand{\sS}{\mathcal{S}}
	\newcommand{\cO}{\mathcal{O}}
	\newcommand{\cX}{\mathcal{X}}
	\newcommand{\cH}{\mathcal{H}}
	\newcommand{\cG}{\mathcal{G}}
	\newcommand{\BV}{\mathrm{BV}}
	\newcommand{\rP}{\mathrm{P}}
	\newcommand{\divv}{\mathrm{div}\,}
	\newcommand{\cE}{\mathcal{E}}
	\newcommand{\M}{\mathcal{M}}
	\newcommand{\B}{\mathcal{B}}
	\newcommand{\X}{\mathcal{X}}
	\newcommand{\F}{\mathcal{F}}
	\newcommand{\e}{\varepsilon}
	\newcommand{\z}{\zeta}
	\newcommand{\s}{\sigma}
	\newcommand{\g}{\gamma}
	\newcommand{\ca}{\operatorname{cap}}
	\newcommand{\diam}{\operatorname{diam}}
	\newcommand{\maxi}{\operatorname{max}}
	\newcommand{\mini}{\operatorname{min}}
	\newcommand{\rad}{\operatorname{rad}}
	\newcommand{\dist}{\operatorname{dist}}
	
	\newcommand{\dx}{\dif x}
	\newcommand{\dy}{\dif y}
	\newcommand{\dt}{\dif t}
	\newcommand{\ds}{\dif s}
	\newcommand{\dbl}{\dif \bl}
	\newcommand{\dble}{\dif \ble}
	\newcommand{\dblt}{\dif \blt}
	\newcommand{\dblh}{\dif \blh}
	\newcommand{\dblet}{\dif \blet}
	\newcommand{\dbleh}{\dif \bleh}
	\newcommand{\inttimet}[1]{\int_{0}^{t} #1 \, \dif s}
	\newcommand{\inttimeT}[1]{\int_{0}^{T} #1 \, \dif s}
	\newcommand{\inttimetra}[1]{\int_{0}^{t} #1 \, \dif \beta_{k}}
	\newcommand{\intsptit}[1]{\int_{0}^{t}\int_{\cO} #1 \, \dif x \dif s}
	\newcommand{\intsptist}[1]{\int_{s}^{t}\int_{\cO} #1 \, \dif x \dif s}
	\newcommand{\intsptiT}[1]{\int_{0}^{T}\int_{\cO} #1 \, \dif x \dif s}
	\newcommand{\intsptitra}[1]{\int_{0}^{t}\int_{\cO} #1 \, \dif x \dif \beta_{k}^{\e}}
	\newcommand{\intsptiTra}[1]{\int_{0}^{T}\int_{\cO} #1 \, \dif x \dif \beta_{k}^{\e}}
	
	\newcommand{\intEwert}[1]{\mathbb{E} \left[ \int_{0}^{t}\int_{\cO} #1 \, \dif x \dif s \right]}
	
	\newcommand{\betr}[1]{\left| #1 \right|}
	\newcommand{\lpr}[2]{L^{#1}(#2)}
	\newcommand{\Hsob}[2]{H^{#1}(#2)}
	\newcommand{\Hsobper}[2]{H_{\text{per}}^{#1}(#2)}
	\newcommand{\sob}[2]{W^{#1}(#2)}
	\newcommand{\li}[2]{\underset{#1 \rightarrow #2 }{\lim}}
	\newcommand{\Ew}[1]{\mathbb{E}\left[  #1  \right]}
	\newcommand{\Ewt}[1]{\ti{\mathbb{E}}\left[  #1  \right]}
	\newcommand{\prob}[1]{\mathbb{P}\left[  #1  \right]}
	\newcommand{\ti}[1]{\tilde{#1}}
	
	\newcommand{\ut}{\ti{u}}
	\newcommand{\pt}{\ti{p}}
	\newcommand{\pe}{{p^{\e}}}
	\newcommand{\pet}{\ti{p}^{\e}}
	\newcommand{\pex}{p_{x}^{\e}}
	\newcommand{\Je}{J^{\e}}
	\newcommand{\Jet}{\tilde{J}^{\e}}
	\newcommand{\Jt}{\tilde{J}}
	\newcommand{\uh}{\hat{u}}
	\newcommand{\ph}{\hat{p}}
	\newcommand{\peh}{\hat{p}^{\e}}
	\newcommand{\Jeh}{\hat{J}^{\e}}
	\newcommand{\Jh}{\hat{J}}
	\newcommand{\ue}{u^{\e}}
	\newcommand{\us}{u_{\sigma}}
	\newcommand{\usx}{(u_{\sigma})_{x}}
	\newcommand{\usxx}{(u_{\sigma})_{xx}}
	\newcommand{\uast}{u_{\sigma}^{\ast}}
	\newcommand{\uen}{u^{\e}_{0}}
	\newcommand{\ua}{u^{\alpha}}
	\newcommand{\uam}{u^{\alpha-1}}
	\newcommand{\uap}{u^{\alpha+1}}
	\newcommand{\ux}{u_{x}}
	\newcommand{\uxx}{u_{xx}}
	\newcommand{\uxxx}{u_{xxx}}
	\newcommand{\uex}{\ue_{x}}
	\newcommand{\uexx}{\ue_{xx}}
	\newcommand{\uexxx}{\ue_{xxx}}
	\newcommand{\un}{u_{0}}
	\newcommand{\uea}{(\ue)^{\alpha}}
	\newcommand{\ueam}{(\ue)^{\alpha-1}}
	\newcommand{\ueap}{(\ue)^{\alpha+1}}
	\newcommand{\lal}{\lambda_{k}}
	\newcommand{\lan}{\lambda_{0}}
	\newcommand{\laml}{\lambda_{-k}}
	\newcommand{\ssl}{\sigma_{s}^{k}}
	\newcommand{\Cstr}{C_{\text{Strat}}}
	\newcommand{\xiast}{\xi^{\ast}}
	\newcommand{\fast}{f^{\ast}}
	
	\newcommand{\uet}{\ti{u}^{\e}}
	\newcommand{\ps}{p_{\sigma}}
	\newcommand{\psx}{(\ps)_{x}}
	\newcommand{\uent}{\ti{u}^{\e}_{0}}
	\newcommand{\uat}{\ti{u}^{\alpha}}
	\newcommand{\uamt}{\ti{u}^{\alpha-1}}
	\newcommand{\uapt}{\ti{u}^{\alpha+1}}
	\newcommand{\uxt}{\ti{u}_{x}}
	\newcommand{\uxh}{\hat{u}_{x}}
	\newcommand{\uxxt}{\ti{u}_{xx}}
	\newcommand{\uxxh}{\hat{u}_{xx}}
	\newcommand{\uext}{\uet_{x}}
	\newcommand{\uexh}{\ueh_{x}}
	\newcommand{\uexxt}{\uet_{xx}}
	\newcommand{\uexxh}{\ueh_{xx}}
	\newcommand{\unt}{\ti{u}_{0}}
	\newcommand{\unh}{\hat{u}_{0}}
	\newcommand{\ueat}{(\uet)^{\alpha}}
	\newcommand{\ueamt}{(\uet)^{\alpha-1}}
	\newcommand{\ueapt}{(\uet)^{\alpha+1}}
	\newcommand{\Gas}{G_{\alpha }}
	\newcommand{\Gmcas}{\mathcal{G}_{\alpha }}
	\newcommand{\wt}{\ti{w}}
	\newcommand{\vt}{\ti{v}}
	\newcommand{\zt}{\ti{z}}
	\newcommand{\ve}{v^{\e}}
	\newcommand{\ze}{z^{\e}}
	\newcommand{\vet}{\ti{v}^{\e}}
	\newcommand{\wet}{\ti{w}^{\e}}
	\newcommand{\zet}{\ti{z}^{\e}}
	\newcommand{\vext}{\ti{v}^{\e}_{x}}
	\newcommand{\zext}{\ti{z}^{\e}_{x}}
	\newcommand{\vexxt}{\ti{v}^{\e}_{xx}}
	\newcommand{\zetxx}{\ti{z}^{\e}_{xx}}
	
	\newcommand{\epsw}[1]{\e\, \cW(#1)}
	\newcommand{\epsws}[1]{\e\, \cW'(#1)}
	\newcommand{\epswss}[1]{\e\, \cW''(#1)}
	\newcommand{\gl}{g_{k}}
	\newcommand{\gml}{g_{-k}}
	\newcommand{\glx}{(g_{k})_{x}}
	\newcommand{\glxx}{(g_{k})_{xx}}
	
	\newcommand{\bl}{\beta_{k}}
	\newcommand{\blt}{\ti{\beta}_{k}}
	\newcommand{\blh}{\hat{\beta}_{k}}
	\newcommand{\ble}{\beta^{\e}_{k}}
	\newcommand{\blet}{\ti{\beta}^{\e}_{k}}
	\newcommand{\W}{F}
	\newcommand{\We}{W^{\e}}
	\newcommand{\Wt}{\ti{W}}
	\newcommand{\Wet}{\ti{W}^{\e}}
	\newcommand{\Wpq}{\mathcal{W}_{pq\text{Strat}}}
	\newcommand{\WStrat}{F_{\text{Strat}}}
	\newcommand{\Tmax}{T}
	\newcommand{\Ot}{\cO_{T}}
	\newcommand{\Otb}{\bar{\cO}_{T}}
	
	\newcommand{\Norm}[1]{\lVert \,  #1 \, \rVert }
	\newcommand{\Normsq}[1]{\lVert \, #1 \, \rVert^{2} }
	\newcommand{\skp}[2]{\left( #1  ,  #2  \right)}
	\newcommand{\iprod}[2]{\langle \,  #1 \, , \,#2 \, \rangle}
	
	\newcommand{\abl}[1]{\frac{\partial}{\partial #1}}
	\newcommand{\ablsq}[1]{\frac{\partial^{2}}{\partial #1^{2}}}
	\newcommand{\ablcu}[1]{\frac{\partial^{3}}{\partial #1^{3}}}
	
	\newcommand{\intt}{{\int_0^{t}}}
	\newcommand{\intT}{{\int_0^{T}}}
	\newcommand{\intsp}{\int_{\cO}}
	\newcommand{\intTs}{{\int_0^{T\wedge T_{\sigma}}}}
	\newcommand{\intts}{{\int_0^{t \wedge T_{\sigma}}}}
	
	\newcommand{\stochbas}{(\Omega, \mathcal{F}, (\mathcal{F}_{t})_{t \ge 0},\mathbb{P})}
	\newcommand{\stochbasti}{(\ti{\Omega}, \ti{\mathcal{F}}, (\ti{\mathcal{F}_{t}})_{t \in [0,T]},\ti{\mathbb{P}})}
	\newcommand{\qv}[1]{\langle #1 \rangle}
	\newcommand{\crossv}[2]{\langle\langle #1,#2 \rangle\rangle}
	
	\newtheorem{sat1}{Theorem}[section]
	\newtheorem{def1}[sat1]{Definition}
	\newtheorem{the1}[sat1]{Theorem}
	\newtheorem{lem1}[sat1]{Lemma}
	\newtheorem{kor1}[sat1]{Corollary}
	\newtheorem{bsp1}[sat1]{Example}
	\newtheorem{bem1}[sat1]{Remark}
	\newtheorem{prop1}[sat1]{Proposition}
	
	\newenvironment{myproof}[1]%
	{\paragraph{\textbf{Proof{#1}:}}}
	{\hfill$\square$\newline} 
	
	\def\Xint#1{\mathchoice
		{\XXint\displaystyle\textstyle{#1}}%
		{\XXint\textstyle\scriptstyle{#1}}%
		{\XXint\scriptstyle\scriptscriptstyle{#1}}%
		{\XXint\scriptscriptstyle\scriptscriptstyle{#1}}%
		\!\int}
	\def\XXint#1#2#3{{\setbox0=\hbox{$#1{#2#3}{\int}$ }
			\vcenter{\hbox{$#2#3$ }}\kern-.6\wd0}}
	\def\ddashint{\Xint=}
	\def\dashint{\Xint-}
	
	\newcommand{\expect}{{\mathbb E}}
	\newcommand{\fbeta}{{F_\beta}}
	\newcommand{\vdm}{{|v_\delta|^{m-1}}}
	\newcommand{\meanball}{{\dashint_{B_\delta(x_0)}}}
	\newcommand{\Rball}{{B_R(x_0)}}
	\newcommand{\rball}{{B_r(x_0)}}
	\newcommand{\ort}{{\mathcal O}}
	\newcommand{\intort}{{\int_\ort}}
	\newcommand{\filt}{{\mathcal F_t}}
	\newcommand{\fils}{{\mathcal F_s}}
	\newcommand{\rhod}{{\rho_\delta\ast}}
	\newcommand{\fdu}{{F_\delta[u]}}
	\newcommand{\probab}{{\mathbb P}}
	\newcommand{\tprobab}{{\tilde\probab}}
	\newcommand{\integrand}{{|u|^{2m}(x,t)dxdt}}
	\newcommand{\intte}{{\int_0^{T_E}}}
	
	\newcommand{\Ih}{\mathcal{I}_h}
	\newcommand{\domain}{\mathcal{O}}
	\newcommand{\Lh}{{L_h}}
	
	\title[]{Existence of nonnegative energy-dissipating solutions to a class of stochastic thin-film equations under weak slippage: Part II -- compactly supported initial data}
	\date{June 11, 2024}
	
	\author[G.~Gr\"un]{G\"unther Gr\"un}
	\address[G.~Gr\"un]{Friedrich--Alexander--Universität Erlangen--Nürnberg, Cauerstraße 11, 91058 Erlangen,~Germany}
	\email{gruen@math.fau.de}
	
	\author[L.~Klein]{Lorenz Klein}
	\address[L.~Klein]{Friedrich--Alexander--Universität Erlangen--Nürnberg, Cauerstraße 11, 91058 Erlangen, Germany}
	\email{lorenz.klein@fau.de}

	\keywords{stochastic thin-film equation; stochastic partial differential equation; conservative nonlinear multiplicative noise; compactly supported initial data}
	\makeatletter
	\@namedef{subjclassname@2020}{%
		\textup{2020} Mathematics Subject Classification}
	\makeatother 
	\subjclass[2020]{60H15, 76A20, 35G20, 35Q35, 35K65, 35R35, 35R37} 
	
	
	\begin{abstract}
		We prove existence of martingale solutions to a class of stochastic thin-film equations for mobility exponents $n \in (2,3)$ and compactly supported initial data.   	
		With the perspective to study free boundary problems related to stochastic thin-film equations in \cite{GruenKlein24FSOP}, we start from the surface-tension driven stochastic thin-film equation with Stratonovich noise and exploit the regime of coefficients in front of the Stratonovich correction term (which is of porous-media-type) for which energy dissipation can be established. By Bernis inequalities, third order spatial derivatives of appropriate powers of the solution are controlled. Analytically, we rely on approximation with $\mathbb{P}$-almost surely strictly positive solutions and compactness methods based on energy-entropy estimates as well as martingale identification techniques. 
	\end{abstract}	
	
	\maketitle

	\section{Introduction}\label{sec:Intro}
	
	This paper is the second in a series of three publications, aiming to prove existence and finite speed of propagation of energy-dissipating martingale solutions for a class of stochastic thin-film equations
	\begin{align}\label{Intro:stfe-Strat-equ1} 
		\dif u + (u^{n}\uxxx)_{x} \dt  
		- 
		(\Cstr+S) (u^{n-2}\ux)_{x} \dt
		= \sum_{k \in \Z} \lal (\gl u^{n/2})_{x} \dbl(t)  
	\end{align}	
	on a bounded spatial domain $\cO := (0,L)$ subject to periodic boundary conditions. 
	The constants $\Cstr$ and $S$ as well as assumptions on the noise will be specified in Section~\ref{sec:Prel}. 
	
	Our main goal here is to prove the existence of energy-dissipating nonnegative solutions to \eqref{Intro:stfe-Strat-equ1} with compactly supported initial data that have regularity properties -- in particular control of first and third order derivatives of certain powers of $u$ -- suitable to prove results on their qualitative behavior such as finite speed of propagation. The latter will be the content of the last paper in this series \cite{GruenKlein24FSOP}. 
	 
	Let us give a brief overview of the literature on the stochastic thin-film equation.
	For a more exhaustive overview of the mathematics of thin-film flow, we refer to the first section of \cite{GruenKlein24Exist} and the references therein.  
	Stochastic versions of thin-film equations have been derived in \cite{DavidovitchStone2005, MeckeRauscher} to investigate the influence of thermal fluctuations on spreading or dewetting of very thin liquid films -- see also \cite{Dur_n_Olivencia_2019} for more recent numerical studies.
	A first existence result for stochastic thin-film equations was given by Fischer and the first author \cite{FischerGruen}.
	In this paper, the case of dewetting was studied for $n = 2$, leading to conservative linear multiplicative noise, with the physical energy given by the sum of a surface tension term and a singular conjoining/disjoining effective interface potential. 
	Purely surface-tension driven flow was studied by Gess and Gnann in \cite{GessGnann2020}. 
	They were the first to study a thin-film equation with Stratonovich noise and they showed existence of nonnegative solutions for $n = 2$, in \cite{DGGG} the range of mobility exponents was extended to $[\frac{8}{3},4)$, however, allowing only initial data with full support. 
	In two space dimensions there are corresponding results to \cite{FischerGruen} and \cite{GessGnann2020} by Metzger and the first author \cite{MetzgerGruen2022} and by Sauerbrey \cite{Sauerbrey2021} available. 
	
	We note also the recent contributions \cite{Dareiotis2023solutions} and \cite{Sauerbrey2024solutions}, where $n \in (2,3)$ and initial data may be compactly supported. 
	These papers follow an approach based solely on $\alpha$-entropy estimates, extending the deterministic estimate
	\begin{align}\label{Intro:Entropy} \notag
			\frac{1}{\alpha(\alpha+1)}\intsp u^{\alpha+1}(t) \dx 
			&+
			c\intt \intsp \left|u^{\frac{n+\alpha+1}{4}}\right|_{x}^{4} 
			+
			\left|u^{\frac{n+\alpha+1}{2}}\right|_{xx}^{2} \dx \ds
			\\
			&\le
			\frac{1}{\alpha(\alpha+1)} \intsp \un^{\alpha+1} \dx \, ,
	\end{align} 
	with $\alpha \in (\frac{1}{2} -n ,2-n) \backslash \{0,-1\}$, cf. also \cite{Beretta1995} and \cite{Bertozzi1996}, to the stochastic setting.
	However, without suitable energy estimates in the spirit of 
	\begin{align}\label{Intro:energy}
		\frac{1}{2} \intsp u_{x}^{2}(t)\dx + \intt \intsp u^{n} u^{2}_{xxx} \dx \ds \le \frac{1}{2} \intsp (\un)_{x}^{2} \dx  
	\end{align}
	at hand,  there is even in the deterministic setting no strategy known
	to show finite speed of propagation in the regime $n \in [2,3)$.\footnote{For $n<2$, the situation is even worse as the Stratonovich correction term in this setting is given by a fast diffusion operator.} 

	With energy estimates, the picture changes completely. In the deterministic setting, there are various approaches -- based on a combination  with differential inequalities \cite{BernisFinite2} or with Stampacchia's lemma \cite{GruenDropletWTP}, see also \cite{DalPassoGiacoGruenWTP} -- to establish finite speed of propagation for $n \in [2,3)$. An important ingredient are (weighted versions of) Bernis' inequalities \cite{BernisInequ1996} (or their multi-dimensional generalization \cite{GruenBernisInequ}, respectively) which read in one dimension as 
	\begin{align}\label{Intro:BernisBasic}
		\intort u^{n-4} \ux^{6} \dx + \intort u^{n-1} |\uxx|^{3} \dx + \intort \left( u^{\frac{n+2}{2}}\right)^{2}_{xxx} \dx \le C \intort u^{n} \uxxx^{2} \dx \, .
	\end{align} 
	For qualitative results on finite speed of propagation in the regime $0 < n<2$, we refer to \cite{BernisFinite} and \cite{BertschDalPassoGarckeGruen}.
	
	In the stochastic setting, Fischer and the first author presented in \cite{FischerGruenFSOP} a recipe to prove finite speed of propagation for stochastic porous-medium equations with linear multiplicative source-type noise by transforming methods from the deterministic setting
	involving Stampacchia's lemma to the stochastic one. 
	In \cite{GruenKlein24FSOP}, we generalize this approach to energy-dissipating solutions of \eqref{Intro:stfe-Strat-equ1} with conservative nonlinear multiplicative noise.
	
	The present work and its predecessor \cite{GruenKlein24Exist} are devoted to establish the existence of the solutions needed in \cite{GruenKlein24FSOP}.
	Especially in this paper, solutions to compactly supported initial data are constructed as singular limits of appropriate strictly positive solutions. 
	Starting point are weak martingale solutions with singular effective interface potential to 
	\begin{align}\label{Intro:stfe-Strat-equ} 
		\dif \ue +  ((\ue)^n(\uexx - \e \W'(\ue))_{x})_{x} \dt  
		- 
		(\Cstr+S) ((\ue)^{n-2}\uex)_{x} \dt
		=
		\sum_{k \in \Z} \lal (\gl(\ue)^{\tfrac{n}{2}})_{x} \dble(t) \,   ,
	\end{align}
	subject to periodic boundary conditions and strictly positive initial data $\un^{\delta} = \un + \delta$. 
	Here, $F(u):=u^{-p}$, $\e \in (0,1)$.
	Such solutions have been shown to exist in \cite{GruenKlein24Exist}. 
	In this paper, we pass successively with $\e$ and then with $\delta$ to the limit zero.
	\\
	For the passage to vanishing potential -- i.e. $\e \to 0$ -- it is important to note that the parameter $\delta$, and hence initial data, remain untouched. 
	This way, we may use entropy estimates like \eqref{Intro:Entropy} for the case $\alpha = 1-n$. The additional regularity that comes along with this estimate has two important consequences. First, the control of a negative power of approximating solutions shows that sets where the limit may touch down have measure zero, cf. Remark~\ref{bem:ut-nullset}. Secondly, control of the Laplacian of solutions leads to spatially continuously differentiable solutions. In turn, we may apply a suitable version of Bernis' inequalities \cite{GiacShish2005} that goes back to \cite{GiacOtto2002}, to infer the aforementioned additional integral estimates, see Proposition~\ref{prop:Bernis}. 
	These estimates persist when passing to the limit of nonnegative initial data $\delta \to 0$ in the second part of this work.
	\\ 
	
	The outline of the paper is as follows. 
	In Section~\ref{sec:Prel}, we fix the setting and we  formulate general assumptions on initial data, noise, and parameters. 
	
	In Section~\ref{sec:LimitPot}, we discuss the passage to the limit $\e \to 0$ of solutions to \eqref{Intro:stfe-Strat-equ}. We divide this section into three subsections. After having formulated the main result of this section in Theorem~\ref{theo:MainResult}, we are concerned with the derivation of an energy-entropy estimate based on It\^o's formula and a stopping time argument. Next, we justify the application of Jakubowski's generalization of the Skorokhod convergence theorem, using the a-priori estimates to prove tightness on suitable path spaces. This way, we pass to the limit on a new probability space $(\tilde{\Omega}, \tilde{\mathcal{F}}, \tilde{\mathbb{P}})$ in the deterministic terms (i.e. the left-hand side of \eqref{Intro:stfe-Strat-equ}).
	Thirdly, we identify the stochastic integral in the limit. 
	Here, we apply techniques introduced in \cite{HofmanovaSeidler2012} and \cite{Brezniak2007} (for applications see for e.g. \cite{BreitHofmanova2016,DebusscheHofmanovaVovelle,HofmanovaWeak2016}) which provide a convenient way of martingale representation based on elementary calculations involving the quadratic variation process of the right-hand side of \eqref{Intro:stfe-Strat-equ} and the quadratic covariation processes of it with respect to standard Brownian motions.    
	Finally, we give the proof of Theorem~\ref{theo:MainResult} as well as the result on additional regularity in Proposition~\ref{prop:Bernis}. 
	Due to its technical nature and hence for the sake of readability, one assertion regarding positivity of solutions is proven in Appendix~\ref{app:Pos}.
	 
	In Section~\ref{sec:LimitInit}, we will pass to the limit $\delta \to 0$, i.e. $\un^{\delta} \to \un$, where $\un$ is nonnegative. We state our main result Theorem~\ref{theo:ConvS} which claims that in the limit we recover a nonnegative energy-dissipating martingale solution of the stochastic thin-film equation that satisfies additional regularity as needed for the proof of finite speed of propagation in \cite{GruenKlein24FSOP}.
	The general approach is the same as in Section~\ref{sec:LimitPot}.  
	As the technical details have many common aspects, we will be rather brief in this section. 
	  
	There are three appendices. Appendix~\ref{app:Aux} summarizes the results on existence of positive approximating solutions to \eqref{eq:stfe-Strat-equ} from \cite{GruenKlein24Exist}. Moreover, we introduce some useful estimates of powers of $\ue$, and we provide a strictly positive lower bound of solutions $\ue$. 
	In Appendix~\ref{appendix:Itoformula}, we justify the application of It\^o's formula as needed in the proof of Theorem~\ref{lem:EnEnt} and Lemma~\ref{lem:EnalphaEnt1}. 
	Finally in Appendix~\ref{app:Pos}, we derive an $\alpha$-entropy estimate and we prove the positivity properties of solutions $\ut$ formulated in Theorem~\ref{theo:MainResult}, (1). 
	\\
	
	\noindent
	\textbf{Notation:} 
	We consider a spatial domain $\cO := (0,L)$ and define the space-time cylinder $\cO_{T} := (0,L) \times [0,T]$.
	By $\rball$ we denote open balls with radius $r>0$ and center $x_{0} \in \cO$.
	The notation $a \wedge b$ stands for the minimum of $a$ and $b$.
	\\
	Regarding function spaces, we use standard notation for Lebesgue and Sobolev spaces.
	For $\gamma, \gamma' \in (0,1)$, $C^{\gamma, \gamma'}(\bar{\cO}_{T})$ denotes the space of functions that are Hölder-continuous in space (respectively, time) with Hölder coefficient $\gamma$  (respectively, $\gamma'$).
	The notation $C_{0}^{\infty}(\cO)$ is used for the class of smooth functions  on $\cO$ that are compactly supported. 
	Whenever we consider functions that are periodic in space, we add the subscript 'per' to the corresponding function space. 
	\\
	The quadratic variation process of a stochastic process $X$ is denoted by $\langle X \rangle$. $\langle X,Y \rangle$ is the quadratic covariation process of $X$ with another process $Y$.
	For two separable Hilbert spaces $U$ and $H$ we denote by $L_{2}(U,H)$ the space of Hilbert-Schmidt operators from $U$ to $H$, i.e. the space of linear bounded operators that satisfy $\Norm{A}_{L_{2}(U,H)}^{2} := \sum_{k\in \Z} \Norm{A e_{k}}^{2}_{H} < \infty$, where $(e_{k})_{k\in \Z}$ is an orthonormal basis on $U$. We abbreviate $L^{0}_{2}(U,H) := L_{2}(Q^{\frac{1}{2}}U,H)$ with the linear operator $Q$ defined in (A2) below.
	
\section{Preliminaries}\label{sec:Prel}
	For $F(u):=u^{-p}$ and $\e \in (0,1)$, we consider on $\Ot := [0,L] \times [0,T]$ the equations
	\begin{align}\label{eq:stfe-Strat-equ} \notag
		\dif \ue =  &- ((\ue)^n(\uexx - \e \W'(\ue))_{x})_{x} \dt  
		+ 
		(\Cstr+S) ((\ue)^{n-2}\uex)_{x} \dt
		\\
		& + \sum_{k \in \Z} \lal (\gl(\ue)^{\tfrac{n}{2}})_{x} \dble(t) \,   
	\end{align}
	and
	\begin{align}\label{eq:stfe-Strat-equ1} 
		\dif u =  -(u^{n}\uxxx)_{x} \dt  
		+ 
		(\Cstr+S) (u^{n-2}\ux)_{x} \dt
		+ \sum_{k \in \Z} \lal (\gl u^{\tfrac{n}{2}})_{x} \dbl(t)  
	\end{align}	
	subject to periodic boundary conditions and the following assumptions:
	\begin{itemize}
		\item[(A1)]
		$n \in (2,3)$.
		\item[(A2)]
		Let $\stochbas$ be a stochastic basis with a complete right-continuous filtration. 
		$W$ is an $\mathcal{F}_{t}$-adapted $Q$-Wiener process on $\Omega$ defined by an operator $Q\colon\lpr{2}{\cO} \rightarrow \lpr{2}{\cO}$ with $Q\gl = \lal^{2}\gl$, $k \in \Z$. Here $(\lal)_{k \in \Z}$ are nonnegative numbers and $(\gl)_{k \in \Z}$ are eigenfunctions of the Laplacian on $\cO$ subject to periodic boundary conditions,
		\begin{align}\label{def:L2-Basis}
			\gl(x) = 
			\begin{cases}
				\sqrt{\frac{2}{L}} \sin(\frac{2 \pi k x }{L})\quad &k >0,
				\\
				\frac{1}{\sqrt{L}}\quad &k=0, 
				\\
				\sqrt{\frac{2}{L}} \cos(\frac{2 \pi k x }{L})\quad &k<0.
			\end{cases}
		\end{align} 
		The noise is colored in the sense that $\sum_{k \in \Z} k^{4} \lal^{2} < \infty$. 
		Moreover, we assume $\lal = \laml$ for all $k\in \N$. 
		\item[(A3)]
		$S$ satisfies $S > C_{Strat} \frac{3}{2^{n-4}}  \frac{n-2}{3-n} $,
		where $\Cstr = \frac{1}{2} \frac{n^{2}}{4} \left(\frac{\lan^{2}}{L} + \sum_{k=1}^{\infty} \frac{2\lal^{2}}{L}\right)$.
		\item[(A3)*]
		$S$ satisfies
		$S > \Cstr \frac{9}{4}\,\frac{(n-2)^{2}}{(3-n)(2n-3)}$, $\Cstr$ as above. 
		\item[(A4)]
		Let $\Lambda^{0}$ be a probability measure on $H^{1}_{\text{per}}(\cO)$ equipped with the Borel $\sigma$-algebra which is supported on the subset of nonnegative functions such that there is a positive constant $C$ with the property that 
		\begin{align*}
			\operatorname{esssup}_{v \in \operatorname{supp \Lambda^{0}}} 
			\left\{		
			\int_{\cO} \frac{1}{2} \abs{v_{x}}^{2}\dx +  \int_{\cO} v \dx    	
			\right\}
			\le C \, .
		\end{align*}
		We assume there exists an $\mathcal{F}_{0}$-measurable random variable $\un$ such that  $\Lambda^{0} = \mathbb{P} \circ \un^{-1}$.
		
		\item[(A4$^{\delta}$)]
		Let $\Lambda^0$ be a probability measure on $H^1_{\text{per}}(\ort)$ satisfying hypothesis (A4). Then, $\Lambda^{\delta}$ is the probability measure on $H^1_{\text{per}}(\ort)$ defined by 
		$\Lambda^{\delta} = \Lambda^{0} \circ S_{\delta}^{-1}$, where \newline $S_{\delta}\colon H^{1}_{\text{per}}(\cO) \rightarrow H^{1}_{\text{per}}(\cO)$ is  
		given by $S_{\delta}(u) = u + \delta$. 
		
	\end{itemize}
	Moreover, we set $\pe := - \uexx + \e \W'(\ue)$ and $\un^{\delta} := \un + \delta$.
	
	\begin{bem1}\label{rem:InitalData}
		Initial data as in (A$4^{\delta}$) also satisfy (H3) in \cite{GruenKlein24Exist}. 
	\end{bem1}
	\begin{bem1}\label{rem:lowerboundS}
		The lower bound of $S$ in (A3) is due to the construction of our approximating functions in Theorem~3.2 in \cite{GruenKlein24Exist}, cf. Definition~\ref{DefinitionMartingaleSolution} and Theorem~\ref{MainResultFG}. There it was needed to account for discrete integration by parts formulae. Formally, in our analysis here the weaker assumption (A3)* is sufficient, cf. the proof of Lemma~\ref{lem:EnEnt}.
		Note that equation \eqref{Intro:stfe-Strat-equ1} is the surface-tension driven stochastic thin-film equation with Hänggi-Klimontovich (or backward It\^o) noise (cf. \cite{Karatzas2014, Kuo2006}) provided $S$ can be chosen equal to $\Cstr$. This is possible for values of $n$ sufficiently close to $2$.
	\end{bem1}
	
	\section{Passage to the limit $\e \rightarrow 0$}\label{sec:LimitPot}
	In this section, we will carry out the first step in our construction of compactly supported solution to \eqref{eq:stfe-Strat-equ1}.
	We start with solutions to \eqref{eq:stfe-Strat-equ} and let $\e$ tend to zero to remove the influence of the stabilizing potential. 
	Subsequently, we prove additional regularity properties of these solutions in Proposition~\ref{prop:Bernis}.   
	We emphasize that in this section the parameter $\delta >0$ remains untouched, i.e. $\delta$ is a uniform lower bound for initial data.
	\\
	Technically, we rely on an energy-entropy estimate, see Lemma~\ref{lem:EnEnt}. To ensure that It\^o's formula is applicable, we impose a suitable stopping time argument -- the details can be found in the appendix, cf. Lemma~\ref{lem:appItoformula}. 
	Using tightness of the laws of the approximating sequences on certain path spaces, cf. Lemma~\ref{sat:Strafheit},
	we may apply Jakubowski's theorem \cite{Jakubowski1994}. On the expense of a change of the underlying probability space, this gives convergence of the sequences under consideration on the chosen path spaces almost surely, see Lemma~\ref{sat:JakubowskiAnwendung}. 
	Recovering equation \eqref{eq:stfe-Strat-equ1} is -- besides standard arguments of pde-analysis -- based on  
	techniques for martingale representation, going back to \cite{BreitHofmanova2016,Brezniak2007,DebusscheHofmanovaVovelle,HofmanovaSeidler2012}.
	The plan of the proof has some features in common with ideas presented in  \cite{FischerGruen}, \cite{GruenKlein22} or \cite{GruenKlein24Exist}. 
	For the reader's convenience, we make all the results needed explicit, but many of the proofs are just cited. When the range of the mobility parameter $n \in (2,3)$ requires new insights, we present the details.
	\\ 

	Theorem~\ref{MainResultFG} guarantees the existence of a family  $(\ue_{\delta})_{\e \in (0,1)}$ of solutions to \eqref{eq:stfe-Strat-equ} in the sense of Definition~\ref{DefinitionMartingaleSolution}, if $n\in (2,3)$ and $S > C_{strat} \frac{3}{2^{n-4}} \frac{n-2}{3-n}$. For the ease of presentation, we assume that every  $\ue_{\delta}$ is defined on the same probability space $(\Omega, \mathcal{F}, \mathbb{P})$.
	\\
	We make the following definitions  
	\begin{align}\label{def:energyeps}
		E^{\e}[u] &:= \frac{1}{2} \intsp \ux^{2} \dx + \e \intsp \W(u) \dx \, ,
		\\
		\mathcal{G}[u] &:= \intort G(u) \dx \, , \label{def:entropy}
	\end{align}
	where
	\begin{align}\label{def:entropy1}
		G(u) := \int_{1}^{u} \int_{1}^{s} r^{-n} \dif r \ds = \frac{1}{(n-1)(n-2)} u^{2-n} + \frac{1}{n-1} u - \frac{1}{n-2}\, .
	\end{align}
	The main result of this section is the following existence theorem for weak martingale solutions to \eqref{eq:stfe-Strat-equ1} with uniformly strictly positive initial data which are constructed as limits of $(u^{\e}_{\delta})_{\e \in (0,1)}$ for $\e \to 0$.
	\begin{the1}\label{theo:MainResult}
		Let (A1)-(A3) and (A$4^{\delta}$) for $\delta > 0$ be satisfied.
		For $T>0$, there exist a stochastic basis $\stochbasti$ with a complete, right-continuous filtration, an $\ti{\mathcal{F}}_{t}$-adapted Q-Wiener process $\ti{W} = \sum_{k \in \Z}\lal\gl \blt$, 
		a random variable $\unt^{\delta} \in \lpr{2}{\ti{\Omega};H^{1}_{\text{per}}(\cO)}$, and
		a continuous $L^{2}(\cO)$-valued, $\ti{\mathcal{F}}_{t}$-adapted process $\ut$ contained in  $ \lpr{2}{\ti{\Omega};L^{2}(0,T;\Hsobper{1}{\cO}})  \cap $ $L^{2}(\ti{\Omega}; C^{\ti{\gamma}, \frac{\ti{\gamma} }{4}}(\Otb))$, $ 0 < \ti{\gamma}<1/2$, with the property $\tilde{u}_{xxx} \in L_{loc}^{2}(\{ \ut > 0 \})$ for $ d\tilde{\mathbb{P}}\otimes dt$-almost all $(\tilde{\omega},t) \in \tilde{\Omega} \times [0,T]$,  
		such that the following holds:
		\begin{enumerate}
			\item \label{theo:MainResult1} 
			$\tilde{\mathbb{P}}$-almost surely $\ut(t, \cdot)$ is strictly positive for almost all $t \in [0,T]$ 
			and we have  
			$
			\Lambda^{\delta} = \mathbb{P} \circ (\un^{\delta})^{-1} = \mathbb{\tilde{P}} \circ (\unt^{\delta})^{-1},
			$
			\item\label{theo:MainResult3}
			for $t \in [0,T]$ and all
			$\phi \in H^{1}_{per}(\ort)$
			\begin{align}\notag
				\label{eq:MainResultWeakForm}
				\intsp{(\ut(t)-\unt^{\delta}) \phi} \dx
				&= \intt \intort \ut^{n} \ut_{xxx} \phi_{x} \dx \ds 
				- (\Cstr + S)\intt \intsp \ut^{n-2}  \uxt \phi_{x} \dx \ds
				\\ 
				&\quad - \sum_{k \in \Z}\int_{0}^{t}\intsp{ \lal\gl\ut^{\frac{n}{2}}\phi_{x}} \dx \dblt 
			\end{align}  
			holds $\ti{\mathbb{P}}$-almost surely,
			
			\item \label{theo:MainResult4}
			for arbitrary $q \ge 1$,
			$\ut$ satisfies the estimate  
			\begin{align}\label{est:MainResult} \notag	 
				&\Ew{\sup_{t\in [0,T]} \left(\frac{1}{2}\intsp \uxt^{2}(t) \dx \right)^{q}  
				}
				+
				\Ew{\left(\intsptiT{\ut^{n} \ut_{xxx}^{2}}\right)^{q}}
				\\ \notag
				&+ 
				\Ew{\left(\intsptiT{\ut^{n-4} \uxt^{4}}\right)^{q}}
				+
				\Ew{\left(\intsptiT{\ut^{n-2}\uxxt^{2}}\right)^{q}}
				\\ 
				&+\Ew{\left(\intT \intort \ut_{xx}^{2} \dx \ds\right)^{q}}
				+
				\Ew{\left(\intT \intort \ut^{-2} \ut_{x}^{2} \dx \ds \right)^{q}}
				\le  C(\unt, \delta,q,T). 
			\end{align}
		\end{enumerate}
	\end{the1}
\begin{bem1}
	Although we implicitly assume in Theorem~\ref{theo:MainResult} that the solutions $\ue_{\delta}$ correspond to (A3), the techniques used in the proof of Theorem~\ref{theo:MainResult} also apply, if the approximating solutions $\ue_{\delta}$ only correspond to (A3)*. Recall that the basic obstacle is still the missing existence result related to (A3)*.
\end{bem1}	
 Using the extension of Bernis' inequalities as presented in \cite{GiacShish2005}, going back to earlier results in \cite{GiacOtto2002}, we obtain the following result. \color{black}
	\begin{prop1}\label{prop:Bernis}
		For solutions $\ut$ of \eqref{eq:stfe-Strat-equ1} as derived by Theorem~\ref{theo:MainResult}, we have
		\begin{align}
			\intsp \z^{6}  \betr{\left(\ut^{\frac{n+2}{6}}\right)_{x} }^{6} \dx 
			&\le
			C \int_{\{\ut>0\}} \z^{6} \ut^{n} \betr{\ut_{xxx}}^{2} \dx + C \intsp \betr{\z_{x}}^{6} \ut^{n+2}\dx
			\\
			\intsp \z^{6} \betr{\left(\ut^{\frac{n+2}{3}}\right)_{xx} }^{3} \dx &\le C \int_{\{\ut>0\}} \z^{6} \ut^{n} \betr{\ut_{xxx}}^{2} \dx + C \intsp \betr{\z_{x}}^{6} \ut^{n+2}\dx
			\\
			\intsp \z^{6} \betr{\left(\ut^{\frac{n+2}{2}}\right)_{xxx} }^{2} \dx &\le C \int_{\{\ut>0\}} \z^{6} \ut^{n} \betr{\ut_{xxx}}^{2} \dx + C \intsp \betr{\z_{x}}^{6} \ut^{n+2}\dx
		\end{align}
		for almost all $t \in [0,T]$, $\ti{\mathbb{P}}$-almost surely, and $\zeta \in W^{1,\infty}_{per}(\cO)$ .
	\end{prop1}
	The proofs of Theorem~\ref{theo:MainResult} and Proposition~\ref{prop:Bernis} will be given in Section~\ref{sec:proofTheo11} below.

	\subsection{Energy-entropy estimate}
	For better readability, for the rest of this section, we will write $\ue$ instead of $\ue_{\delta}$ for solutions to \eqref{eq:stfe-Strat-equ} with initial data $\un^{\delta}$.
	The crucial ingredient for the proof of Theorem~\ref{theo:MainResult} is an integral estimate for energy and entropy, cf. \eqref{def:energyeps} and \eqref{def:entropy}, which will be established combining It\^o's formula with a stopping time argument. 
	We consider  
	\begin{align}\label{def:stoppingTimes}
		T_{\sigma}:= \inf\{t \in [0,T] \vert E^{\e}(\ue) \ge \sigma^{-1}\}
	\end{align}
	for positive parameters $\sigma$. For a verification that these random times are indeed stopping times, we refer to the discussion in Section~5 of \cite{GruenKlein22}.
	To comply with the assumptions needed in It\^o's formula,
	we further introduce the following cut-off versions of solutions to \eqref{eq:stfe-Strat-equ}
	\begin{align}\label{def:ModifiedSolutions}
		\us(\cdot, t) :=
		\begin{cases}
			\ue(\cdot, t) \quad & t \in [0, T_{\sigma}]
			\\
			\ue(\cdot, T_{\sigma}) \quad & t \in (T_{\sigma},\Tmax] \, ,
		\end{cases}
	\end{align}	
	where we skip the index $\e$. Moreover, we recall $\pe := -\uexx + \e \W'(\ue)$ and set $\ps := -\usxx + \e \W'(\us)$.
	We note that by Lemma~\ref{lem:Minimum-ue}, 
	$
		\ue(\cdot, t) \ge \bar{C}_{p} \e^{\frac{1}{p-2}} \sigma^{\frac{2}{p-2}}
	$
	for $t \in [0,T_{\sigma}]$,
	where the constant $\bar{C}_{p}$ depends only on $p$. 
	We have the following lemma which can be proven in a similar way as Lemma~5.5 in \cite{FischerGruen}.
	\begin{lem1}\label{lem:ConvergenceStoppingTimes}
		It holds $\Liml_{\sigma \rightarrow 0} T_{\sigma} = \Tmax$  $\mathbb{P}$-almost surely.
	\end{lem1}
	
	We will show the following uniform estimate.
	\begin{lem1}\label{lem:EnEnt}
		Let $\delta > 0$ be given and
		let $\ue = \ue_{\delta}$, $\e \in (0,1)$, be a solution to  \eqref{eq:stfe-Strat-equ} as considered in Theorem~\ref{theo:MainResult}.
		There is a positive constant $C(\un,\delta,q,T,\eta)$ independent of $\e >0$ such that for $q\ge 1$
		\begin{align}
			\label{est:EnEnt} \notag	 
			&\Ew{\sup_{t\in [0,T]} E^{\e}[\ue(t)]^{q} +  \sup_{t \in [0,T]}\mathcal{G}[\ue(t)]^{q}} 
			+
			\Ew{\left(\intsptiT{(\ue)^{n} (\pex)^{2}}\right)^{q}}
			\\ \notag
			&+ 
			c_{1}^{q}
			\Ew{\left(\intsptiT{(\ue)^{n-4} (\uex)^{4}}\right)^{q}}
			+
			c_{2}^{q}
			\, \Ew{\left(\intsptiT{(\ue)^{n-2}(\uexx)^{2}}\right)^{q}}
			\\ \notag
			&+
			S^{q} \, \Ew{\left(\intsptiT{ \e \W''(\ue) (\ue)^{n-2}(\uex)^{2}}\right)^{q}}
			\\ \notag
			&+
			\Ew{\left( \intsptiT{(\uexx)^{2}}\right)^{q}}
			+
			S^{q} \,\Ew{\left( \intsptiT{(\ue)^{-2}(\uex)^{2}}\right)^{q}}
			\\ 
			&+
			\Ew{\left(\intsptiT{\e \W''(\ue) (\uex)^{2}}\right)^{q} }
			\le  C(\un, \delta,q,T,\eta) \, ,
		\end{align}
	where $c_{1} = \mu S \frac{|(n-2)(n-3)|}{3}- \eta$ and 
	$c_{2} = S + S (1-\mu) \frac{3(n-2)}{3-n} - \Cstr \frac{9(n-2)^{2}}{4(n-3)^{2}}$
	are positive for $\mu > 0$ and $\eta > 0$ chosen sufficiently small.	
	\end{lem1}
	\proof
		We will first establish an estimate like \eqref{est:EnEnt} for stopped solutions $\us$. 
		Then, we argue as follows.
		From Lemma~\ref{lem:Minimum-ue} we infer that $\ue$ is strictly positive on $\ort\times[0,T]$ $\probab$-almost surely. In combination with Lemma~\ref{lem:ConvergenceStoppingTimes}, we find that the sets $A_\sigma:=\{\omega\in\Omega:\ue(\cdot,\omega)\geq \sigma \text{ on } \ort\times[0,T]\}$ tend for $\sigma\to 0$ to $\Omega$ up to a set of measure zero. Hence, using nonnegativity of the terms on the left-hand side of \eqref{est:EnEnt} (formulated for $\us$) and monotone convergence, Theorem~\ref{lem:EnEnt} is proven.
		\\
		Let us now derive the estimate with respect to $\us$.
		Referring to equation \eqref{eq:ItoG1} in the appendix, It\^o's formula yields for $t \in [0,T]$
		\begin{align}\notag \label{eq:ItoEntropy}
			\mathcal{G}[\us(t)] =\, &\mathcal{G}[\un^{\delta}] - \intts \intsp  \us^{n}(-\usxx + \e \W'(\us))_{x} (G'(\us))_{x} \dx \ds
			\\ \notag
			&- \left(\Cstr +S\right) \intts \intsp \us^{n-2}\usx (G'(\us))_{x} \dx \ds
			\\ \notag
			&+ \frac{1}{2}  \intts \sum_{k \in \Z}\lal^{2} \intsp (\us^{n/2}\gl)_{x}^{2} G''(\us) \dx \ds
			\\ 
			&+  \intts \sum_{k \in \Z}\lal\intsp (\us^{n/2}\gl )_{x} G'(\us) \dx \dble 
			=: \mathcal{G}[\un^{\delta}] + R_{1} +\dots + R_{4} \, ,
		\end{align} 
		where $G'(\us) =- \tfrac{1}{n-1} \us^{1-n} + \tfrac{1}{n-1}$ and $G''(\us) = \us^{-n} $.
		Due to the assumptions on the initial data, cf. (A4$^{\delta}$) and (A4), we have
		\begin{align}\label{entropy:initialdata} \notag
			\Ew{\mathcal{G}(\un^{\delta})} &= \Ew{ \intort 
				\frac{1}{(n-1)(n-2)} (\un+\delta)^{2-n} + \frac{1}{n-1} (\un+{\delta}) - \frac{1}{n-2} \dx} 
			\\
			&\le C(\delta) < \infty \, 
		\end{align}
		independently of $\e$.
		
		\noindent
	For  $R_{1}$ and  $R_{2}$, we obtain
	\begin{align}
		R_{1} 
		=
		-\intts \intsp \usxx^{2} + \e \W''(\us) \usx^{2} \dx \ds \, 
	\end{align}
	and
	\begin{align}
		R_{2} 
		= -\left(\Cstr +S\right) \intts \intsp \us^{-2}\usx^{2} \dx \ds \, ,
	\end{align}
	respectively.
	
	\noindent
	Ad $R_{3}$:
	Combining (A2) and Lemma~\ref{ONBrelations} gives
	\begin{align}\notag \label{entropy1}
		R_{3}
		&=
		\frac{1}{2}   \intts \sum_{k \in \Z} \lal^{2} \intsp (\us^{\frac{n-2}{2}} \usx \gl + \us^{n/2}\glx )^{2} \us^{-n} \dx \ds
		\\
		&\le
		\Cstr \intts \intsp \us^{-2} \usx^{2}  \dx \ds + C((\lal)_{k}, T ,L) \, .
	\end{align}
	The first term on the right-hand side of \eqref{entropy1} can be absorbed in $R_{2}$.
	
	\noindent
	Thus, for the entropy we find -- taking powers, suprema, and expectation for $q\ge1$ --
	\begin{align}\notag \label{entropy2}
		&\Ew{\sup_{t\in [0,T]}\mathcal{G}[\us(t)]^{q}} 
		+\Ew{\left(\int_{0}^{T_{\s}} \intsp \usxx^{2}  \dx \ds\right)^{q}} + 	\Ew{\left(\int_{0}^{T_{\s}}\intsp \e \W''(\us) \usx^{2} \dx \ds\right)^{q}}
		\\ \notag 
		&+ S^{q}\Ew{ \left(\int_{0}^{T_{\s}} \intsp \us^{-2}\usx^{2} \dx \ds\right)^{q}}
		\\
		&\le C(\un, \delta,T,q) 
		+
		\mathbb{E} \Bigg[\left|\sup_{t\in [0,T]}\intts \sum_{k \in \Z}\lal\intsp \left(\us^{n/2}\gl \right)_{x} G'(\us) \dx \dble \right|^{q} \Bigg]\, .
	\end{align}
	To get an estimate of the second term on the right-hand side of \eqref{entropy2}, we apply the Burkholder-Davis-Gundy inequality
	and obtain for any $\eta > 0$
		\begin{align}\notag
			&\Ew{\sup_{t \in [0,T]}
				\left|
				\sum_{k \in \Z} \int_{0}^{t \wedge T_{\s}} \intsp ((\us)^{\frac{n}{2}}\gl \lal)_{x} G'(\us) \dx \dble 
				\right|^{q} } 
			\\ \notag
			&\le
			C_{BDG} C((\lal)_{k},n,L,T,q)\, \Ew{ \left(\int_{0}^{T_{\s}}
				\intort \us^{2-n} \dx\ds
				\right)^{q/2} }	
			\\ 
			&\le
			\eta \, \Ew{ \sup_{t\in [0,T]}
			\left(\intort \us^{2-n} \dx 
				\right)^{q} } + C((\lal)_{k},n,L,T,q, \eta)	\, .	
		\end{align}
		The first term may be absorbed, the second one is constant. 
		Hence, we find
		\begin{align}\notag
			&\Ew{ \sup_{t\in [0,T]}\mathcal{G}[\us(t)]^{q}} 
			+\Ew{\left(\int_{0}^{T_{\s}} \intsp \usxx^{2}  \dx \ds\right)^{q}} + \Ew{\left(\int_{0}^{T_{\s}} \intsp \e \W''(\us) \usx^{2} \dx \ds \right)^{q}}
			\\ 
			&+ S^{q}\Ew{ \left( \int_{0}^{T_{\s}}  \intsp \us^{-2}\usx^{2} \dx \ds \right)^{q}}
			\le C(\un,\delta,(\lal)_{k},n,L,T,q,\eta) \,
		\end{align}
		uniformly in $\e$.	
		Next, we apply It\^o's formula with respect to the functional \eqref{def:energyeps}, cf. Lemma~\ref{lem:appItoformula} for details, and get
		\begin{align}\notag \label{eq:ItoEnergy}
			E^{\e}[\us(t)] = &E^{\e}[\un^{\delta}]
			+
			\intts \intsp \us^{n} \psx ((\us)_{xxx} -\e (\W'(\us))_{x} ) \dx \ds
			\\ \notag
			&+
			\left(\Cstr+S\right)\intts \intsp \us^{n-2} \usx (\us)_{xxx} \dx \ds
			\\ \notag
			&-
			 \left(\Cstr+S\right)\intts \intsp \us^{n-2} \usx \e (\W'(\us))_{x} \dx \ds
			\\ \notag
			&+
			\frac{1}{2} \sum_{k \in \Z} \lal^{2}
			\intts \intsp \e \W''(\us) (\us^{\frac{n}{2}} \gl)_{x}^{2} \dx \ds
			\\ \notag
			&+
			\sum_{k \in \Z} \lal
			\intts \intsp \e \W'(\us) (\us^{\frac{n}{2}} \gl)_{x} \dx \dble
			\\ \notag
			&+
			\frac{1}{2} \sum_{k \in \Z} \lal^{2}
			\intts \intsp (\us^{\frac{n}{2}} \gl)_{xx}^{2} \dx \ds
			\\
			&+
			\sum_{k \in \Z} \lal
			\intts \intsp \usx (\us^{\frac{n}{2}} \gl)_{xx} \dx \dble
			=: E^{\e}[\un^{\delta}] + T_{1} + \dots + T_{7} \, .
		\end{align}
		We have
		\begin{align}\label{energy:initialdata}
			\Ew{E^{\e}[\un^{\delta}]} &= 	\Ew{\frac{1}{2} \intort (\un)_{x}^{2} \dx + \e \intort (\un + \delta)^{-p} \dx} 
			\le C(\delta)  \, ,
		\end{align}
		where we used $\Lambda^{\delta} = \mathbb{P} \circ (\un^{\delta})^{-1}$, cf. (A4$^{\delta}$), as well as (A4).
		Let us now discuss the terms $T_{1}, \dots , T_{7}$.
		
		\noindent
For $T_{1}$, $T_{2}$, and $T_{3}$, we get
		\begin{align}
			T_{1} = - \intts \intsp \us^{n} \psx^{2}\dx\ds \, ,
		\end{align}
		\begin{align}\notag 	\label{EnEnt2}
			T_{2}
			&=
			- (\Cstr + S) \Bigg(\intts \intsp \tfrac{|(n-2)(n-3)|}{3} (\us)^{n-4} \usx^{4} \dx\ds 
			\\
			&\qquad \qquad \qquad \qquad
			+ \intts \intsp \us^{n-2} \usxx^{2} \dx\ds\Bigg)  \, ,
		\end{align}
and
		\begin{align}\label{EnEnt3}
			T_{3} 
			=
			- (\Cstr+S)\intts \intsp\us^{n-2} \usx^{2} \e \W''(\us) \dx \ds \, ,
		\end{align}
		respectively.
		
		\noindent
		Ad $T_{4}$:
		Using Lemma~\ref{ONBrelations}, we have
		\begin{align}\notag 	\label{EnEnt4}
			& \frac{1}{2} \sum_{k \in \Z} \lal^{2}
			\intts \intsp \e\W''(\us) (\us^{\frac{n}{2}} \gl)_{x}^{2} \dx \ds
			\\ 
			&=
			\Cstr
			\intts \intsp \e \W''(\us) \us^{n-2} \usx^{2}  \dx \ds 
			+
			\frac{1}{2} \frac{8 \pi^{2}}{L^{3}} \sum_{k=1}^{\infty}\lal^{2}k^{2} C \intts \intsp \e\us^{n-p-2}
			\dx \ds \, .
		\end{align}
		The first term in the last line of \eqref{EnEnt4} cancels out against the corresponding term in $T_{3}$, the second one will be a Gronwall term with respect to the potential, as for $n-p-2 \in (-p,0)$, the estimate
		\begin{align}
			\e \intts \intsp \us^{n-p-2}
			\dx \ds 
			&\le  C_{1} \intts \intsp 	\e\, \us^{-p} 
			\dx \ds  + C_{2}(T,L) 
		\end{align}
		suggests.
		
		\noindent
		Ad $T_{6}$: Using the same identity as in $T_{2}$ combined with integration by parts, as well as the relations \ref{ONBrelations} in the appendix, we find
		\begin{align}\notag \label{EnEnt6}
			T_{6} &= \Cstr \intts \intsp \us^{n-2} \usxx^{2} \dx \ds
			\\ \notag
			&\quad+ \Cstr \left(\tfrac{(n-2)^{2}}{4}+\tfrac{1}{3} |(n-3)(n-2)|\right) \intts \intsp \us^{n-4} \usx^{4} \dx \ds
			\\ 
			&\quad+
			C \intts \intsp \us^{n-2}\usx^{2} \dx \ds + C \intts \intsp \us^{n}\dx \ds 
			=:
			\operatorname{I} + \dots + \operatorname{IV}\, ,
		\end{align}
		cf. Lemma~\ref{lem:adT6} for the computation.
		The term $\operatorname{I}$ can be absorbed in $T_{2}$, where we do not need the additional part given by $S>0$. 
		For $\operatorname{III}$, we find for arbitrary $\eta>0$
		\begin{align}\label{EnEnt7}
			\operatorname{III} \le \eta  \intts \intsp \us^{n-4} \usx^{4} \dx \ds + C_{\eta}\intts \intsp \us^{n} \dx \ds\, .
		\end{align}
		The last term in \eqref{EnEnt7} and $\operatorname{IV}$ will be estimated via Gagliardo-Nirenberg's inequality, Young's inequality and mass conservation, cf. Lemma~\ref{app:estun}, 
		\begin{align}\label{EnEnt8}
			\intts \intsp \us^{n} \dx \ds \le \eta \intt \intsp \usx^{2} \dx \ds + C\left(M, T, \eta\right)  \,.
		\end{align} 
		Here $M := \intsp \un \dx$. The first term on the right-hand side of \eqref{EnEnt8} will be a Gronwall term. 
		It remains to estimate the term $II$ in \eqref{EnEnt6} and the first term in \eqref{EnEnt7}. 
		We have  
		\begin{align}
			T_{2} + I + II =  \notag
				&-S \intts \intsp \us^{n-2} \usxx^{2} \dx \ds 
				\\
				&+ \left(-S \frac{|(n-2)(n-3)|}{3} 
				+ \Cstr \frac{(n-2)^{2}}{4}\right) \intts \intsp \us^{n-4} \usx^{4} \dx \ds .
		\end{align}
		and the estimate
		\begin{align}\notag \label{est:II}
			\intts &\intsp  \us^{n-4} \usx^{4} \dx \ds 
			= - \frac{3}{n-3} \intts \intsp  \us^{n-3} \usx^{2} \usxx \dx \ds  
			\\ 
			&\le 
			\frac{3 \kappa}{2(3-n)} 
			\intts \intsp  \us^{n-4} \usx^{4} \dx \ds + \frac{3}{2\kappa(3-n)} \intts \intsp  \us^{n-2} \usxx^{2} \dx \ds 
		\end{align}
		for $\kappa > 0$ chosen arbitrarily. Minimizing the prefactor in \eqref{est:II} with respect to $\kappa$ gives $\kappa = \frac{3-n}{3}$. Then
		\begin{align}\notag \label{est:optS1}
			&\left(-S \frac{|(n-2)(n-3)|}{3} + \Cstr \frac{(n-2)^{2}}{4}\right) \intts \intsp \us^{n-4} \usx^{4} \dx \ds
			\\ 
			&\qquad\le
			\left(-S \frac{3(n-2)}{3-n} + \Cstr \frac{9(n-2)^{2}}{4(3-n)^{2}}\right) \intts \intsp \us^{n-2} \usxx^{2} \dx \ds \, . 
		\end{align}
		Hence, to end up with a good term and to keep  a portion of the $\usx^{4}$-term in $T_{2}$ as a good term for future use, we have to choose $S$ such that the condition
		\begin{align*}
			&S + S (1-\mu) \frac{3(n-2)}{3-n} - \Cstr \frac{9(n-2)^{2}}{4(n-3)^{2}} > 0 
			\\
			\Leftrightarrow &S > \Cstr \frac{9(n-2)^{2}}{4(3-n)} \frac{1}{(2n-3) - 3\mu (n-2)} \, ,
		\end{align*}
		where $0 < \mu \ll 1$, is satisfied. For suitable $\mu$ this is met by
		\begin{align}
			S > \Cstr \frac{9}{4} \frac{(n-2)^{2}}{(3-n)(2n-3)} \, 
		\end{align} 
		which is assumption $(A3)^{*}$.
		Hence, we may absorb $II$ in \eqref{EnEnt6} and the first term in \eqref{EnEnt7}, as $\eta$ may be chosen sufficiently small (depending on $\mu$).

		Collecting all terms from the energy, rearranging, combining the constants, applying $q$-th powers, $q\ge1$, 
		and taking suprema as well as expectation,
		we get for arbitrary $t' \in [0,T]$
		\begin{align}\notag \label{EnEntr9}
			&\Ew{\sup_{t \in [0,t' \wedge T_{\sigma}]} E^{\e}[\us(t)]^{q}} 
			+
			\Ew{\left(\int_{0}^{t'\wedge T_{\sigma}} \intsp \us^{n} \psx^{2}\dx \ds \right)^{q}}
			\\ \notag 
			&\quad+	
			(\mu S \frac{|(n-2)(n-3)|}{3}- \eta)^{q}
			\Ew{ \left(\int_{0}^{t'\wedge T_{\sigma}} \intsp \us^{n-4}\usx^{4} \dx\ds \right)^{q}}
			\\ \notag 
			&\quad+
			\left(S + S (1-\mu) \frac{3(n-2)}{3-n} - \Cstr \frac{9(n-2)^{2}}{4(n-3)^{2}}\right)^{q}
			\Ew{ \left(\int_{0}^{t'\wedge T_{\sigma}} \intsp \us^{n-2} \usxx^{2} \dx\ds\right)^{q} }
			\\ \notag
			&\quad+
			S^{q}\Ew{ \left(
				\int_{0}^{t'\wedge T_{\sigma}} \intsp \us^{n-2} \usx^{2} \e \W''(\us) \dx \ds \right)^{q}}
			\\ \notag
			&\le
			C(\un, \delta, q,T,\eta) + C\Ew{ \int_{0}^{t'\wedge T_{\sigma}} E^{\e}(\us)^{q}  \ds }
			\\ \notag
			&\quad+
			\mathbb{E} \Bigg[ \sup_{t \in [0,t' \wedge T_{\sigma}]} \left|
			\sum_{k \in \Z} \lal
			\int_{0}^{t\wedge T_{\sigma}} \intsp F'(\us) \left(\us^{\frac{n}{2}} \gl\right)_{x} \dx \dble \right.
			\\
			&\qquad \qquad
			+
			\left.
			\sum_{k \in \Z} \lal
			\int_{0}^{t\wedge T_{\sigma}} \intsp \usx \left(\us^{\frac{n}{2}} \gl\right)_{xx} \dx \dble
			\right|^{q} \Bigg] \, .
		\end{align}
		
		For an application of the Burkholder-Davis-Gundy inequality, we define the operator
		\begin{align}\notag
			\tau_{E}(s)[v] := &\intort \usx(s) \Big(\us^{\frac{n}{2}}(s) \sum_{\ell \in \Z} \skp{g_{\ell}}{v}_{\lpr{2}{\cO}} g_{\ell} \Big)_{xx} \dx
			\\ \notag
			&\quad+ \e  \intort \W'(\us)(s) (\us^{\frac{n}{2}}(s) \sum_{\ell\in \Z} \left(v, g_{\ell}\right)_{\lpr{2}{\cO}}g_{\ell} )_{x} \dx
			\\ 
			= 	
			&-\intort  \us^{\frac{n}{2}}(s) \psx(s)  \sum_{\ell \in \Z} \skp{g_{\ell}}{v}_{\lpr{2}{\cO}}g_{\ell} \dx \, .
		\end{align}
		Then 
		\begin{align}
			||\tau_{E}(s)||^{2}_{L^{0}_{2}(\lpr{2}{\cO},\R)}
			&=
			\sum_{k\in \Z} \lal^{2} \left|
			-\intort  \us^{\frac{n}{2}} \psx \gl \dx
			\right|^{2}
			&\le
			C((\lal)_{k}, (\gl)_{k})
			\intort 
			\us^{n} \psx^{2}
			\dx \, ,
		\end{align}
		where we used $||\gl||^{2}_{\lpr{\infty}{\cO}} < C$ for all $k \in \Z$ and $\sum_{k \in \Z} \lal^{2} < \infty$ in the last step, cf. (A2).
		Consequently, for positive constants $\eta$ and $C_{\eta}$ arising from Young's inequality
		\begin{align}
			\Ew{\sup_{t \in [0,t' \wedge T_{\sigma}]} \left|T_{5}+T_{7} \right|^{q}}
			\le
			\eta  \Ew{\left( \int_{0}^{t'\wedge T_{\sigma}} \intort 
				\us^{n} \psx^{2}  
				\dx\right)^{q}}
			+ C_{\eta} \, , 
		\end{align}
		where we may absorb the first term.
		Finally, to obtain an estimate of the second term on the right-hand side of \eqref{EnEntr9}, we use
		\begin{align}\notag 
			\mathbb{E}\left[\sup_{t \in [0, t']} E^{\e}[\us(t\wedge T_{\sigma} )]^{q}\right]
			&\le C(\un, \delta, q,T,\eta) + C \, \mathbb{E}\left[\int_{0}^{t'\wedge T_{\sigma}} E^{\e}[\us(s)]^{q} \ds\right]
			\\ 
			&\le C(\un, \delta,q,T,\eta) + C \, \mathbb{E}\left[\int_{0}^{t'} \sup_{t\in [0,s]}E^{\e}[\us(t\wedge T_{\sigma})]^{q}\ds\right] \, .
		\end{align}
		Fubini's theorem and Gronwall's lemma applied with respect to the mapping \\ $t' \mapsto \mathbb{E} \left[ \sup_{t\in [0,t']} E^{\e}[\us(t\wedge T_{\sigma})]^{q} \right]$ then give the desired result,  since $t'\in [0,T]$ was arbitrary.
	\proofend
	
	\subsection{Compactness}\label{sec:Compactness}
	We set $\Je := (\ue)^{n} \pex$.
	To get suitable compactness results, we will use Jakubowski's theorem, cf. \cite{Jakubowski1994}. 
	For this, we define for $\ti{\gamma} < \gamma$ the path spaces
	\begin{align*}
		\cX_{u} &:=  C^{\tilde{\gamma}, \frac{\tilde{\gamma}}{4}}(\Otb) ,
		\quad
		\cX_{\ux} :=  L^{2}(0,T;L^{2}(\cO))_{weak}  ,
		\quad
		\cX_{\uxx} :=  L^{2}(0,T;L^{2}(\cO))_{weak}  ,
		\\
		\cX_{J} &:=  L^{2}(0,T;L^{2}(\cO))_{weak}   ,
		\quad
		\cX_{W} :=  C([0,T];L^{2}(\cO)) ,
		\quad
		\cX_{\un^{\delta}} := H^{1}_{\text{per}}(\cO)
	\end{align*} 
	and set $\cX := \cX_{u} \times \cX_{\ux} \times \cX_{\uxx}
	\times \cX_{J} \times \cX_{W} \times \cX_{\un^{\delta}}$. By standard techniques, cf. for example \cite{FischerGruen} Lemma~5.2, we get
	\begin{lem1}\label{sat:Strafheit}
		The laws  $\mu_{\ue},\mu_{\uex},\mu_{\uexx},
		\mu_{J^{\e}},\mu_{W^{\e}}, \mu_{\un^{\delta}}$ of the corresponding random variables are tight on the path spaces $\cX_{u},\cX_{\ux}, \cX_{\uxx}, 
		\cX_{J}, \X_{W}$, and $\cX_{\un^{\delta}}$, respectively.
	\end{lem1}
	Consequently, with \cite{Jakubowski1994} we infer 
	\begin{prop1}\label{sat:JakubowskiAnwendung}
		For subsequences of $\ue, \uex, \uexx, J^{\e} , \un^{\delta}$, and $W^{\e}$, there exist a probability space $(\ti{\Omega}, \ti{\mathcal{F}}, \ti{\mathbb{P}})$, a sequence of $L^{2}(\cO)$-valued stochastic
		processes $\ti{W}^{\e}$ on $(\ti{\Omega}, \ti{\mathcal{F}}, \ti{\mathbb{P}})$, sequences of random variables
		$
		\uet\colon \ti{\Omega} \rightarrow C^{\ti{\gamma},\frac{\ti{\gamma}}{4}}(\Otb )$, 
		$\wet\colon \ti{\Omega} \rightarrow L^{2}(0,T;L^{2}(\cO))$, 
		$\vet\colon \ti{\Omega} \rightarrow L^{2}(0,T;L^{2}(\cO))$, 
		$\tilde{J}^{\e}\colon \ti{\Omega} \rightarrow \lpr{2}{0,T;\lpr{2}{\cO}}$, 
		$(\unt^{\delta})^{\e}\colon \ti{\Omega} \rightarrow H^{1}_{\text{per}}(\cO)$, 
		random variables
		$
		\ut \in L^{2}(\ti{\Omega};C^{\ti{\gamma},\frac{\ti{\gamma}}{4}}(\Otb ))$, 
		$\wt \in  L^{2}(\ti{\Omega};L^{2}(0,T;L^{2}(\cO)))$, 
		$\vt \in  L^{2}(\ti{\Omega};L^{2}(0,T;L^{2}(\cO)))$,
		\\
		$\tilde{J} \in L^{2}(\ti{\Omega};\lpr{2}{0,T;\lpr{2}{\cO}})$, 
		$\unt^{\delta} \in L^{2}(\ti{\Omega};H^{1}_{\text{per}}(\cO))$,
		as well as an $L^{2}(\cO)$-valued process $\ti{W}$ such that
		\begin{enumerate}
			\item
			for all $\e \in  (0,1)$  the law of 
			$(\uet,\wet, \vet, 
			\Jet,\ti{W}^{\e}, (\unt^{\delta})^{\e})$
			on $\cX$ with respect to  the measure $\ti{\mathbb{P}}$ equals the law of 
			$(\ue, \uex, \uexx,
			J^{\e},W^{\e}, \un^{\delta})$
			with respect to  $\mathbb{P}^{\e}$.
			\item 
			as $\e \rightarrow 0$, the sequence 
			$(\uet, \wet, \vet, 
			\Jet,\ti{W}^{\e}, (\unt^{\delta})^{\e})$ 
			converges
			$\ti{\mathbb{P}}$-almost surely to \newline
			$(\ut, \wt, \vt,
			\tilde{J},\ti{W}, \unt^{\delta})$
			in the topology of $\cX$.		
		\end{enumerate}
	\end{prop1}
	
	\begin{bem1}
		The uniform estimates for members of $(\ue)_{\e \in (0,1)}$ established so far hold true for the corresponding expressions in $(\uet)_{\e\in(0,1)}$ by continuity and the equality of the laws of $\ue$ and $\uet$ for all $\e \in (0,1)$. 
	\end{bem1}
	Using the coincidence of the laws in Proposition~\ref{sat:JakubowskiAnwendung}, we have the following identifications. 
	\begin{lem1}\label{lem:IdentifikationFormel}
		We have 
		$
		\wet = \uext 
		$,
		$
		\vet = \uexxt 
		$,
		as well as
		$
		\Jet = (\uet)^{\frac{n}{2}} (-\uexxt + \e \W'(\uet))_{x}
		$
		$\mathbb{\ti{P}}$-almost surely.
	\end{lem1}
	
	We consider normal filtrations $(\ti{F}_{t})_{t \ge 0 }$ and $(\ti{F}^{\e}_{t})_{t \ge 0}$ with
	\begin{align}\label{def:Ft}
		\ti{F}_{t} := \sigma\left(\sigma(r_{t}\ut,r_{t}\Wt) \cup \{N \in \ti{F}:\ti{\mathbb{P}}(N)=0
		\} \cup \sigma(\unt^{\delta})\right)
	\end{align}  
	and
	\begin{align}\label{def:F}
		\ti{F}^{\e}_{t} := \sigma\left(\sigma(r_{t}\uet,r_{t}\Wet) \cup \{N \in \ti{F}:\ti{\mathbb{P}}(N)=0
		\} \cup \sigma((\unt^{\delta})^{\e}) \right) \, .
	\end{align}  
	Here, $r_{t}$ is the restriction of a mapping on $[0,T]$ to the time interval $[0,t]$, $t \in [0,T]$.
	The proof of the next lemma can be found in \cite{FischerGruen} Lemma~5.7.
	
	\begin{lem1}\label{lem:WienerProcesses}
		The processes $\Wet$ and $\Wt$ are Q-Wiener processes which are adapted to the filtrations $(\ti{F}^{\e}_{t})_{t \ge 0}$ and $(\ti{F}_{t})_{t \ge 0 }$, respectively. We have
		\begin{align}\label{lem:WienerProcesses-h1}
			\Wet(t) = \sum_{k \in \Z} \lal \blet(t)\gl
		, \quad 
			\Wt(t) = \sum_{k \in \Z} \lal \blt(t)\gl
		\end{align}
		with families $(\blet)_{k \in \Z}$ and $(\blt)_{k \in \Z}$ of i.i.d. standard Brownian motions with respect to  $(\ti{F}^{\e}_{t})_{t\ge 0}$ and $(\ti{F}_{t})_{t\ge 0}$, respectively.
	\end{lem1}
	Now, we will identify the limits in Proposition~\ref{sat:JakubowskiAnwendung}. 
	In particular, we show that the corresponding terms in the equation \eqref{eq:stfe-Strat-equ} converge to the terms stipulated by \eqref{eq:stfe-Strat-equ1}. 
	First, we consider terms that do not involve a stochastic integral. 
	We start with the limits $\wt$, $\vt$ and $\Jt$ and get easily the following identities.
	\begin{lem1}\label{lem:Identwt}
		We have $\mathbb{\ti{P}}$-almost surely
		$
		\wt =  \uxt 
		$
		and 
		$
		\vt =  \uxxt.
		$	
	\end{lem1}

	\begin{bem1}\label{bem:ut-nullset}
		Similarly as before, we may apply Jakubowski's theorem to infer the weak convergence of $(\uet)^{\frac{2-n}{2}}$ in $L^{2}(\Ot)$ $\tilde{\mathbb{P}}$-almost surely and thus also the uniform boundedness of $||(\uet)^{\frac{2-n}{2}}||_{L^{2}(\Ot)}$ $\tilde{\mathbb{P}}$-almost surely. Using additionally the uniform convergence $\uet \to \ut$ $\tilde{\mathbb{P}}$-almost surely for $\e \to 0$, it follows easily that $|\{\ut = 0\}| = 0$ on $\Ot$ $\tilde{\mathbb{P}}$-almost surely.
	\end{bem1}
		
	Regarding $\Jt$, we have the following result. The proof is in the spirit of Proposition~5.7 of \cite{DGGG}.
	\begin{lem1}\label{lem:weakconvdiss}
		The weak derivative $\ut_{xxx}$ satisfies $\ut_{xxx} \in \lpr{2}{\{\ut > r\}}$ for any $r>0$, $\tilde{\mathbb{P}}$-almost surely.
		Moreover, we have $\Jet \chi_{\{\ut > 0\}} = (\uet)^{\frac{n}{2}} \pet_{x} \chi_{\{\ut > 0\}}\rightharpoonup - \ut^{\frac{n}{2}} \ut_{xxx} \chi_{\{\ut > 0\}}$ weakly in $\lpr{2}{\Ot}$ $\tilde{\mathbb{P}}$-almost surely, and
		\begin{align}\label{weakconvdiss1}
			\intT \intort (\uet)^{n} \pet_{x} \phi \dx\ds \rightarrow -\intT  \intort (\ut)^{n} \ut_{xxx} \phi \dx\ds \,
		\end{align}
		for all $\phi \in \lpr{2}{\Ot}$, $\tilde{\mathbb{P}}$-almost surely.  
	\end{lem1}
	\proof
	Let us fix an arbitrary $r>0$.
	We will first derive an $\e$-independent bound of $\uet_{xxx}$ in $\lpr{2}{\{\ut>r\}}$, $\tilde{\mathbb{P}}$-almost surely. To do so, we write 
	\begin{align}\notag
		&\intT \int_{\{\ut(s, \cdot) >r \}} (\uet)^{n} (\pet_{x})^{2} \dx \ds 
		\\ 
		&=
		\intT \int_{\{\ut(s, \cdot) >r \}} (\uet)^{n} (\uet_{xxx})^{2} - 2 \e (\uet)^{n-p-2} \uet_{xxx} \uet_{x} + \e^{2} (\uet)^{n-2p-4} (\uet)^{2}_{x} \dx \ds 
	\end{align}
	and obtain
	\begin{align}\notag
		\intT \int_{\{\ut(s, \cdot) >r \}} (\uet)^{n} (\uet_{xxx})^{2} \dx \ds 
		\le
		&\intT \intort (\uet)^{n} (\pet_{x})^{2} \dx \ds 
		\\
		&+ 2  \e \intT \int_{\{\ut(s, \cdot) >r \}} (\uet)^{n-p-2} \uet_{xxx} \uet_{x} \dx \ds  \,.
	\end{align} 
	By Young's inequality and Proposition~\ref{sat:JakubowskiAnwendung}, we have $\tilde{\mathbb{P}}$-almost surely
	\begin{align}\notag \label{unifbounddiss}
		&\frac{1}{2} \intT \int_{\{\ut(s, \cdot) >r \}} (\uet)^{n} (\uet_{xxx})^{2} \dx \ds 
		\\ \notag
		&\le 
		\intT \intort  (\uet)^{n} (\pet_{x})^{2} \dx \ds + C \e^{2} \intT \int_{\{\ut(s, \cdot) >r \}} (\uet)^{n-2p-4} (\uet_{x})^{2} \dx \ds
		\\ 
		&\le 
		C(\omega) + C \left(\e \sup_{\{\ut >r \}} (\uet)^{-p} \right) \e \intT \int_{\{\ut(s, \cdot) >r \}} (\uet)^{n-p-4} (\uet_{x})^{2} \dx \ds
	\end{align}
	Using the uniform convergence of $\uet \rightarrow \ut$, cf. Proposition~\ref{sat:JakubowskiAnwendung}, we may choose $\e = \e(\omega)$ such that 
	\begin{align}\label{epssmall}
		||\ut - \uet||_{L^{\infty}(\Ot)} < \frac{r}{2}\, ,
	\end{align}
	$\tilde{\mathbb{P}}$-almost surely.
	Under this assumption the integral term in the last line of \eqref{unifbounddiss} is controlled $\tilde{
		\mathbb{
			P}}$-almost surely via the uniform bound on $\uext$ implied by Proposition~\ref{sat:JakubowskiAnwendung} for $\e$ sufficiently small depending on $\omega$:
	\begin{align}
		\intT \int_{\{\ut(s, \cdot) >r \}} (\uet)^{n-p-4} (\uet_{x})^{2} \dx \ds \le 
		\left(\frac{r}{2}\right)^{n-p-4} \sup_{t\in [0,T]}\intort (\uet_{x}(t))^{2} \dx\, \le C(r, \omega) .
	\end{align}
	Therefore, from \eqref{unifbounddiss}, we deduce the following bound independent of $r$
	\begin{align}\label{limsupdisbound}
		\limsup_{\e \rightarrow 0} \intT \int_{\{\ut(s, \cdot) >r \}} (\uet)^{n} (\uet_{xxx})^{2} \dx \ds \le C < \infty \,  
	\end{align}
	$\tilde{\mathbb{P}}$-almost surely with a constant $C=C(\omega)$.  
	
	Combining \eqref{epssmall} with \eqref{unifbounddiss} and \eqref{limsupdisbound}, we find 
	\begin{align}
		\intT \int_{\{\ut(s, \cdot) >r \}} (\uet_{xxx})^{2} \dx \ds  \le \frac{2^{n}}{r^{n}} \intT \int_{\{\uet(s, \cdot) >\frac{r}{2} \}} (\uet)^{n} (\uet_{xxx})^{2} \dx \ds \le C(r, \omega)
	\end{align}
	$\tilde{\mathbb{P}}$-almost surely with a constant $C(r, \omega)$ independent of $\e$. This implies for a subsequence
	\begin{align}\label{weakconvuxxx}
		\uet_{xxx} \chi_{\{\ut >r \}} \rightharpoonup \eta^{r} \chi_{\{\ut >r \}}
	\end{align}
	in $\lpr{2}{\Ot}$ $\tilde{\mathbb{P}}$-almost surely, where $\eta^{r} \in \lpr{2}{\Ot}$ depends on $\omega$.
	For the identification $\eta^{r} = \ut_{xxx}$ on $\{\ut >r \}$ $\tilde{\mathbb{P}}$-almost surely, we refer to \cite{DGGG} Proposition~5.7.
	
	Let us now identify $\Jt$ and show \eqref{weakconvdiss1}.
	By Proposition~\ref{sat:JakubowskiAnwendung} and Lemma~\ref{lem:Identwt} we have
	\begin{align}\label{weakconvJe}
		\Jet = (\uet)^{n/2}\pet_{x} = (\uet)^{n/2}\left(-\uet_{xx} + \e \W'(\uet)\right)_{x} \rightharpoonup \Jt \, 
	\end{align} 
	weakly in $\lpr{2}{\Ot}$, $\tilde{\mathbb{P}}$-almost surely.
	We aim to identify $\Jt \in \lpr{2}{\Ot}$ on the sets
	$\{\ut > r \}$
	where $r >0$ is as before arbitrary but fixed.  
	We have
	\begin{align}\notag
		&\left|\intT \intort \Jet \chi_{\{\ut >r \}} \tilde{\phi} \dx \ds - \intT \intort - \ut^{n/2} \ut_{xxx} \chi_{\{\ut >r \}} \tilde{\phi} \dx \ds \right| 
		\\ \notag
		&\le 
		\left|\intT \int_{\{\ut(s, \cdot) >r \}} \left( - (\uet)^{n/2} + \ut^{n/2}\right) \uet_{xxx} \tilde{\phi} \dx \ds\right|
		\\ \notag
		&\quad+\left| 
		\intT \int_{\{\ut(s, \cdot) >r \}} \left(- \uet_{xxx} + \ut_{xxx}\right) (\ut)^{n/2} \tilde{\phi} \dx \ds\right|
		\\ 
		&\quad+\left|
		\intT \int_{\{\ut(s, \cdot) >r \}}(\uet)^{n/2}  \e (\uet)^{-p-2} \ue_{x} \tilde{\phi} \dx \ds\right|  
		=: I_{1} + I_{2} + I_{3}\, , 
	\end{align}    
	where $\ti{\phi} \in L^{\infty}(\Ot)$.
	The uniform convergence $\uet \rightarrow \ut$ $\tilde{\mathbb{P}}$-almost surely as well as \eqref{weakconvuxxx} show $\lim_{\e \rightarrow 0} I_{1} = \lim_{\e \rightarrow 0} I_{2} = 0$. For $I_{3}$ we have, again assuming $\e$ to be small enough, 
	\begin{align}\notag
		I_{3} 
		&\le 
		\left(\e \intT \int_{\{\uet(s,\cdot) >\frac{r}{2} \}} (\uet)^{-p} \dx \ds\right)^{\frac{1}{2}} 
		\left(\e \intT \int_{\{\uet(s,\cdot) >\frac{r}{2} \}} (\uet)^{n-p-4} (\uet)^{2}_{x}  \dx \ds\right)^{\frac{1}{2}} ||\ti{\phi}||_{\lpr{\infty}{\Ot}} 
		\\ 
		&\le 
		\left(\e \intT \int_{\{\uet(s,\cdot) >\frac{r}{2} \}} (\uet)^{-p} \dx \ds\right)^{\frac{1}{2}} 
		\left(C(r) \e\sup_{t\in [0,T]} \intort (\uet)^{2}_{x}  \dx \ds\right)^{\frac{1}{2}} ||\ti{\phi}||_{\lpr{\infty}{\Ot}}
		\rightarrow 0
	\end{align}
	as the second term tends to zero while the first one is uniformly bounded. Hence, we have identified $\Jt \chi_{\{\ut >r\}} = - \ut^{n/2}\ut_{xxx} \chi_{\{\ut >r\}}$  $\tilde{\mathbb{P}}$-almost surely for $r>0$, and since $r$ was arbitrary, we find the identification of $\Jt$ on $\{\ut > 0 \}$ as well. 
	As $(\uet)^{\frac{n}{2}}$ converges $\tilde{\mathbb{P}}$-almost surely uniformly on $\Ot$, we finally obtain $(\uet)^{\frac{n}{2}} \Jet \rightharpoonup \ut^{\frac{n}{2}} \Jt$ in $\lpr{2}{\{\ut >0 \}}$ $\tilde{\mathbb{P}}$-almost surely.
	This implies \eqref{weakconvdiss1}, as the set ${\{\ut = 0 \}} \subset \Ot$ has measure zero $\tilde{\mathbb{P}}$-almost surely due to the bound of $\mathcal{G}(\uet)$ in Lemma~\ref{lem:EnEnt} and $\uet \rightarrow \ut$ uniformly on $\Ot$ $\tilde{\mathbb{P}}$-almost surely, cf. Remark~\ref{bem:ut-nullset}.
	\proofend
	Using the convergence properties of $\uet \to \ut$ and $\uext \rightharpoonup \uxt$, we establish easily:
	\begin{lem1}\label{lem:weakconvcorr}
		We have 
		\begin{align}
			\intT \intsp (\uet)^{n-2} \uext \phi \dx\ds\rightarrow \intT \intsp (\ut)^{n-2} \uxt \phi\dx\ds
		\end{align}
		$\tilde{\mathbb{P}}$-almost surely for any $\phi \in L^{2}(\Ot)$.
	\end{lem1}
	 Finally, we show that the term involving the effective interface potential vanishes for $\e \to 0$.
	\begin{lem1}
		We have for $\phi \in W^{2,\infty}_{per}(\ort)$
		\begin{align}
			\Ew{\left|\intT \intort (\uet)^{n} \e (\W'(\ue))_{x} \phi_{x} \dx \ds\right|} \rightarrow 0
		\end{align}
		$\ti{\mathbb{P}}$-almost surely for $\e \rightarrow  0 $.
	\end{lem1}
	\proof
	Using Young's inequality, we find for $\eta > 0$ and $p \in (2,\infty)$ with $\W(x) = x^{-p}$
	\begin{align}\notag
		&\abs{\intsptiT{\e(\uet)^{n}\uext (\uet)^{-p-2} \phi_{x}  }}
		=\abs{ \frac{1}{n-p-1} \intsptiT{\e(\uet)^{\frac{-p}{2}} (\uet)^{\frac{2n-p-2}{2}} \phi_{xx}}} 
		\\ \notag
		&\le \frac{C}{4}\intsptiT{\e^{2} \e^{\eta -1} (\uet)^{-p} \phi^{2}_{xx} } 
		 +C  \intsptiT{\e^{1-\eta}	(\uet)^{2n-p-2}} 
		=: \operatorname{I} + \operatorname{II}.
	\end{align}
	Due to the boundedness of $\phi_{xx}$ and \eqref{est:EnEnt}, we have for $ \operatorname{I}$
	\begin{align}\notag
		C \e^{\eta +1} \Ew{\intsptiT{\W (\uet) \phi_{xx}^{2}}}
		\le
		\e^{\eta +1} C \,  \Ew{\sup_{t\in [0,T]}\intsp \W(\uet(t)) \dx} 
		\le 
		\e^{\eta} C \rightarrow 0 
	\end{align}
	$\mathbb{\ti{P}}$-almost surely for $\e \to 0$.
	For $\operatorname{II}$ we argue with $\delta > 0$ as follows: 
	\begin{align}\notag
		& \intsptiT{\e^{1-\eta} (\uet)^{2n-p-2}	} 
		\\ \notag
		&=\intT \int_{[\uet<\e^{\frac{\eta+\delta}{2n-2}}]}{\e^{1-\eta} (\uet)^{2n-p-2}} \dx \ds
		+ \intT \int_{[\uet \ge \e^{\frac{\eta+\delta}{2n-2}}]} \e^{1-\eta} (\uet)^{2n-p-2} \dx \ds
		\\ \notag
		&\le \e^{1-\eta+\eta + \delta} \intsptiT{\W(\uet)} 
		+ \intT \int_{[\uet \ge \e^{\frac{\eta+\delta}{2n-2} }]}{\e^{1-\eta+(2n-p-2)\left(\frac{\eta+\delta}{2n-2}\right)}	}\dx \ds 
		\\ \notag
		&= 
		\e^{1+\delta} \intsptiT{\W(\uet)} 
		+ C \e^{1+\delta -\frac{(\eta+\delta)p}{2n-2}} \,.
	\end{align}
	Choosing for example $\eta = \delta = \frac{2n-2}{2p}>0$, we find $1+\delta -\frac{(\eta+\delta)p}{2n-2}=\delta$. Then,
	$
	\lim_{\e \to 0}\Ew{\operatorname{II}}
	=0
	$
	follows.
	\proofend
	\subsection{Identification of the stochastic integral}
	\label{sec:IdentificationStochInt}
	To show that under the already established convergence properties also the stochastic integral converges to the correct expression, we will use the fact that under regularity assumptions met in our setting the stochastic integral in \eqref{eq:stfe-Strat-equ} -- and hence also the corresponding version on the new probability space $(\tilde{\Omega}, \tilde{\mathcal{A}}, \tilde{\mathbb{P}})$ given by Jakubowski's theorem -- are martingales, see Lemma~\ref{sat:martingaleigenschaftTilde}. 
	Proving that also the limit expression remains to be a martingale, cf. Lemma~\ref{lem:M0Martingal} and Corollary~\ref{kor:M0Martingal}, for the identification we apply the ideas from \cite{BreitHofmanova2016,Brezniak2007,DebusscheHofmanovaVovelle,HofmanovaSeidler2012}. 
	For  $\phi \in H^{1}_{\text{per}}(\cO)$ arbitrary but fixed, we consider the processes 
	$\M_{\e,\phi}: \Omega \times [0,T] \rightarrow \mathbb{R} $ defined by
	\begin{align}\label{def:Meps}
		\M_{\e,\phi}(t) := 
		&\intsp{(\ue(t)-\un^{\delta})\phi} \dx \notag
		 - \intt \int_{\{\ut(s,\cdot) >0\}}(\ue)^{n} \pex \phi_{x} \dx \ds \notag	 
		\\
		&- (\Cstr+S) \intsptit{ (\ue)^{n-2} \uex \phi_{x} } \, . 
	\end{align}
	We then have  
	\begin{align}\label{DarstellungMeItoIntegral}
		\M_{\e,\phi}(t) 	&=  \sum_{k \in \Z}  \lal \int_{0}^{t} \int_{\cO} (\gl(\ue)^{\frac{n}{2}})_{x}\phi \dx \dble 
	\end{align}
	for $t \in [0,T]$, i.e. $\M_{\e,\phi}$ is a continuous, square integrable $\mathcal{F}_{t}^{\e}$-martingale.
	We will need the following results which can be proven similarly as Lemma~5.10 and Lemma~5.12 in \cite{GruenKlein24Exist}: 
	\begin{align}\label{eq:QuadrVarMeps}
		\qv{\M_{\e,\phi}} 
		= \int_{0}^{\cdot}{ \sum_{k \in \Z} \lal^{2} 
			\left(
			\intsp ((\ue)^{n/2} \gl)_{x}\phi \dx
			\right)^{2}	} \ds ,
	\end{align}
	\begin{align}\label{est:QuadrVarMeps}
		\qv{\M_{\e,\phi}} \le C \Norm{\phi}^{2}_{H^{1}_{\text{per}}} \inttimeT{\Norm{(\ue)^{n/2}}^{2}_{L^{2}(\cO)}},
	\end{align}
	and for $k \in\N$ 
	\begin{align}\label{QuadrCovarMbeta}
		\qv{\M_{\e,\phi}, \ble} = 
		\lal \int^{\cdot}_{0}\int_{\cO}{ ((\ue)^{n/2} \gl )_{x} \phi} \dx \ds .
	\end{align}
	With these results at hand we can establish 
	\begin{kor1}\label{kor:Martingal1}
		Let $k \in\N$ and $\e \in (0,1)$. The processes
		\begin{align}
			\M_{\e,\phi}^{2} - \int_{0}^{\cdot}{ \sum_{k \in \Z} \lal^{2} 
				\left(
				\intsp ((\ue)^{n/2} \gl)_{x}\phi \dx
				\right)^{2}	} \ds 
		\end{align}
		and 
		\begin{align}
			\M_{\e,\phi}\ble -	\lal \int_{0}^{\cdot}\intsp ((\ue)^{n/2} \gl )_{x} \phi \dx \ds
		\end{align}	
		are continuous $\mathcal{F}^{\e}_{t}$-martingales. 
	\end{kor1}
	By the equality of laws stated in Proposition~\ref{sat:JakubowskiAnwendung}, we get the same results for  
	\begin{align*}
		\ti{\M}_{\e,\phi}(t) := 
		&\intsp (\uet(t)-(\ut_{0}^{\delta})^{\e})\phi \dx \notag
		- \intt \int_{\{\ut(s,\cdot) >0 \}} (\uet)^{n} \pex \phi_{x} \dx \ds \notag	 
		\\
		&- (\Cstr+S) \intsptit{ (\uet)^{n-2} \uext \phi_{x} }
		\, .
	\end{align*}

	\begin{lem1}\label{sat:martingaleigenschaftTilde} 
		For $k \in\N$ and $\e \in (0,1)$,
		\begin{align}\label{sat:martingaleigenschaftTildehh1}
			&\ti{\M}_{\e,\phi}\, , 
			\\ 
			\label{sat:martingaleigenschaftTildehh2}
			&\ti{\M}_{\e,\phi}^{2} - \int_{0}^{\cdot}{ \sum_{k \in \Z} \lal^{2} 
				\left(
				\intsp ((\uet)^{n/2} \gl)_{x}\phi \dx
				\right)^{2}	} \ds  \, ,
			\\ \label{sat:martingaleigenschaftTildehh3}
			&\ti{\M}_{\e,\phi}\ti{\ble} -	\lal \int_{0}^{\cdot}\int_{\cO} ((\uet)^{n/2} \gl )_{x} \phi \dx \ds 
		\end{align}
		are continuous $\ti{\mathcal{F}}_{t}^{\e}$-martingales. Moreover, on $[0,T]$ we have
		\begin{align}\label{sat:martingaleigenschaftTildeh1}
			&\qv{\ti{\M}_{\e,\phi}}_{t} =  \int_{0}^{t}{ \sum_{k \in \Z} \lal^{2} 
				\left(
				\intsp ((\uet)^{n/2} \gl)_{x}\phi \dx
				\right)^{2}	} 
			\ds  \, , 
			\\\label{sat:martingaleigenschaftTildeh2}
			&\qv{\ti{\M}_{\e,\phi},\ti{\ble}}_{t} = 	\lal \intsptit{ ((\uet)^{n/2} \gl )_{x} \phi}.
		\end{align}
	\end{lem1}
	The next step is to show that the martingale property is preserved in the limit. We show that for $\phi \in H^{1}_{\text{per}}(\cO)$
	\begin{align}\label{def:M0t}
		\ti{\M}_{0,\phi}(t) :=
		&\intsp (\ut(t)-\unt^{\delta})\phi \dx \notag
		- \intt \int_{\{\ut(s,\cdot) >0 \}} \ut^{n} \ut_{xxx} \phi_{x} \dx \ds \notag	 
		\\
		& - (\Cstr+S) \intsptit{ \ut^{n-2} \uxt \phi_{x} }
	\end{align}
	has the martingale property.
	\begin{lem1}\label{lem:M0Martingal}
		For $s,t \in [0,T]$ with $s \le t$ and
		for all continuous functions \newline $\Psi \colon C^{\ti{\gamma} ,\frac{\ti{\gamma}}{4}}(\bar{\cO}\times[0,s]) \times C([0,s];L^{2}(\cO)) \rightarrow [0,1]$, we have 
		\begin{align}\label{sat:M0Martingalh1}
			\Ew{\Psi(r_{s}\ut, r_{s} \ti{W}) \left(\ti{\M}_{0,\phi}(t)-\ti{\M}_{0,\phi}(s)\right)} = 0. 
		\end{align}
	\end{lem1}  
	\proof 
	We treat the terms inside the expectation in (\ref{sat:M0Martingalh1}) one by one. The continuity of $\Psi$ as well as the convergence of $\uet$ to $\ut$ and of $\ti{W}^{\e}$ to $\ti{W}$ in $C(0,T;L^{2}(\cO))$ show
	\begin{align}\label{KonvergenzPsi}
		\underset{\e \rightarrow 0}{\lim} \Psi(r_{s}\uet, r_{s} \ti{W}^{\e}) =  \Psi(r_{s}\ut, r_{s} \ti{W})
	\end{align}
	$\ti{\mathbb{P}}$-almost surely on $[0,1]$. To see the convergence of the expected values, we use Vitali's theorem. Therefore, since $\Psi$ is bounded, it suffices to show uniform boundedness of moments of the integral-terms in (\ref{sat:M0Martingalh1}) and use the convergence results already established.
	\\
	By the strong convergence of $\uet$ in $ C^{\ti{\gamma} ,\frac{\ti{\gamma}}{4}}(\Otb)$, cf. Proposition~\ref{sat:JakubowskiAnwendung}, we have $\ti{\mathbb{P}}$-almost surely 
	\begin{align}\label{Konvergenzue}
		\underset{\e \rightarrow 0}{\lim}	\intsp{(\uet(t)-\uet(s))\phi} \dx =  \intsp{(\ut(t)-\ut(s))\phi} \dx \, .
	\end{align}    
	The energy-entropy estimate \eqref{est:EnEnt}  and Poincar\'e's lemma give
	\begin{align}\label{est:u_Linfty}
		\Ew{\underset{t \in [0,T]}{\operatorname{esssup}} 	||\uet(t, \cdot)||^{2q}_{\lpr{\infty}{\cO}}} < C
	\end{align}	
	uniformly for an arbitrary $q > 1$ and thus, the uniform integrability of a $q$-th absolute moment of $\intsp{(\uet(t)-\uet(s))\phi} \dx$.
	\\
	Next, we consider $(\uet)^{n} \pet_{x}$. Weak convergence in $\lpr{2}{\Ot}$ $\mathbb{\ti{P}}$-almost surely, and therefore convergence of $\intsptist{(\uet)^{n} (\pet_{x})\phi_{x} }$  $\mathbb{\ti{P}}$-almost surely, has been established in Lemma~\ref{lem:weakconvdiss}.
	Thus it remains to find a uniform bound of a $q$-th moment. Using Young's inequality, $\phi \in H^{1}_{per}(\cO)$, as well as Lemma~\ref{lem:EnEnt}, we find for $q>1$ 
	\begin{align}\notag
		&\Ew{\abs{\intsptist{(\uet)^{n} (\pet_{x})\phi_{x} } }^{q } }
		\\ 
		&\le C(q) \Ew{\left(\frac{1}{2}\intsptist{(\uet)^{n} (\pet_{x})^{2} }\right)^{q} + 
			\left( \frac{1}{2}\intsptist{\phi_{x}^{2} }
			\right)^{q} ||(\uet)||^{nq}_{\lpr{\infty}{\Ot}} } 
		\le C \, ,
	\end{align} 
	where we used \eqref{est:u_Linfty} in the last step.
	Similarly, for the Stratonovich correction term we find 
	\begin{align}\notag
		&\Ew{\abs{\intsptist{(\uet)^{n-2}(\uext)\phi_{x}} }^{q}  }
		\\ 
		&\le 
		C \, \Ew{\abs{\intsptist{ ((\uet)^{n-2} (\uext)^{2}  } }^{q}  }	
		+
		\, C \Ew{\abs{\intsptist{(\uet)^{n-2}\phi_{x}^{2} } }^{q} }
		\le 
		C 
	\end{align}
	which we combine with Lemma~\ref{lem:weakconvcorr}.
	\proofend
	\noindent
	Dynkin's lemma, cf. \cite{Karatzas2014} S.49,
	in combination with Lemma~\ref{lem:M0Martingal} implies the martingale property for $\ti{\M}_{0,\phi}$, cf. for example \cite{HofmanovaSeidler2012}. Hence, we have the following.
	\begin{kor1}\label{kor:M0Martingal}
		$\ti{\M}_{0,\phi}$ is a continuous $\ti{\F_{t}}$-martingale.
	\end{kor1}
	By similar arguments as in Lemmas~5.14 and 5.15 in \cite{FischerGruen}, using in particular \eqref{est:QuadrVarMeps}, and Lemma~\ref{app:estun}, we can show 
	\begin{lem1}For 
		$0 \le s \le t \le T$ and $\Psi$ as in Lemma~\ref{lem:M0Martingal}
		\begin{align}\label{sat:Martingalea1}
			\Ew{ \Psi(r_{s}\ut,r_{s}\ti{W})  \left( \ti{\M}_{0,\phi}^{2}(t)- \ti{\M}_{0,\phi}^{2}(s) 
				- \int_{s}^{t}{ \sum_{k \in \Z} \lal^{2} 
					\left(
					\intsp (\ut^{n/2} \gl)_{x}\phi \dx
					\right)^{2}	} \ds \right)} = 0
		\end{align}
		and 
		\begin{align}\label{sat:Martingalea2}
			\Ew{ \Psi(r_{s}\ut,r_{s}\ti{W}) \left( (\ti{\M}_{0,\phi}\blt)(t) - (\ti{\M}_{0,\phi}\blt)(s) -	\lal \intsptist{ (\ut^{n/2} \gl )_{x} \phi}\right)  } = 0 \, .
		\end{align}
	\end{lem1}
	As in Corollary~\ref{kor:Martingal1}, we find that 
	\begin{kor1}\label{kor:Martingale1}
		\begin{align}
			\ti{\M}_{0,\phi}^{2}
			- \int_{0}^{\cdot} \sum_{k \in \Z} \lal^{2} 
			\left(
			\intort (\ut^{n/2} \gl)_{x}\phi \dx
			\right)^{2}	 \ds 
		\end{align}
		and
		\begin{align}
			(\ti{\M}_{0,\phi}\blt) -	\lal \int_{0}^{\cdot} \intort (\ut^{n/2} \gl )_{x} \phi \dx \ds
		\end{align}
		are continuous $\ti{\mathcal{F}}_{t}$-martingales.
	\end{kor1}
	Consequently, we have the following identities
	\begin{kor1}\label{kor:Martingale}
		For $t \in [0,T]$ 
		\begin{align}
			\qv{\ti{\M}_{0,\phi} }_{t} =
			\int_{0}^{t}{ \sum_{k \in \Z} \lal^{2} 
				\left(
				\intsp{(\ut^{n/2} \gl)_{x}\phi} \dx
				\right)^{2}	} \ds,
		\end{align}
		\begin{align}
			\qv{\ti{\M}_{0,\phi},\ti{\beta}_{k}}_{t} = \lal \intsptit{ (\ut^{n/2} \gl )_{x} \phi} ,
		\end{align}
		\begin{align}
			\qv{\int_{0}^{\cdot} \intsp \sum_{k \in \Z} \lal (\gl \ut^{n/2})_{x}\phi \dx\dif \ti{\beta}_{k} }_{t} = \int_{0}^{t}{ \sum_{k \in \Z} \lal^{2} 
				\left(
				\intsp{(\ut^{n/2} \gl)_{x}\phi} \dx
				\right)^{2}	} \ds
		\end{align}
		holds.
	\end{kor1}
	\proof
	This follows from Corollary~\ref{kor:Martingale1} and the same calculation as for (\ref{eq:QuadrVarMeps}), cf. \cite{FischerGruen} Lemma~5.10 and Lemma~5.12.
	\proofend
	Following the argumentation of Lemma~5.16 in \cite{FischerGruen}, the identification of the stochastic term is achieved. 
	\begin{lem1}\label{sat:IdentM0}
		It holds $\mathbb{\ti{P}}$-almost surely
		\begin{align}\label{sat:IdentM0-h3}
			\ti{\mathcal{M}}_{0,\phi} =  \sum_{k \in \Z}\int_{0}^{\cdot}\intsp{ \lal (\ut^{n/2} \gl)_{x}\phi} \dx \dblt \, .
		\end{align}	
	\end{lem1}
	
	\proof
	We show
	\begin{align}\label{sat:IdentM0-h0}
		\qv{\ti{\M}_{0,\phi}- \int_{0}^{\cdot} \intsp{\sum_{l \in \Z} \lal (\ut^{n/2} \gl)_{x}\phi} \dx \dif\ti{\beta}_{k} } = 0 .
	\end{align}
	From the definition of the quadratic variation we get
	\begin{align}\notag
		\qv{\ti{\M}_{0,\phi} - \int_{0}^{\cdot} \intsp{\sum_{k \in \Z} \lal (\ut^{n/2} \gl)_{x}\phi} \dx \dif\ti{\beta}_{k} }_{t} 
		=
		&\qv{\ti{\M}_{0,\phi} }_{t} 
		+
		\qv{\int_{0}^{\cdot} \intsp \sum_{k \in \Z} \lal (\ut^{n/2} \gl)_{x}\phi  \dx \dif \ti{\beta}_{k} }_{t} \\ \label{sat:IdentM0-h1}
		&- 2 \qv{\ti{\M}_{0,\phi},\int_{0}^{\cdot} \intsp{\sum_{k \in \Z} \lal (\ut^{n/2} \gl)_{x}\phi} \dx \dif\ti{\beta}_{k}}_{t} \,.
	\end{align}
	For the third term in (\ref{sat:IdentM0-h1}) we use the cross-variation formula, cf.  \cite{Karatzas2014} Lemma~2.16 in section 3.2. We get
	\begin{align}\label{sat:IdentM0-h2}
		\qv{\ti{\M}_{0,\phi},\int_{0}^{\cdot} \intsp{\sum_{k \in \Z} \lal (\ut^{n/2} \gl)_{x}\phi} \dx \dif\ti{\beta}_{k}	}_{t}   
		=
		\int_{0}^{t} \intsp{\sum_{k \in \Z} \lal (\ut^{n/2} \gl)_{x}\phi} \dx \dif \, \qv{\ti{\M}_{0,\phi},\ti{\beta}_{k}}
	\end{align}
	for all $t \in [0,T]$.
	Corollary~\ref{kor:Martingale} states
	\begin{align}
		\qv{\ti{\M}_{0,\phi},\ti{\beta}_{k}}_{t} = \lal \intsptit{ ( \ut^{n/2} \gl)_{x} \phi} \, .
	\end{align}
	Since $s \mapsto \qv{\ti{\M}_{0,\phi},\ti{\beta}_{k}}_{s}$ is absolutely continuous, we conclude
	\begin{align}
		\dif \,	\qv{\ti{\M}_{0,\phi},\ti{\beta}_{k}}_{s} = \lal \intsp{(\ut^{n/2}(s) \gl )_{x} \phi} \dx \ds \, .
	\end{align}
	Together with (\ref{sat:IdentM0-h2}) 
	\begin{align}
		\qv{\ti{\M}_{0,\phi},\int_{0}^{\cdot} \intsp{\sum_{k \in \Z} \lal (\ut^{n/2} \gl)_{x}\phi} \dx \dif \ti{\beta}_{k}}_{t} 
		=
		\int_{0}^{t} 	\sum_{k \in \Z} \lal^{2} \left( \intsp (\ut \gl  )_{x} \phi \dx \right)^{2}	 \ds \, 
	\end{align}
	follows.
	Thus, (\ref{sat:IdentM0-h1}) and the identities from Corollary~\ref{kor:Martingale} imply (\ref{sat:IdentM0-h0}). This proves (\ref{sat:IdentM0-h3}).
	\proofend 
	
	\subsection{Proof of Theorem~\ref{theo:MainResult}}\label{sec:proofTheo11}
	Having established Lemma~\ref{sat:IdentM0}, together with the results on convergence in Section~\ref{sec:Compactness} we may give the proof of Theorem~\ref{theo:MainResult}.
	\begin{myproof}{ of Theorem~\ref{theo:MainResult}}
		From Proposition~\ref{sat:JakubowskiAnwendung}, Lemma~\ref{lem:WienerProcesses}, and the identifications in Lemma~\ref{lem:Identwt} and Lemma~\ref{lem:weakconvdiss}, the existence of a stochastic basis, a $Q$-Wiener process $\ti{W}$, as well as random variables $\ut \in \lpr{2}{\Omega;C^{\ti{\gamma},\frac{\ti{\gamma}}{4}}(\bar{\cO}_{T})} \cap  \lpr{2}{\Omega; \lpr{2}{[0,T];H^{1}(\ort)}}$, $\Jt \in \lpr{2}{\Omega; \lpr{2}{\Ot}}$, and $\unt^{\delta}$ with 
		\begin{align}
			\Jt = \ut^{\frac{n}{2}}\ut_{xxx}  \qquad \ti{\mathbb{P}}\text{-almost surely in } \{ \ut>0 \} \, 
		\end{align}
		follows.
		Moreover, we have $\mathbb{\tilde{P}} \circ (\unt^{\delta})^{-1} = \mathbb{P} \circ (\un^{\delta})^{-1} = \Lambda^{\delta}$. In particular $\unt^{\delta}$ is strictly positive $\ti{\mathbb{P}}$-almost surely.		
		\newline
		\noindent
		Due to Lemma~\ref{sat:IdentM0} we have for all $\phi \in H^{1}_{per}(\ort)$
		\begin{align}\notag
			\intsp{(\ut(t)-\unt^{\delta}) \phi} \dx
			&= \intt \intort \ut^{n} \ut_{xxx} \phi_{x} \dx \ds 
			+ (\Cstr + S)\intt \intsp \ut^{n-2}  \uxt \phi_{x} \dx \ds
			\\ 
			&\quad- \sum_{k \in \Z}\int_{0}^{t}\intsp{ \lal\gl\ut^{\frac{n}{2}}\phi_{x}} \dx \dblt \, 
		\end{align}  
		which is the weak formulation \eqref{eq:MainResultWeakForm}.
		The claimed positivity properties of $\ut$ are shown in appendix \ref{app:Pos}.
		\\
		To establish estimate \eqref{est:MainResult}, we rewrite \eqref{est:EnEnt} for $q\ge1$ in the form 
		\begin{align}\notag	 \label{mainproofEnEntEst}
			&\Ew{\sup_{t\in [0,T]} \left(\frac{1}{2}\intort  (\uext)^{2} \dx \right)^{q} 
				} 
			+
			\Ew{\left(\intsptiT{\chi_{[\ut > 0]} (\uet)^{n} (\ti{p}^{\e}_{x})^{2}}\right)^{q}}
			\\ \notag
			&+ 
			\Ew{\left(\intsptiT{(\uet)^{n-4} (\uext)^{4}}\right)^{q}}+
			\Ew{\left(\intsptiT{(\uet)^{n-2}(\uexxt)^{2}}\right)^{q}}
			\\ 
			&+
			\Ew{\left(\intsptiT{(\uexxt)^{2}}\right)^{q}}
			\le  C(\unt,\delta,q,T)\, ,
		\end{align}
		where we neglected terms involving the potential. 
		Using Fatou's lemma, we have
		\begin{align}\notag	 \label{mainprooffatou}
			&\mathbb{E} \Bigg[ \liminf_{\e \rightarrow 0} 
			\Bigg\{
			\sup_{t\in [0,T]} \left(\frac{1}{2}\intort  (\uext)^{2} \dx \right)^{q}
			+
			\left(\intsptiT{ \chi_{[\ut > 0]} (\Jet)^{2}}\right)^{q}
			\\ \notag
			&\quad+ 
 \left(\intsptiT{\left((\uet)^{\frac{n}{4}}\right)_{x}^{4}} \right)^{q}
			\\ \notag			
			&\quad+
			 \left(\intsptiT{(\uet)^{n-2}(\uexxt)^{2}}\right)^{q}
			+
			\left(\intsptiT{(\uexxt)^{2}} \right)^{q} \Bigg\} \Bigg]
			\\ \notag 
			&\le
			\liminf_{\e \rightarrow 0}  \mathbb{E} \Bigg[ 
			\sup_{t\in [0,T]} \left(\frac{1}{2}\intort  (\uext)^{2} \dx \right)^{q}  
			+
			\left(\intsptiT{\chi_{[\ut > 0]} (\Jet)^{2}}\right)^{q}
			\\ \notag
			&\quad+ 
			 \left(\intsptiT{\left((\uet)^{\frac{n}{4}}\right)_{x}^{4}} \right)^{q}
			\\ 		
			&\quad+ 
			\left(\intsptiT{(\uet)^{n-2}(\uexxt)^{2}}\right)^{q}
			+
			\left(\intsptiT{(\uexxt)^{2}}\right)^{q}\Bigg]
			\le  C(\unt, \delta,q,T)\,.
		\end{align}
		Using the lower semi-continuity of the $\lpr{\infty}{0,T;\lpr{2}{\ort}}$ norm with respect to  
		convergence in the sense of distributions 
		(for $\uext \rightharpoonup \uxt$)
		and
		the lower semi-continuity of the $\lpr{2}{\Ot}$ norm or the $\lpr{4}{\Ot}$ norm for $\Jet$ and $\uexxt$ or for $\left((\uet)^{\frac{n}{4}}\right)_{x}$, respectively, we find
		\begin{align}\label{fatou1} \notag	 
			&\Ew{\sup_{t\in [0,T]} \left(\frac{1}{2}\intsp \uxt^{2} \dx \right)^{q} 
			} 
			+
			\Ew{\left(\intsptiT{ \ut^{n} \ut_{xxx}^{2}}\right)^{q}}
			\\ 
			&+ 
			\Ew{\left(\intsptiT{\left(\ut^{\frac{n}{4}}\right)_{x}^{4}}\right)^{q}}
			+\Ew{ \left(\intT \intort \ut_{xx}^{2} \dx \ds \right)^{q}}
			\le  C(\unt, \delta,q,T)\, .
		\end{align}
		To obtain also the estimate
		\begin{align}
			\Ew{\left(\intT \intort \ut^{n-2} \ut_{xx}^{2}\dx \ds\right)^{q}} \le C(\unt, \delta,q,T) \, ,
		\end{align}
		we use the uniform bound of $\Ew{\left(\intT \intort (\uet)^{n-2} (\uet)_{xx}^{2}\dx \ds\right)^{\ti{q}}} $ for any $\ti{q} \ge 1$, and the weak convergence of  $\ut^{\frac{n-2}{2}} \ut_{xx}$ in $\lpr{2\ti{q}}{\Omega; \lpr{2}{\Ot}}$. The identification of the limit follows by the weak convergence of $\uexxt$ in $\lpr{2}{\Omega;\lpr{2}{\Ot}}$ and the strong convergence of $\uet$ in $\lpr{2}{\Omega; C^{\ti{\gamma}, \frac{\ti{\gamma}}{4}}(\bar{\Ot}) }$.
		\\ \vphantom{I}
	\end{myproof}
	
	To conclude the first part of this paper, we prove Proposition~\ref{prop:Bernis}. The additional regularity resulting from it will be crucial to obtain qualitative results like finite speed of propagation in \cite{GruenKlein24FSOP}.
	\\
	\begin{myproof}{ of Proposition~\ref{prop:Bernis}}
		We check the assumptions of Lemma~B.1. in \cite{GiacShish2005}.
		For a given $\omega \in \Omega$ we take a representative $\ut(\omega, \cdot , \cdot) \in C^{\tilde{\gamma}, \frac{\ti{\gamma}}{4}}(\bar{\Ot})$.
		For $s\in [0,T]$ we further consider sets $K \subset\subset \{ \ut(\omega,s,\cdot)>0 \}$. As $\ut(\omega, \cdot , \cdot) \in C^{\tilde{\gamma}, \frac{\ti{\gamma}}{4}}(\bar{\Ot})$, we have $K \subset\subset \{ \ut(\omega,t,\cdot)>0 \}$ for $t$ sufficiently close to $s$. This implies
		\begin{align}\label{bernis1alt}
			\int_{K} \ut_{xxx}^{2} \dx  		
			\le
			C(K)
			\int_{ \{\ut(\omega,\tau,\cdot) > 0\}}
			\ut^{n} \ut_{xxx}^{2} \dx 
			< \infty
		\end{align}
		for almost all $\tau \in (s-\delta, s+\delta)$,
		as otherwise, this would contradict $\intt \intort \ut^{n} (\ut_{xxx})^{2} < \infty$. 
		Thus, we have found that $\ut \in H^{3}_{loc}(\{ \ut > 0
		\})$ for almost all $t$ and $\ti{\mathbb{P}}$-almost surely.
		Nonnegativity of $\ut$ follows 
		from the strict positivity of $\uet$, $\e \in (0,1)$, $\tilde{\mathbb{P}}$-almost surely, combined with the uniform convergence $\uet \to \ut$, $\tilde{\mathbb{P}}$-almost surely,
		$\uxt(\omega, s, ) \in \lpr{2}{\cO}$ $\tilde{\mathbb{P}}$-almost surely and almost everywhere in $[0,T]$ follows from \eqref{est:MainResult}.
		Finally, 
		since  \eqref{fatou1} implies $||\ut||_{H^{2}_{per}(\ort)} < \infty$ $\ti{\mathbb{P}}$-almost surely and almost everywhere in $[0,T]$, using Sobolev embedding, 
		we have $\ut(\omega,s,\cdot) \in C^{1}(\ort)$. 
		\\ \vphantom{I}
	\end{myproof}

	\section{Passage to the limit $\delta \rightarrow 0$}\label{sec:LimitInit}
	This section is devoted to the passage to nonnegative initial data. We note that it was the shift of the initial data that enabled the entropy estimate and thus, to establish Proposition~\ref{prop:Bernis}. Losing the positivity of initial data, we do not have the entropy estimate at our disposal any longer. Thus, we will rely solely on the energy estimate 
	for the passage $\delta \to 0$.  
	As the techniques we use are very similar to those in the first part of the paper, we will be much more concise here. 
	In the following we associate with each $\delta\in (0,1)$ a martingale solution $((\ti{\Omega}^{\delta}, \ti{\mathcal{F}}^{\delta}, (\ti{\mathcal{F}^{\delta}_{t}})_{t \ge 0}, \ti{\mathbb{P}}^{\delta}),\ut^{\delta}, \ti{W}^{\delta} )$ of \eqref{eq:stfe-Strat-equ1} as constructed in Theorem~\ref{theo:MainResult} with initial data given by the law of $\unt^{\delta} := \unt + \delta$, where $\unt \in L^{2}(\Omega; H^{1}_{per}(\cO))$ is nonnegative almost surely. In particular, these laws satisfy (A4$^{\delta}$).
	Let us now give the precise statement of the second theorem in this paper. 
	Once more, we emphasize that the subsequent result remains true if the sequence $(\ut^{\delta})_{\delta\in (0,1)}$ were constructed via Theorem~\ref{theo:MainResult} replacing (A3) by the weaker assumption (A3)$^*$ in Theorem~\ref{MainResultFG}.
	\begin{the1}\label{theo:ConvS}
		Let (A1)-(A3) and (A4) be satisfied and let $T>0$ be given. 
		There exist a stochastic basis $(\hat{\Omega},\hat{\F},(\hat{\F}_{t})_{t\in[0,T]},\hat{\mathbb{P}})$ with a complete, right-continuous filtration, an $\hat{\mathcal{F}}_{t}$-adapted Q-Wiener process $\hat{W} = \sum_{k \in \Z}\lal\gl \blh$, 
		a random variable 
		$\unh \in \lpr{2}{\hat{\Omega};H^{1}_{\text{per}}(\cO)}$, and
		a continuous $L^{2}(\cO)$-valued, $\hat{\mathcal{F}}_{t}$-adapted process 
		$\uh $ contained in $ \lpr{2}{\hat{\Omega};L^{2}(0,T;\Hsobper{1}{\cO})} \cap$ $ L^{2}(\hat{\Omega}; C^{\hat{\gamma}, \frac{\hat{\gamma}}{4}}(\Otb))$, $\hat{\gamma}<1/2$, with the property $\hat{u}_{xxx} \in L_{loc}^{2}(\{ \ut > 0 \})$ for $ d\hat{\mathbb{P}}\otimes dt$-almost all $(\hat{\omega},t) \in \hat{\Omega} \times [0,T]$,
		such that the following holds:
		\begin{enumerate}
			\item \label{theo:ConvS1}
			$\uh$ and $\unh$ are $\hat{\mathbb{P}}$-almost surely nonnegative and we have  
			$
			\Lambda^{0} = \mathbb{\hat{P}} \circ \hat{u}_{0}^{-1}
			$,
			\item\label{theo:ConvS2}
			for $t \in [0,T]$ and all
			$\phi \in H^{1}_{\text{per}}(\cO)$
			\begin{align}\notag \label{theo:ConvS3}
				\intsp{(\uh(t)-\unh) \phi} \dx
				&= \intt \int_{\{\uh(s,\cdot) >0 \}} \uh^{n} \uh_{xxx} \phi_{x} \dx \ds 
				- (\Cstr+S)\intt \intsp \uh^{n-2} \uxh \phi_{x} \dx \ds
				\\ 
				&\quad - \sum_{k \in \Z}\int_{0}^{t}\intsp{ \lal\gl\uh^{\frac{n}{2}}\phi_{x}} \dx \dblh
			\end{align}  
			holds $\hat{\mathbb{P}}$-almost surely,
			\item \label{theo:ConvS4}
			for arbitrary $q \ge 1$, $\uh$ satisfies the estimate  
			\begin{align}\label{est:ConvS1} \notag	 
				&\Ew{\sup_{t\in [0,T]} \left(\frac{1}{2} \intsp \uxh^{2}(t) \dx \right)^{q} 
				} 
				+
				\Ew{\left(\intT \intsp \betr{\left(\uh^{\frac{n+2}{6}}\right)_{x} }^{6} \dx \ds \right)^{q}}
				\\ 
				&
				+
				\Ew{\left(\intT \intsp \betr{\left(\uh^{\frac{n+2}{2}}\right)_{xxx} }^{2} \dx \ds \right)^{q}}
				+
				\Ew{\left(\intsptiT{\left(\uh^{\frac{n}{4}}\right)_{x}^{4}}\right)^{q}}
				\le  C(\unh,q,T).
			\end{align}
		\end{enumerate}
	\end{the1}
	
	The proof of Theorem~\ref{theo:ConvS} will be given below in Section~\ref{sec:proofMain2}.
	\subsection{Energy estimate}
	Let us establish an energy estimate. We will reuse the result from Lemma~\ref{lem:EnEnt}, where the estimates with respect to the energy were actually $\delta$-independent.  
	We define
	\begin{align}\label{def:energy}
		E[u] &:= \frac{1}{2} \intsp \ux^{2} \dx \,.
	\end{align}
	The following lemma holds.
	\begin{lem1}\label{lem:EnalphaEnt}
		Let $\delta > 0$ be given and let $\ut^{\delta}$ be a solution to \eqref{eq:stfe-Strat-equ1} as constructed in Theorem~\ref{theo:MainResult}.
		There is a positive constant $C(\unt,q,T)$ independent of $\delta >0$ such that for any  
		$q>1$ 
		\begin{align}
			\label{est:EnMomentsdelta} \notag	 
			&\Ew{\sup_{t\in [0,T]} \left(E[\ut^{\delta}(t)]\right)^{q}} 
			+
			\Ew{\left(\intsptiT{\left((\ut^{\delta})^{\frac{n+2}{6}}\right)^{6}_{x}}\right)^{q}}
			\\ \notag
			&+
			\Ew{\left(\intsptiT{\left((\ut^{\delta})^{\frac{n+2}{2}}\right)^{2}_{xxx}}\right)^{q}}
			+ 
			\Ew{\left(\intsptiT{(\ut^{\delta})^{n-4} (\uxt^{\delta})^{4}}\right)^{q} }
			\\
			&+
			\Ew{\left(\intsptiT{(\ut^{\delta})^{n-2}(\uxxt^{\delta})^{2} }\right)^{q} }
			\le  C(\unt,q,T).
		\end{align}
	\end{lem1}
	\proof
	We observe that from Theorem~\ref{theo:MainResult} and Proposition~\ref{prop:Bernis} the estimate follows with a $\delta$-dependent right-hand side. 
	Neglecting the terms derived by the entropy estimate and observing that \eqref{energy:initialdata} is $\delta$-independent after the passage to the limit $\e \to 0$, we may derive an $\delta$-independent estimate mimicking the arguments in the proof of Theorem~\ref{theo:MainResult}. 	
	Now we may proceed in the spirit of the proof of Theorem~\ref{theo:MainResult}, cf. \eqref{mainprooffatou}, to find
	\begin{align}\notag	 
		&\Ew{\sup_{t\in [0,T]} E[\ut^{\delta}(t)]^{q} } 
		+
		\Ew{\left(\intsptiT{ (\ut^{\delta})^{n} (\ut^{\delta})_{xxx}^{2}}\right)^{q}}
		\\ 
		&+ 
		\Ew{\left(\intsptiT{\left((\ut^{\delta})^{\frac{n}{4}}\right)_{x}^{4}}\right)^{q}}
		+
		\Ew{\left(\intsptiT{(\ut^{\delta})^{n-2}(\uxxt^{\delta})^{2}}\right)^{q}}
		\le  C(\unt,q,T)\, .
	\end{align}
	Invoking additionally Proposition~\ref{prop:Bernis} with $\zeta \equiv 1$ and using monotonicity, the result follows.
	\proofend
	
	\subsection{Compactness}
	In this section, we establish convergence -- as before we might have to change the underlying probability space. The setup and the arguments are very similar to those in Section~\ref{sec:Compactness} and we refer to the corresponding results there for more details.
	We recall $\Jt^{\delta} = - (\ut^{\delta})^{n} \ut^{\delta}_{xxx}$ and define for $\hat{\gamma} < \ti{\gamma}$
	\begin{align*}
		\cX_{\ut} &:=  C^{\hat{\gamma}, \frac{\hat{\gamma}}{4}}(\Otb) ,
		\quad
		\cX_{\uxt} :=  L^{2}(0,T;L^{2}(\cO))_{weak} ,
		\quad 
		\cX_{\Jt} :=  L^{2}(0,T;L^{2}(\cO))_{weak} ,
		\\
		\cX_{\Wt} &:=  C([0,T];L^{2}(\cO)) ,
		\quad
		\cX_{\unt} := H^{1}_{\text{per}}(\cO).
	\end{align*} 
	Moreover, $\tilde{\cX} := \cX_{\ut} \times \cX_{\uxt} \times \cX_{\Jt} \times \cX_{\ti{W}} \times \cX_{\unt}$.
	Establishing tightness by similar arguments as in \cite{FischerGruen} Lemma~5.2., we may apply Jakubowski's theorem \cite{Jakubowski1994} as in Proposition~\ref{sat:JakubowskiAnwendung}.
	\begin{prop1}\label{sat:JakubowskiAnwendungdelta}
		For subsequences of $\ut^{\delta}, \uxt^{\delta}, \Jt^{\delta} , \unt^{\delta}$, and $\Wt^{\delta}$, there exist a probability space $(\hat{\Omega}, \hat{\mathcal{F}}, \hat{\mathbb{P}})$, a sequence of $L^{2}(\cO)$-valued stochastic
		processes $\hat{W}^{\delta}$ on $(\hat{\Omega}, \hat{\mathcal{F}}, \hat{\mathbb{P}})$, sequences of random variables
		$\uh^{\delta} \colon \hat{\Omega}\rightarrow C^{\hat{\gamma},\frac{\hat{\gamma}}{4}}(\Otb )$, 
		$\hat{w}^{\delta}\colon \hat{\Omega} \rightarrow L^{2}(0,T;L^{2}(\cO))$,
		$\Jh^{\delta}\colon \hat{\Omega} \rightarrow \lpr{2}{0,T;\lpr{2}{\cO}}$, 
		$\unh^{\delta}\colon \ti{\Omega} \rightarrow H^{1}_{\text{per}}(\cO)$, 
		random variables
		$\uh \in L^{2}(\ti{\Omega};C^{\hat{\gamma},\frac{\hat{\gamma}}{4}}(\Otb ))$, 
		\\
		$\hat{w} \in  L^{2}(\ti{\Omega};L^{2}(0,T;L^{2}(\cO)))$, 
		$\Jh \in L^{2}(\ti{\Omega};\lpr{2}{0,T;\lpr{2}{\cO}})$, 	
		$\unh \in L^{2}(\ti{\Omega};H^{1}_{\text{per}}(\cO))$,
		as well as an $L^{2}(\cO)$-valued process $\hat{W}$ such that
		\begin{enumerate}
			\item
			for all $\delta \in  (0,1)$  the law of 
			$(\uh^{\delta},\hat{w}^{\delta}, \Jh^{\delta},\hat{W}^{\delta}, \unh^{\delta})$
			on $\cX$ with respect to  the measure $\hat{\mathbb{P}}$ equals the law of 
			$(\ut^{\delta}, \uxt^{\delta}, \Jt^{\delta},\Wt^{\delta}, \unt^{\delta})$
			with respect to  $\tilde{\mathbb{P}}^{\delta}$.
			\item 
			as $\delta \rightarrow 0$, the sequence 
			$(\uh^{\delta},\hat{w}^{\delta}, \Jh^{\delta},\hat{W}^{\delta}, \unh^{\delta})$ 
			converges
			$\hat{\mathbb{P}}$-almost surely to \newline
			$(\uh, \hat{w}, \Jh, \hat{W}, \unh)$
			in the topology of $\tilde{\cX}$.		
		\end{enumerate}
	\end{prop1}
	
	We identify the random variables and processes from Proposition~\ref{sat:JakubowskiAnwendungdelta} as follows. 
	
	\begin{lem1}\label{lem:IdentifikationFormeldelta}
		We have 
		$
		\hat{w}^{\delta} = \uh^{\delta}_{x} 
		$,
		as well as
		$
		\Jh^{\delta} =- (\uh^{\delta})^{\frac{n}{2}} \uh^{\delta}_{xxx} 
		$
		$\mathbb{\hat{P}}$-almost surely.
	\end{lem1}
	
	We consider normal filtrations $(\ti{F}_{t})_{t \ge 0 }$ and $(\ti{F}^{\delta}_{t})_{t \ge 0}$ given by
	\begin{align}\label{def:Ftdelta}
		\hat{F}_{t} := \sigma(\sigma(r_{t}\uh,r_{t}\hat{W}) \cup \{N \in \hat{F}:\hat{\mathbb{P}}(N)=0
		\} \cup \sigma(\unh))
	\end{align}  
	and
	\begin{align}\label{def:Fdelta}
		\hat{F}^{\delta}_{t} := \sigma(\sigma(r_{t}\ut^{\delta},r_{t}\hat{W}^{\delta}) \cup \{N \in \hat{F}:\hat{\mathbb{P}}(N)=0
		\} \cup \sigma(\unt^{\delta})) \, .
	\end{align}  
	As in Lemma~\ref{lem:WienerProcesses} we have
	\begin{lem1}\label{lem:WienerProcessesdelta}
		The processes $\hat{W}^{\delta}$ and $\hat{W}$ are Q-Wiener processes which are adapted to the filtrations $(\hat{F}^{\delta}_{t})_{t \ge 0}$ and $(\hat{F}_{t})_{t \ge 0 }$, respectively. We have
		\begin{align}\label{lem:WienerProcesses-h1delta}
			\hat{W}^{\delta}(t) = \sum_{k \in \Z} \lal \hat{\beta}_{k}^{\delta}(t)\gl
		, \quad 
			\hat{W}(t) = \sum_{k \in \Z} \lal \hat{\beta}_{k}(t)\gl
		\end{align}
		with families $(\hat{\beta}_{k}^{\delta})_{k \in \Z}$ and $(\hat{\beta}_{k})_{k \in \Z}$ of i.i.d. standard Brownian motions with respect to  $(\hat{F}^{\delta}_{t})_{t\ge 0}$ and $(\hat{F}_{t})_{t\ge 0}$, respectively.
	\end{lem1}
	The proof of the next lemma is identical to the proof of the corresponding result in Lemma~\ref{lem:Identwt}.
	\begin{lem1}\label{lem:Identwh}
		We have $\mathbb{\hat{P}}$-almost surely
		$
		\hat{w} =  \uh_{x} .
		$		
	\end{lem1}   
	While technically the identification of the pseudo-flux $\Jh$ is in line with Lemma~\ref{lem:weakconvdiss} (we therefore omit the proof of Lemma~\ref{lem:weakconvdissdelta}), we note that the identification is only possible on the set where solutions are strictly positive. However, as we do not have the entropy estimate at our disposal, the set $\{ \uh = 0\}$ is not necessarily a null set in $\Ot$, cf. also Remark~\ref{bem:ut-nullset}. On the other hand, we have the subsequent Lemma~\ref{lem:dissnullsetdelta} to cover the convergence behavior of $\Jh^{\delta}$ in a weak sense on $\{ \uh = 0\}$.
	\begin{lem1}\label{lem:weakconvdissdelta}
		The weak derivative $\uh_{xxx}$ satisfies $\uh_{xxx} \in \lpr{2}{\{\uh > r}$ for any $r>0$, $\hat{\mathbb{P}}$-almost surely.
		Moreover, we have $\Jh^{\delta} \chi_{\{\uh > 0\}} = -(\uh^{\delta})^{\frac{n}{2}} \hat{u}^{\delta}_{xxx} \chi_{\{\uh > 0\}}\rightharpoonup -\uh^{\frac{n}{2}} \uh_{xxx} \chi_{\{\uh > 0\}}$ weakly in $\lpr{2}{\Ot}$ $\tilde{\mathbb{P}}$-almost surely, and
		\begin{align}\label{weakconvdiss1delta}
			\intT \int_{\{\uh(s,\cdot) >0\}} (\ut^{\delta})^{n} \Jh^{\delta}_{x} \phi \dx\ds \rightarrow -\intT \int_{\{\uh(s, \cdot) >0\}} \uh^{n} \uh_{xxx} \phi \dx\ds \,,
		\end{align}
		for all $\phi \in \lpr{2}{\Ot}$, $\hat{\mathbb{P}}$-almost surely.  
	\end{lem1}
	
	\begin{lem1}\label{lem:dissnullsetdelta}
		For $\phi \in \Hsob{1}{\cO}$ we have
		\begin{align}
			\intT \int_{\{\uh(s,\cdot) = 0\}} (\uh^{\delta})^{n} \uh^{\delta}_{xxx} \phi_{x} \dx\ds \rightarrow 0 \,
		\end{align}
		for $\delta \rightarrow 0$ $\hat{\mathbb{P}}$-almost surely.
	\end{lem1}
	\proof
	We use Hölder's inequality to obtain
	\begin{align}\notag
		\intT \int_{\{\uh(s,\cdot) = 0\}} |(\uh^{\delta})^{n} \uh^{\delta}_{xxx} \phi_{x}| \dx\ds =	\intT \int_{\{\uh(s,\cdot) = 0\}}  |(\uh^{\delta})^{n/2} (\uh^{\delta})^{n/2} \uh^{\delta}_{xxx} \phi_{x}| \dx\ds
		\\ 
		\le 
		\left(\intT \intort (\uh^{\delta})^{n} (\uh^{\delta}_{xxx})^{2} \dx\ds \right)^{\frac{1}{2}} \left(\intT \int_{\{\uh(s,\cdot)=0 \}} (\uh^{\delta})^{n} \phi_{x}^{2} \dx\ds \right)^{\frac{1}{2}} \, .
	\end{align}
	Since the first term is uniformly bounded and the second one converges to zero as $\delta \rightarrow 0$ $\hat{\mathbb{P}}$-almost surely, the result follows.
	\proofend
	
	By the same reasoning as in Lemma~\eqref{lem:weakconvcorr}, we get -- using $n \in (2,3)$ --
	\begin{lem1}
		We have
		\begin{align}
			\intT \intsp (\uh^{\delta})^{n-2} \uh^{\delta}_{x} \phi_{x} \dx\ds\rightarrow \intT \intsp \uh^{n-2} \uxh \phi_{x} \dx\ds
		\end{align}
		$\hat{\mathbb{P}}$-almost surely for any $\phi \in \Hsobper{1}{\cO}$.
	\end{lem1}
	
	\subsection{Identification of the stochastic integral}		
	The arguments of section \ref{sec:IdentificationStochInt} apply with minor changes in the proofs to the situation at hand. We note in particular that the reasoning there was purely based on the energy estimate and the entropy estimate was not used in the course of the proofs. Thus, we have 
	\begin{lem1}\label{sat:IdentM0delta}
		It holds $\mathbb{\ti{P}}$-almost surely
		\begin{align}\notag
			\intsp{(\uh(t)-\unh) \phi} \dx \ds
			&= \intt \int_{\{\uh(s,\cdot) >0 \}} \uh^{n} \uh_{xxx} \phi_{x} \dx \ds 
			+ (\Cstr+S)\intt \intsp \uh^{n-2} \uxh \phi_{x} \dx \ds
			\\ 
			&\quad- \sum_{k \in \Z}\int_{0}^{t}\intsp{ \lal\gl\uh^{\frac{n}{2}}\phi_{x}} \dx \dblh \, 
		\end{align} 
		for all $\phi \in H^{1}_{per}(\cO)$.
	\end{lem1}
	
	\subsection{Proof of the main result}
	\label{sec:proofMain2}
	Combining the results established in Section~\ref{sec:LimitInit} so far, we can prove our main result.
	\begin{myproof}{ of Theorem~\ref{theo:ConvS}}
		From Proposition~\ref{sat:JakubowskiAnwendungdelta}, Lemma~\ref{lem:WienerProcessesdelta}, and the identifications in Lemma~\ref{lem:weakconvdissdelta} and Lemma~\ref{lem:Identwh}, the existence of a stochastic basis, a $Q$-Wiener process $\hat{W}$, as well as random variables $\uh \in \lpr{2}{\Omega;C^{\hat{\gamma},\frac{\hat{\gamma}}{4}}(\bar{\cO}_{T})} \cap \lpr{2}{\Omega;\lpr{2}{[0,T];H^{1}(\ort)}}$, $\Jh \in \lpr{2}{\Omega; \lpr{2}{\Ot}}$, and $\unh$ with 
		\begin{align}
			\Jh = -\uh^{\frac{n}{2}}\uh_{xxx}  \qquad \hat{\mathbb{P}}\text{-almost surely in } \{ \uh>0 \} \, ,
		\end{align}
		follows. Regarding initial data, we have 
		$\ut_{0}^{\delta} \rightarrow \ut_{0} $ pointwise and thus, using Slutzki's theorem,
		$$\mathbb{\hat{P}} \circ (\uh_{0}^{\delta})^{-1} = \mathbb{\ti{P}} \circ (\ut_{0}^{\delta})^{-1}  \rightarrow \mathbb{\tilde{P}} \circ \ut_{0}^{-1} = \mathbb{P} \circ u_{0}^{-1} = \Lambda^{0} \, .$$
		As also $\mathbb{\hat{P}} \circ (\uh_{0}^{\delta})^{-1} \rightarrow \mathbb{\hat{P}} \circ \uh_{0}^{-1}$, cf. Proposition~\ref{sat:JakubowskiAnwendungdelta}, by the uniqueness of weak limits in Polish spaces, we find 
		\begin{align}
			\mathbb{\hat{P}} \circ \uh_{0}^{-1} = \Lambda^{0}\, .
		\end{align}
		\noindent
		Due to Lemma~\ref{sat:IdentM0delta}, we have for all $\phi \in H^{1}_{per}(\ort)$
		\begin{align}\notag 
			\intsp{(\uh(t)-\unh) \phi} \dx \ds
			&= \intt \int_{\{\uh(s,\cdot) >0 \}} \uh^{n} \uh_{xxx} \phi_{x} \dx \ds 
			+ (\Cstr+S)\intt \intsp \uh^{n-2} \uxh \phi_{x} \dx \ds
			\\ 
			&\quad - \sum_{k \in \Z}\int_{0}^{t}\intsp{ \lal\gl\uh^{\frac{n}{2}}\phi_{x}} \dx \dblh
		\end{align}  
		which is the weak formulation \eqref{theo:ConvS3}.
		To establish \eqref{est:ConvS1}, we use \eqref{est:EnMomentsdelta} 		
		to get
		\begin{align}\notag	 \label{mainproof1delta}
			&\Ew{\sup_{t\in [0,T]} E[\uh^{\delta}(t)]^{q} } 
			+
			\Ew{\left(\intT \intsp \left((\uh^{\delta})^{\frac{n+2}{6}}\right)_{x}^{6} \dx \ds \right)^{q}}
			\\ 
			&+
			\Ew{\left(\intT \intsp \left((\uh^{\delta})^{\frac{n+2}{2}}\right)_{xxx}^{2} \dx \ds \right)^{q}}
			+  
			\Ew{\left(\intsptiT{\left( (\uh^{\delta})^{\frac{n}{4}}\right)_{x}^{4} } \right)^{q}}
			\le  C(\unh,q,T)\, .
		\end{align}
		Fatou's lemma then gives
		\begin{align}\notag	 
			&\mathbb{E} \Bigg[ \liminf_{\delta \rightarrow 0} 
			\Bigg\{
			\sup_{t\in [0,T]} E[\uh^{\delta}(t)]^{q}  
			+
			\left(\intsptiT{\left((\uh^{\delta})^{\frac{n+2}{6}}\right)_{x} ^{6}}\right)^{q}
			\\ \notag
			&\quad+
			\left(\intsptiT{\left((\uh^{\delta})^{\frac{n+2}{2}}\right)_{xxx}^{2}}\right)^{q}
			+ 
			\left(\intsptiT{ \left( (\uh^{\delta})^{\frac{n}{4}}\right)_{x}^{4} }\right)^{q} \Bigg\} \Bigg]
			\\ \notag			
			&\le
			\liminf_{\delta \rightarrow 0}  \mathbb{E} \Bigg[ 
			\sup_{t\in [0,T]} E[\uh^{\delta}]^{q}  
			+
			\left(\intsptiT{\left((\uh^{\delta})^{\frac{n+2}{6}}\right)_{x} ^{6}}\right)^{q}
			\\ 
			&\quad+
			\left(\intsptiT{\left((\uh^{\delta})^{\frac{n+2}{2}}\right)_{xxx}^{2}}\right)^{q}
			+ 
			\left(\intsptiT{\left( (\uh^{\delta})^{\frac{n}{4}}\right)_{x}^{4}}\right)^{q} \Bigg]
			\le
			C(\unh,q,T)\, .
		\end{align}
		Using the lower semi-continuity of the $\lpr{\infty}{0,T;\lpr{2}{\ort}}$ norm with respect to  
		convergence in the sense of distributions 
		(for $\uh_{x}^{\delta} \rightharpoonup \uh_{x}$) and the lower semi-continuity of the $\lpr{4}{0,T;\lpr{4}{\ort}}$ norm with respect to weak convergence, we find
		\begin{align}\label{fatou1delta} \notag	 
			&\Ew{\sup_{t\in [0,T]} \left(\frac{1}{2}\intsp \uxh^{2}(t) \dx \right)^{q} } 
			+ 
			\Ew{\left(\intsptiT{\left(\uh^{\frac{n}{4}}\right)_{x}^{4}}\right)^{q}}
			\le  C(\unh,q,T)\, .
		\end{align}
		For the estimate
		\begin{align}
			\Ew{\left(\intsptiT{\left(\uh^{\frac{n+2}{6}}\right)_{x} ^{6}}\right)^{q} 
				+
				\left(\intsptiT{\left(\uh^{\frac{n+2}{2}}\right)_{xxx}^{2}}\right)^{q}}
			\le C(\unh,q,T) \, 
		\end{align}
		we may argue similarly.
		From \eqref{mainproof1delta}, weak convergence of $\left( (\uh^{\delta})^{\frac{n+2}{6}}\right)_{x}$ in $\lpr{6\ti{q}}{\Omega; \lpr{6}{\Ot}}$ and of $\left((\uh^{\delta})^{\frac{n+2}{2}}\right)_{xxx}$ in $\lpr{2\ti{q}}{\Omega; \lpr{2}{\Ot}}$, for $\ti{q}>1$ follows. The identification of the limit functions is readily established. Therefore, we may once again use lower-semi continuity of the corresponding $L^{p}$-norms with respect to weak convergence to conclude.
		\\ \vphantom{I}
	\end{myproof}
		
	\appendix
	
	\section{Auxiliary results}\label{app:Aux}
	For the readers convenience, we will formulate the main result of \cite{GruenKlein24Exist}, i.e. the existence result of strictly positive solutions to \eqref{eq:stfe-Strat-equ} for positive initial data, in a form that fits to our setting. The notion of solution is the following. 
	
	\begin{def1}
		\label{DefinitionMartingaleSolution}
		Let $\Lambda$ be a probability measure on $H^1_{per}(\ort)$. We will call a triple
		\newline
		$((\Omega,\mathcal{F},(\mathcal{F}_t)_{t\geq 0},\mathbb{P}), u, W)$ a martingale solution to the stochastic thin-film equation \eqref{eq:stfe-Strat-equ} with initial data $\Lambda$ on the time interval $[0,T]$ provided
		\begin{itemize}
			\item[i)] $(\Omega,\mathcal{F},(\mathcal{F}_t)_{t\geq 0},\mathbb{P})$ is a stochastic basis with a complete, right-continuous filtration, 
			\item[ii)] $W$ satisfies (A2) with respect to $(\Omega,\mathcal{F},(\mathcal{F}_t)_{t\geq 0},\mathbb{P})$, 
			\item[iii)]
			the continuous $L^{2}(\cO)$-valued, $\ti{\mathcal{F}}_{t}$-adapted process
			$u\in L^2( \Omega;L^2((0,T); H^2_{per}(\ort)))\cap L^2(\Omega; C^{\gamma, \gamma/4}(\bar{\cO}_{T}))$ with $\gamma\in (0,1/2)$ is positive $\probab$-almost surely
			and $u_{xxx} \in L^{2}(\ort \times [0,T])$ $\probab$-almost surely,
			\item[iv)] there exists an $\mathcal F_0$-measurable $H^1_{per}(\ort)$-valued random variable $u_0$ such that $\Lambda=\probab\circ u_0^{-1}$, and the equation
			\begin{align}
				\label{eq:weakformula}
				\begin{split}
					\int_\domain u(t) \phi\,dx =&\int_\domain u_0 \phi\,dx
					+\int_0^t \int_\domain u^{n} (u_{xx}-\e F'(u))_x \phi_x \, dx \,ds
					\\&
					- \int_{0}^{t} \intort u^{n-2}u_{x} \phi_{x} \, dx \, ds
					\\&
					-\sum_{k \in \Z} \lambda_k \int_0^t \int_\domain u^{\frac{n}{2}}g_k \phi_x \,dx \,d\beta_k
				\end{split}
			\end{align}
			holds $\probab$-almost surely for all $t\in [0,T]$ and all $\phi\in H^1_{per}(\domain)$.
		\end{itemize}
	\end{def1}
	
	Note that for brevity we will also call $u$ as in Definition~\ref{DefinitionMartingaleSolution} a solution to the stochastic thin-film equation and $\un$ the initial data. 
	The existence result reads
	\begin{sat1}
		\label{MainResultFG}
		Let $\e >0$, $\delta >0$, $T>0$, and $0<\gamma<1/2$ be given.
		Furthermore, let the assumptions (A1), (A2), (A3), and (A4$^{\delta}$) be satisfied. 
		Then there exist a stochastic basis $(\Omega^{\e},\mathcal{F}^{\e},(\mathcal{F}^{\e}_t)_{t\geq 0},\mathbb{P}^{\e})$, a process $u^{\e}\in L^2(\Omega^{\e};L^2([0,T];H^1_{per}(\mathcal{O})))$, 
		as well as a Q-Wiener process $W^{\e}$
		such that $((\Omega^{\e},\mathcal{F}^{\e},(\mathcal{F}^{\e}_t)_{t\geq 0},\mathbb{P}^{\e}), u^{\e}, W^{\e})$ is a martingale solution to the stochastic thin-film equation \eqref{eq:stfe-Strat-equ} with initial data $\Lambda^{\delta}$ in the sense of Definition~\ref{DefinitionMartingaleSolution}. With  $\pe := -\uexx + \e \W'(\ue)$, we have for any $q \geq 1$ the additional bound
		\begin{align}\label{mainresultest}
			\nonumber
			&\mathbb{E}\Bigg[ \sup_{t\in [0,T]} E[\ue(t)]^{q}
			+ \left(\int_0^{T}  \int_\domain (\ue)^{n} | p^{\e}_x|^2\,dx \,ds \right)^{q} 
			\\&\quad \nonumber
			+\left(\int_0^{T}  \intort |\Delta \ue|^2 \, dx \,ds\right)^{q} 
			+\left(\int_0^{T}  \intort (\ue)^{-p-2}  |\ue_x|^2  \,dx \,ds\right)^{q} 
			\\&\quad \nonumber
			+\left(\int_0^{T} \intort (\ue)^{n-4}\left|\ue_{x}\right|^4\, dx \,ds\right)^{q}
			+\left(\int_0^{T} \intort (\ue)^{n-2}\left|\Delta \ue\right|^2\, dx \,ds\right)^{q}
			\\&\quad
			+\left(\int_0^{T} \intort (\ue)^{n-p-4}\left|\ue_{x}\right|^2\, dx \,ds\right)^{q}
			+\left(\int_0^{T} \intort (\ue)^{-2}\left|\ue_{x}\right|^2\, dx \,ds\right)^{q}\Bigg]
			\le C(u_{0}, \delta ,q, T) \, .
		\end{align}
	\end{sat1}
	
	In the course of the proof of the a-priori estimates, cf. Lemma~\ref{lem:EnEnt}, we will frequently use the following relations that come as a consequence of the choice of our basis for $\lpr{2}{\cO}$ in $\eqref{def:L2-Basis}$ and the assumption $\lal = \laml$, $k \in \N$.
	\begin{lem1}\label{lem:ONBrelations}
		Let assumption (A2) be satisfied. Then
		\begin{align}\notag \label{ONBrelations}
			\sum_{k \in \Z} \gl^{2} \lal^{2}
			=
			\frac{\lan^{2}}{L}
			+
			\sum_{k =1}^{\infty} \frac{2\lal^{2}}{L}\, , \quad
			\sum_{k \in \Z} \glx^{2} \lal^{2} 
			=
			\sum_{k =1}^{\infty} \frac{\lal^{2}k^{2}8\pi^{2}}{L^{3}}\, ,
			\\ \notag 
			\sum_{k \in \Z} \glx\gl \lal^{2} 
			=
			0 \, , \quad
			\sum_{k \in \Z} \glxx^{2} \lal^{2} 
			=
			\sum_{k =1}^{\infty} \frac{\lal^{2}k^{4}32 \pi^{4}}{L^{5}} \, ,
			\\ 
			\sum_{k \in \Z} \glx\glxx \lal^{2} 
			=
			0 \, ,\quad
			\sum_{k \in \Z} \glxx \gl \lal^{2} 
			=
			-
			\sum_{k =1}^{\infty} \frac{\lal^{2}k^{2}8\pi^{2}}{L^{3}} \, 
		\end{align}
		holds.
		In particular, we have 
		\begin{align}
			\frac{1}{2} \frac{n^{2}}{4} \sum_{k \in \Z} \lal^{2} \gl^{2} = \frac{1}{2} \frac{n^{2}}{4} \left(\frac{\lan^{2}}{L} + \sum_{k=1}^{\infty} \frac{2\lal^{2}}{L}\right) = \Cstr \, ,
		\end{align}
		cf. (A3).
	\end{lem1}
	
	Combining Young's inequality with basic calculus and Lemma~\ref{lem:ONBrelations}, we get the following estimate on the It\^o-term.
	\begin{lem1}
		We have for positive constants $C_{1}$ and $C_{2}$
		\begin{align}\notag \label{lem:adT6}
			&\frac{1}{2} \sum_{k \in \Z} \lal^{2}
			\intt \intsp \left((\ue)^{\frac{n}{2}} \gl\right)_{xx}^{2} \dx \ds
			\le
			\Cstr \intt \intsp (\ue)^{n-2} (\uexx)^{2} \dx \ds
			\\ \notag
			&+ \Cstr \left(\tfrac{(n-2)^{2}}{4}+\tfrac{1}{3} |(n-3)(n-2)|\right) \intt \intsp (\ue)^{n-4} (\uex)^{4} \dx \ds
			\\ 
			&+
			C_{1} \intt \intsp (\ue)^{n-2}(\ue_{x})^{2} \dx \ds + C_{2} \intt \intsp (\ue)^{n}\dx \ds \,.
		\end{align}
	\end{lem1}
	
	The following result is a straight forward consequence of the Gagliardo-Nirenberg inequality.
	\begin{lem1}\label{app:estun}
		Let $u \in \Hsob{1}{\Ot}$, $n<4$, and assume mass to be conserved. Then,
		\begin{align}
			\intT \intort u^{n} \dx \ds 
			\le
			C_{1} \e \intT\intsp u_{x}^{2} \dx \ds
			+ C_{2}(M)
			\, ,
		\end{align}
		where $M = \intsp u(0) \dx$ is the initial mass.
	\end{lem1}	

	The next result provides a positive lower bound of solutions $\ue$ given that the energy does not exceed a certain threshold. It was introduced in \cite{GruenKlein22} Lemma~4.5 and can be proven in a similar way as Lemma~4.1 in \cite{FischerGruen}.
	\begin{lem1}\label{lem:Minimum-ue}
		Let $\e \in (0,1)$ and consider for $u\in H^1(\ort,\R^+)$ the functional
		\begin{align}
			H_ {\e}[u] := \frac{1}{2}\int_{\cO}  \betr{\ux}^{2} + \e (u)^{-p} \dx \, .
		\end{align}
		Then there is a positive constant $C_p$ independent of $\e$ such that  
		\begin{align}
			\underset{x \in \cO}{\sup}  (u)^{-1} \le \left(\fint_{\cO} u \dx \right)^{-1} + C_{p} \e^{\frac{1}{2-p}} H_{\e}[u]^{\frac{2}{p-2}}.
		\end{align}
		If in addition $H_{\e}[u] \le \sigma^{-1}$ for $ \sigma \in (0,1)$, there is a positive constant $\bar C_p$ independent of $\e$ such that 
		\begin{align}
			\min_{x \in \cO} u \ge \bar{C}_{p} \, \e^{\frac{1}{p-2}} \sigma^{\frac{2}{p-2}} \, .
		\end{align}
	\end{lem1}	
	
	\section{It\^o's formula}\label{appendix:Itoformula}
	Using Lemma~\ref{lem:Minimum-ue} and a stopping time argument we may rigorously apply It\^o's Formula to suitable functionals corresponding to \eqref{def:energyeps}, \eqref{def:entropy}, and \eqref{def:a-Entropy-Operator}. We have the following. 
	\begin{lem1}\label{lem:appItoformula}
		Let $\delta >0$ be given and let $\sigma > 0$ be such that $ \delta \ge \bar{C}_{p} \, \e^{\frac{1}{p-2}} \sigma^{\frac{2}{p-2}}$. For $\us$ as defined in \eqref{def:ModifiedSolutions}  
		we have for $t \in [0,\Tmax]$
		\begin{align}\label{eq:ItoE1a}\notag
			\frac{1}{2}\int_{\cO} (\us(t))^{2}_ {x} \dx
			&=
			\frac{1}{2}\int_{\cO} (\un^{\delta})^{2}_{x} \dx
			+
			\sum_{ k \in \Z} \int_{0}^{t\wedge T_{\sigma}} \int_{\cO} (\us)_{x}   \lal (\us^{\frac{n}{2}}\gl)_{xx} 
			\dx \dble 
			\\ \notag
			&+
			\int_{0}^{t\wedge T_{\sigma}} \int_{\cO} ((\us)^{n} \psx + (\Cstr + S) \us^{n-2} (\us)_{x}) (\us)_{xxx} \dx \ds
			\\ 
			&+
			\frac{1}{2} \int_{0}^{t\wedge T_{\sigma}}  \sum_{ k \in \Z} \lal^{2}\int_{\cO}(\us^{\frac{n}{2}} \gl)^{2}_{xx} \dx \ds \, ,
		\end{align}  
		\begin{align}\label{eq:ItoE2a}\notag
			\e \int_{\cO} \W(\us(t)) \dx  
			&=
			\e \int_{\cO} \W(\un^{\delta}) \dx
			+
			\sum_{ k \in \Z} \int_{0}^{t\wedge T_{\sigma}} \int_{\cO} \e \W'(\us)  \lal (\us^{\frac{n}{2}} \gl)_{x} \dx \dble
			\\ \notag
			&+
			\int_{0}^{t\wedge T_{\sigma}} \int_{\cO}
			(-(\us)^{n} \psx - (C_{Strat}+S) \us^{n-2} (\us)_{x}) (\e \W'(\us))_{x}
			\dx \ds
			\\
			&+
			\frac{1}{2} \int_{0}^{t\wedge T_{\sigma}}  \sum_{ k \in \Z} \lal^{2}\int_{\cO} \e \W''(\us)(\us^{\frac{n}{2}} \gl)^{2}_{x} \dx \ds \, ,
		\end{align} 	
		\begin{align}\notag
			\label{eq:ItoG1}
			\mathcal{G} [\us(t)] 
			=
			\mathcal{G} [\un^{\delta}] 
			+
			&\sum_{ k \in \Z} \int_{0}^{t\wedge T_{\sigma}}\int_{\cO} 
			G'(\us) 
			\lal(\us^{\frac{n}{2}} \gl)_{x} 
			\dx \dble
			\\ \notag
			+
			&\int_{0}^{t\wedge T_{\sigma}} \int_{\cO} 
			(-(\us)^{n} \psx - (C_{Strat}+S) \us^{n-2} (\us)_{x})
			\left(G'(\us)\right)_{x}
			\dx \ds
			\\ 
			+
			&\frac{1}{2} \int_{0}^{t\wedge T_{\sigma}} \sum_{ k \in \Z} \lal^{2}\int_{\cO}G''(\us)(\us^{\frac{n}{2}} \gl)^{2}_{x} \dx \ds \,, 
		\end{align}
		and
		\begin{align}\notag
			\label{eq:ItoGa1}
			\Gmcas [\us(t)] 
			=
			\Gmcas [\un^{\delta}] 
			+
			&\sum_{ k \in \Z} \int_{0}^{t\wedge T_{\sigma}}\int_{\cO} 
			\left(
			\frac{1}{\alpha} (\us)^{\alpha} - \frac{1}{\alpha}
			\right)
			\lal(\us^{\frac{n}{2}} \gl)_{x} 
			\dx \dble
			\\ \notag
			+
			&\int_{0}^{t\wedge T_{\sigma}} \int_{\cO} 
			(-(\us)^{n} \psx - (C_{Strat}+S) \us^{n-2} (\us)_{x})
			\left(\frac{1}{\alpha} \us^{\alpha} - \frac{1}{\alpha}
			\right)_{x}
			\dx \ds
			\\ 
			+
			&\frac{1}{2} \int_{0}^{t\wedge T_{\sigma}} \sum_{ k \in \Z} \lal^{2}\int_{\cO}\us^{\alpha-1}(\us^{\frac{n}{2}} \gl)^{2}_{x} \dx \ds \,
		\end{align}
		for $\mathcal{G}$ as defined in \eqref{def:entropy} and $\mathcal{G}_{\alpha}$ as defined in \eqref{def:a-Entropy}.	
		\end{lem1}
		\proof
		First, we show that all the assumptions 
		of 
		Theorem~3.1 in \cite{Krylov} are satisfied in our setting.
		Let us first derive weak formulations as in (3.1) of \cite{Krylov} in the spaces  $V_{1}:=H^{2}_{\text{per}}(\cO)$, $V_{2}:=H^{1}_{\text{per}}(\cO)$, and $H:=\lpr{2}{\cO}$. 
		By Lemma~\ref{lem:Minimum-ue}, we see that for $t \in [0,T_{\sigma}]$
		\begin{align}\label{est:ue-positive}
			\ue(\cdot, t) \ge \bar{c}_{\e} \sigma^{\frac{2}{p-2}},
		\end{align}
		where $\bar{c}_{\e} := \bar{C}_{p} \, \e^{\frac{1}{p-2}}$.
		For the functions $\us$ as introduced in \eqref{def:ModifiedSolutions}, we get for all $\phi \in H^{1}_{\text{per}}(\cO)$
		\begin{align}\label{eq:ModifiedWeakFormulation}\notag
			\int_{\cO} \us (t)\phi \dx 
			=
			&\int_{\cO}\un^{\delta} \phi \dx + \int_{0}^{t \wedge T_{\sigma}} \int_{\cO} \left(-\us^{n}\, \psx - (\Cstr+S) \us^{n-2} \usx \right) \phi_{x} \dx \ds
			\\
			&- 
			\sum_{k \in \Z} \int_{0}^{t \wedge T_{\sigma}} \int_{\cO} \us^{\frac{n}{2}} \gl \phi_{x} \dx \dbl
		\end{align}
		on $[0,\Tmax]$.
		\\
		The mapping $-\us^{n} \psx - (\Cstr+S)  \us^{n-2} \usx$ satisfies 
		\begin{align}\notag \label{est:unpxL2}
			&\Ew{\intT\intort \us^{2n} \psx^{2} - (\Cstr+S)^{2}  \us^{2n-4} \usx^{2} \dx \ds} 
			\\ \notag 
			&\le 
			\Ew{\sup_{t\in [0,T]} ||\us(t)||^{2n}_{H^{1}(\ort)}}^{\frac{1}{2}} \Ew{\intT\intort \us^{n} \psx^{2} \dx \ds}^{\frac{1}{2}}
			+ C\Ew{   \intT \intort \us^{n-4} \usx^{4} \dx \ds} 
			\\ 
			&\quad+ C \Ew{\intT \intort \us^{3n-4} \dx \ds}
			\le C <\infty \, , 
		\end{align}
		i.e. $-\us^{n} \psx - (\Cstr+S)  \us^{n-2} \usx \in \lpr{2}{\Omega \times \Ot}$.
		\noindent
		Consequently, we may define the mapping
		\begin{align}
			\xi \in \lpr{2}{\Omega \times [0,\Tmax]; (H^{2}_{\text{per}}(\cO))'}
		\end{align} 
		by 
		\begin{align}
			v \mapsto \Ew{\int_{0}^{\Tmax}\int_{\cO} \left(-\us^{n} \psx - (\Cstr+S)  \us^{n-2} \usx \right)_{xx} v \dx \ds} \, .
		\end{align}
		Hence, by Riesz' representation theorem there is $\xiast \in \lpr{2}{\Omega \times [0,\Tmax]; H^{2}_{\text{per}}(\cO)}$ such that
		\begin{align}\label{RieszProjH2}
			\int_{\Omega}
			\int_{0}^{\Tmax}
			{}_{(H^{2}_{\text{per}}(\cO))'}
			\langle 
			\xi, \phi  
			\rangle_{H^{2}_{\text{per}}(\cO)}
			\ds
			\dif\mathbb{P}^{\e}
			=
			\int_{\Omega}
			\int_{0}^{\Tmax}
			\skp{\xiast}{\phi}_{\Hsob{2}{\cO}}
			\ds	
			\dif\mathbb{P}^{\e}
		\end{align} 
		holds for all $\phi \in \lpr{2}{\Omega \times [0,\Tmax]; H^{2}_{\text{per}}(\cO)}$.
		Similarly, we introduce $\fast$ with respect to  $H^{1}_{\text{per}}(\cO)$, i.e. $\fast$ solves 
		\begin{align}\label{RieszProjH1}
			\int_{\Omega}
			\int_{0}^{\Tmax}
			{}_{(H^{1}_{\text{per}}(\cO))'}
			\langle 
			f, \phi  
			\rangle_{H^{1}_{\text{per}}(\cO)} 
			\ds
			\dif\mathbb{P}^{\e}
			=
			\int_{\Omega}
			\int_{0}^{\Tmax}
			\skp{\fast}{\phi}_{\Hsob{1}{\cO}}
			\ds	
			\dif\mathbb{P}^{\e}
		\end{align}
		for every $\phi \in \lpr{2}{\Omega \times [0,\Tmax]; H^{1}_{\text{per}}(\cO)}$,
		where $f \in \lpr{2}{\Omega \times [0,\Tmax]; (H^{1}_{\text{per}}(\cO))'}$ is  defined via
		\begin{align}
			v \mapsto \Ew{\int_{0}^{\Tmax}\int_{\cO} \left(-\us^{n} \psx - (\Cstr+S)  \us^{n-2} \usx\right) v_{x} \dx \ds} \, . 
		\end{align}
		Riesz' representation theorem also shows
		\begin{equation}
			\Ew{
				\int_{0}^{\Tmax}\left(\Norm{\us}^{2}_{H^{1}(\cO)} 
				+
				\Norm{\fast}^{2}_{H^{1}(\cO)}	\right) \ds 
			} 
			\le 
			C(\e,T) \, ,
		\end{equation}
		due to \eqref{est:unpxL2}.
		Owing to the regularity of solutions, cf. Definition~\ref{DefinitionMartingaleSolution}, and the definition of the stopping times $T_{\sigma}$ in \eqref{def:stoppingTimes},   
		\begin{align}
			\Ew{\int_{0}^{\Tmax}
				\left(
				\Norm{\usx}^{2}_{H^{2}(\cO)} 
				+
				\Norm{\xiast}^{2}_{H^{2}(\cO)}	
				\right) \ds } 
			\le
			C(\e,T) 
		\end{align}
		holds as well. 
		\\
		Let us set $\ssl := (\lal \us^{\frac{n}{2}}(s, \cdot) \gl)_{x}$.
		By standard convolution arguments, we find the processes  $\xiast$, $\fast$, $\ssl$, and $(\ssl)_{x}$ to be predictable.
		Thus, we may rewrite \eqref{eq:ModifiedWeakFormulation} by means of \eqref{RieszProjH2} and obtain
		\begin{align}\label{weakformF1} \notag
			\skp{(\us)_{x}(t)}{\phi}_{\lpr{2}{\cO}} 
			=
			\skp{(\un^{\delta})_{x}}{\phi}_{\lpr{2}{\cO}}
			&+
			\int_{0}^{t \wedge T_{\sigma}} \skp{-\xi^{\ast} }{\phi}_{H^{2}(\cO)} \ds
			\\
			&+
			\sum_{ k \in \Z} \int_{0}^{t\wedge T_{\sigma}}  \skp{(\ssl)_{x}}{\phi}_{\lpr{2}{\cO}}  \dbl \, , 
		\end{align}
		where we have multiplied with $\phi_{x}$ for $\phi \in \Hsobper{2}{\cO}$ and integrated by parts. On the other hand, with \eqref{RieszProjH1} we get 
		\begin{align}\label{weakformF2} \notag
			\skp{\us(t)}{\phi}_{\lpr{2}{\cO}} 
			=
			\skp{\un^{\delta}}{\phi}_{\lpr{2}{\cO}}
			&+
			\int_{0}^{t \wedge T_{\sigma}} \skp{\fast}{\phi}_{H^{1}(\cO)} \ds
			\\ 
			&+
			\sum_{ k \in \Z} \int_{0}^{t\wedge T_{\sigma}}  \skp{\ssl}{\phi}_{\lpr{2}{\cO}} \dbl
		\end{align}
		for all $\phi \in \Hsobper{1}{\cO}$. Both formulations \eqref{weakformF1} and \eqref{weakformF2} hold for all $t \in [0,T]$. 
		Concerning the assumption
		\begin{equation}\label{Krylov:StochInt}
			\sum_{k\in \Z} \Ew{\int_{0}^{T} \Norm{\ssl}^{2}_{H} \ds } < \infty, 
		\end{equation}
		we get with Lemma~\eqref{ONBrelations} and Young's inequality 
		\begin{align}\notag
			\sum_{k \in \Z} \int_{\cO} \lal^{2} \abs{((\ue)^{\frac{n}{2}} \gl)_{x}}^{2} \dx
			&=
			\sum_{k \in \Z} \int_{\cO} \lal^{2} \abs{\frac{n}{2}(\ue)^{\frac{n-2}{2}} \uex \gl + (\ue)^{\frac{n}{2}} \glx}^{2} \dx
			\\ \notag 
			&=
			2\Cstr \int_{\cO} (\ue)^{n-2} (\uex)^{2} \dx + \frac{8 \pi^{2}}{L^{3}} \sum_{k=1}^{\infty}k^{2} \lal^{2} \int_{\cO} (\ue)^{n}  \dx \, ,
			\\ \notag 
			&\le
			2\Cstr \int_{\cO} (\ue)^{n-4} (\uex)^{4} \dx + C  \int_{\cO} (\ue)^{n}  \dx \, ,
		\end{align}
		which implies \eqref{Krylov:StochInt} due to the regularity shown for $\ue$ in \cite{GruenKlein24Exist} and Lemma~\ref{app:estun}.
		For $(\ssl)_{x}$ we use similar arguments and estimate
		\begin{align}\notag
			&\sum_{k \in \Z} \int_{\cO} \lal^{2} \abs{((\ue)^{\frac{n}{2}} \gl)_{xx}}^{2} \dx
			=
			\sum_{k \in \Z} \int_{\cO} \lal^{2} \abs{\left(\frac{n}{2}(\ue)^{\frac{n-2}{2}} \uex \gl + (\ue)^{\frac{n}{2}} \glx\right)_{x}}^{2} \dx
			\\ \notag 
			&=
			\sum_{k \in \Z} \int_{\cO} \lal^{2} \left|\frac{n}{2}\frac{n-2}{2}(\ue)^{\frac{n-4}{2}} (\uex)^{2} \gl 
			+ \frac{n}{2}(\ue)^{\frac{n-2}{2}} \uexx  \gl 
			+ \frac{n}{2}(\ue)^{\frac{n-2}{2}} \uex \glx
			\right.
			\\ \notag
			&\qquad\qquad\left.  
			+ \frac{n}{2}(\ue)^{\frac{n-2}{2}} \uex \glx 
			+ (\ue)^{\frac{n}{2}} (\gl)_{xx}\right|^{2} \dx
			\\  
			&\le
			C \left( \intsp  (\uex)^{n-4} (\uex)^{4} \dx 
			+ \int_{\cO} (\ue)^{n-2} (\uexx)^{2}  \dx
			+ \int_{\cO} (\ue)^{n-2} (\uex)^{2}  \dx
			+ \int_{\cO} (\ue)^{n}  \dx \right) \, 
		\end{align}
		to get
		\begin{align}\label{Krylov:StochInt1} \notag
			\sum_{k\in \Z} \Ew{\int_{0}^{T} \Norm{(\ssl)_{x}}^{2}_{H} \ds } 
			&\le 
			C \mathbb{E} \Bigg[ \intT \intsp  (\uex)^{n-4} (\uex)^{4} \dx \ds
			+  \intT\int_{\cO} (\ue)^{n-2} (\uexx)^{2}  \dx \ds 
			\\ 
			&+  \intT\int_{\cO} (\ue)^{n-2} (\uex)^{2}  \dx \ds
			+  \intT\int_{\cO} (\ue)^{n}  \dx \ds  \Bigg]
			\le 
			C(T,u_{0},\delta) \, .
		\end{align}
		Next, we state the operators we work with.
		To guarantee well-posedness on $H$ and continuity of corresponding Fr\'echet derivatives on $H \times H$, we use a cut-off function $\vartheta \in C^{\infty}(\R)$ 
		such that for an appropriate $\delta_{1} > 0$ 
		\begin{align}\label{def:CutoffFunc}
			\vartheta(x) =
			\begin{cases}
				\abs{x} &\qquad \text{for} \quad \abs{x} \ge \bar{c}_{\e}\sigma^{\frac{2}{p-2}}
				\\
				\in \R^{+}  &\qquad \text{for} \quad \abs{x} \in (\bar{c}_{\e}\sigma^{\frac{2}{p-2}} -\delta_{1}, \bar{c}_{\e}\sigma^{\frac{2}{p-2}})
				\\
				\bar{c}_{\e}\sigma^{\frac{2}{p-2}} - \delta_{1} &\qquad \text{for} \quad \abs{x} \le \bar{c}_{\e}\sigma^{\frac{2}{p-2}} - \delta_{1}
			\end{cases}
		\end{align}
		and 
		\begin{align}\label{est:CutoffFuncDer}
			\abs{\vartheta^{(s)}(x)} \le C(s) \delta_{1}^{-s} \quad, \, s\in (1,2) \, .
		\end{align}
		Note that \eqref{est:ue-positive} implies $\vartheta(\ue) =\ue$, $\vartheta'(\ue) =1$, and $\vartheta''(\ue) = 0$ on $[0,T_{\s}]$. By our assumptions on $\sigma$ we also have $\vartheta(\un^{\delta}) = \un^{\delta}$.
		We define for $u \in L^{2}(\cO)$
		\begin{align}\label{def:Energy-Operator}
			&E_{1}: u \mapsto \frac{1}{2}\int_{\cO} u^{2} \dx \, ,
			\\ \label{def:Potential-Operator}
			&E_{2}: u \mapsto \e \int_{\cO} \W (\vartheta(u)) \dx \, ,
			\\ \label{def:Entropy-Operator}
			&\mathcal
			G^{\vartheta}: u \mapsto \int_{\cO} G(\vartheta(u)) \dx \, ,
			\\ \label{def:a-Entropy-Operator}
			&\Gmcas^{\vartheta}: u \mapsto \int_{\cO} G_{\alpha}(\vartheta(u)) \dx \, ,
		\end{align}
		where $G$ is defined as in \eqref{def:entropy1}
		and 
		\begin{align}\label{def:a-Entrop-Operator1}
			G_{\alpha}(u) := \frac{1}{\alpha(\alpha-1)}u^{\alpha+1} - \frac{1}{\alpha}u+\frac{1}{\alpha+1} \, ,
		\end{align}
		cf. \eqref{def:a-Entropy}.
		For $E_{1}$, $E_{2}$, and $\mathcal{G}$ the Fr\'echet derivatives are readily computed.
		Due to the cutoff $\vartheta$, the assumptions i) to iv) of Theorem~3.1 in \cite{Krylov} are readily checked for the space $H$ and its dense subsets $V_{1}$ (in the case of $E_{1}$) and $V_{2}$ (in the case of $E_{2}$ and $\mathcal{G}^{\vartheta}$, and $\Gmcas^{\vartheta}$, respectively).
		Hence, we may choose the operators $E_{1}, E_{2}$, $\mathcal{G}^{\vartheta}$, and $\Gmcas^{\vartheta}$ to apply Theorem~3.1 in \cite{Krylov}, with respect to the weak formulation \eqref{weakformF1} in the first case and with respect to  \eqref{weakformF2} in the other cases.  
		The result follows.
		\proofend
					
		\section{Positivity properties} \label{app:Pos}
		We will show the positivity properties of solutions to \eqref{eq:stfe-Strat-equ1} that were stated in Theorem~\ref{theo:MainResult}. The proof is based on an argument given in Theorem~4.1 of \cite{Beretta1995}. The techniques used there require additional regularity, i.e. an $\alpha$-entropy estimate. 
		Let us define for $u \in L^{2}(\cO)$
		\begin{align}\label{def:a-Entropy}
			\Gmcas[u] :=  \int_{\cO} G_{\alpha}(u) \dx \, 	,
		\end{align}
		where 
		$G_{\alpha}(u)$ is defined in \eqref{def:a-Entrop-Operator1}.
		We have the following lemma.
		\color{black}
		\begin{lem1}\label{lem:EnalphaEnt1}
			Let $\delta > 0$ and $\ue$ be a solution to \eqref{eq:stfe-Strat-equ} in the sense of Definition~\ref{DefinitionMartingaleSolution} with initial data $\Lambda^{\delta}$.
			Then for $\alpha \in (\frac{1}{2}-n, 2-n)\backslash\{-1\}$ there is a positive constant $C(\un, \delta, \alpha, T)$ independent of $S >0$ and $\e >0$ such that for $t \in [0,T]$ 
			\begin{align} \notag \label{eq:EnalphaEnt1}
				\mathbb{E} \Bigg[
				&\mathcal{G}_{\alpha}[\ue(t)] 
				+ S\int_{0}^{T} \intort (\ue)^{n+\alpha-3}  (\uex)^{2} \dx \ds
				\\ \notag 
				&+
				\int_{0}^{T} \intort (\ue)^{n+\alpha-1}  \W''(\ue)(\uex)^{2} \dx \ds
				+ 
				\frac{1}{\bar{\gamma}^{2}} \intT \intort (\ue)^{n+1-2\bar{\gamma} +\alpha} ((\ue)^{\bar{\gamma}})_{xx}^{2} \dx \ds 
				\\ \notag	
				&+(-(\bar{\gamma}-1)^{2} + \frac{1}{3} (2\bar{\gamma} -n - \alpha -1)(n + \alpha -2) )
				\intT \intort (\ue)^{n + \alpha -3} (\uex)^{4} \dx \ds
				\Bigg]
				\\
				\le
				& C(\un, \delta, \alpha, T) \, ,  
			\end{align} 
		where $\bar{\gamma}$ satisfies 
		\begin{align}\label{def:alphaEntrGamma}
			\frac{1}{3} (t + 1 - \sqrt{(t-2)(1-2t)}) \le \bar{\gamma} \le \frac{1}{3} (t+1 + \sqrt{(t-2)(1-2t)})\,, \quad t = \alpha +n \, .
		\end{align}
		\end{lem1}
		\proof
		We apply It\^o's formula to stopped solutions $\us$ to \eqref{eq:stfe-Strat-equ} as defined in
		\ref{def:ModifiedSolutions}. By Lemma~\ref{lem:appItoformula}, we get for $t \in [0,T]$
		\begin{align}\notag
			\label{eq:ItoGa}
			\Gmcas [\us(t)] 
			=
			&\Gmcas [\un^{\delta}] 
			+
			\sum_{ k \in \Z} \int_{0}^{t\wedge T_{\sigma}}\int_{\cO} 
			\left(
			\frac{1}{\alpha} \us^{\alpha} - \frac{1}{\alpha}
			\right)
			\lal(\us^{\frac{n}{2}} \gl)_{x} 
			\dx \dbl
			\\ \notag
			+
			&\int_{0}^{t\wedge T_{\sigma}} \int_{\cO} 
			(-\us^{n} \psx - (C_{Strat}+S) \us^{n-2} (\us)_{x})
			\left(\frac{1}{\alpha} \us^{\alpha} - \frac{1}{\alpha}
			\right)_{x}
			\dx \ds
			\\ 
			+
			&\frac{1}{2} \int_{0}^{t\wedge T_{\sigma}} \sum_{ k \in \Z} \lal^{2}\int_{\cO}\us^{\alpha-1}(\us^{\frac{n}{2}} \gl)^{2}_{x} \dx \ds 
			=:  \Gmcas [\un^{\delta}]  + A_{1} + A_{2} + A_{3}.
		\end{align}
		
		Ad $A_{2}$: 
		\begin{align} \notag
			A_{2} &= - \intts \intort \us^{n}(-\usxx + \W'(\us) )_{x} \us^{\alpha-1} \usx \dx \ds
			\\ 
			&\quad \, \, - (\Cstr + S) \intts \intort \us^{n+\alpha-3} \usx^{2} \dx \ds
			=: A_{21} + A_{22}.
		\end{align}
		Here $A_{22}$ is a good term. For $A_{21}$, we find using integration by parts
		\begin{align}\label{EnalphaEnt1} \notag
			&- \intts \intort \us^{n}(-\usxx + \W'(\us) )_{x} \us^{\alpha-1} \usx \dx \ds
			\\ \notag
			&= 
			- \intts \intort \us^{n+\alpha-1} \usxx^{2}    \dx \ds
			- \intts \intort (n+\alpha-1) \us^{n+\alpha-2} \usx^{2} \usxx  \dx \ds
			\\\notag
			&\quad- \intts \intort \us^{n+\alpha-1}  \W''(\us)\usx^{2} \dx \ds
			\\ \notag
			&= 
			- \intts \intort \us^{n+\alpha-1} \usxx^{2}    \dx \ds
			+ \intts \intort \frac{(n+\alpha-1)(n+\alpha-2)}{3} \us^{n+\alpha-3} \usx^{4} \dx \ds
			\\
			&\quad- \intts \intort \us^{n+\alpha-1}  \W''(\us)\usx^{2} \dx \ds \, .
		\end{align}
		We note that the penultimate term is a good term if $\alpha \in (1-n, 2-n)$. To extend the range of admissible  $\alpha$, we rewrite
		\begin{align}\label{eq:subst_uxx}
			\usxx = \frac{1}{\bar{\gamma}^{2}} \us^{2-2\bar{\gamma}} (\us^{\bar{\gamma}})_{xx}^{2} - (\bar{\gamma} -1)^{2} \us^{-2} \usx^{4} + 2(\bar{\gamma} -1 ) \us^{-1} \usx^{2} \usxx   \, ,
		\end{align}
		where $\bar{\gamma}$ is a positive constant to be determined later on.  
		This idea has been presented in \cite{Beretta1995}, Proposition~2.1.
		Inserting \eqref{eq:subst_uxx} in \eqref{EnalphaEnt1} and neglecting the term involving the potential, gives 
		\begin{align}\notag \label{EnalphaEnt2}
			&- \intts \intort \us^{n+\alpha-1} \usxx^{2}    \dx \ds
			- \intts \intort (n+\alpha-1) \us^{n+\alpha-2} \usx^{2} \usxx  \dx \ds
			\\\notag
			&= 
			-\frac{1}{\bar{\gamma}^{2}} \intts \intort \us^{n+1-2\bar{\gamma} +\alpha} (\us^{\bar{\gamma}})_{xx}^{2} \dx \ds 
			+ (\bar{\gamma} -1)^{2} \intts \intort \us^{n + \alpha -3} \usx^{4} \dx \ds
			\\ \notag 
			&\quad + (2\bar{\gamma} -n -\alpha -1) \intts \intort \us^{n + \alpha -2} \frac{1}{3} ((\usx)^{3})_{x} \dx \ds
			\\ \notag 
			&= -\frac{1}{\bar{\gamma}^{2}} \intts \intort \us^{n+1-2\bar{\gamma} +\alpha} (\us^{\bar{\gamma}})_{xx}^{2} \dx \ds 
			\\ 
			&\quad + ((\bar{\gamma}-1)^{2} - \frac{1}{3} (2\bar{\gamma} -n - \alpha -1)(n + \alpha -2) )
			\intts \intort \us^{n + \alpha -3} \usx^{4} \dx \ds \, .
		\end{align} 
		The number $(\bar{\gamma}-1)^{2} - \frac{1}{3} (2\bar{\gamma} -n - \alpha -1)(n + \alpha -2)$ turns out to be negative if $\alpha \in (\frac{1}{2}-n, 2-n)$ and \eqref{def:alphaEntrGamma} holds.
		Thus, the terms derived by $A_{21}$ are good terms. 
		\\
		Ad $A_{3}$:
		\begin{align}\notag \label{EnalphaEnt3}
			&\frac{1}{2} 
			\int_{0}^{t\wedge T_{\sigma}} \int_{\cO} 
			\sum_{ k \in \Z} \lal^{2} \us^{\alpha-1} (\us^{\frac{n}{2}} \gl)^{2}_{x} \dx \ds 
			\\ \notag
			&=
			\frac{1}{2} 
			\int_{0}^{t\wedge T_{\sigma}} \int_{\cO} 
			\sum_{ k \in \Z} \lal^{2} \us^{\alpha-1} \left(\us^{\frac{n-2}{2}} \usx \gl + \us^{\frac{n}{2}} \glx \right)^{2} \dx \ds 
			\\ 
			&=
			\frac{1}{2} 
			\int_{0}^{t\wedge T_{\sigma}} \int_{\cO} 
			\sum_{ k \in \Z} \lal^{2} \us^{\alpha-1} \left(\us^{n-2} \usx^{2} \gl^{2} + \us^{n} \glx^{2}\right)\dx \ds 
			=: \operatorname{I} + \operatorname{II}\, . 
		\end{align}
		To absorb $\operatorname{I}$, we may use Lemma~\ref{ONBrelations} and the term $A_{22}$. We have for $n \in (2,3)$  
		\begin{align}
			n+\alpha -1 
			\begin{cases}
				\in (0,2)\, , \quad \, &\alpha > 1-n 
				\\ 
				\in (\alpha +1, 0]   \, , \quad \,  &\alpha \le 1-n \, . 
			\end{cases}
		\end{align}
		Thus, in the first case,  we may estimate $\operatorname{II}$ by means of Young's and Poincar\'e's inequalities to get for positive constants $c$ and $C$
		\begin{align}\notag
			\int_{0}^{t\wedge T_{\sigma}} \int_{\cO} 
			\sum_{ k \in \Z} \lal^{2} \us^{n+\alpha-1} \glx^{2} \dx \ds 
			&\le
			\sum_{k =1}^{\infty} \frac{\lal^{2}k^{2}8\pi^{2}}{L^{3}}
			\int_{0}^{t\wedge T_{\sigma}} \int_{\cO} 
			c \us^{2} + C\dx \ds  
			\\ 
			&\le
			C_{\text{Poinc}} C
			\int_{0}^{t\wedge T_{\sigma}} \int_{\cO} 
			\usx^{2} \dx \ds + C(T,L) \, .
		\end{align}
		Here, the first term is uniformly bounded with respect to $\e$ by the energy estimate Lemma~\ref{lem:EnEnt}. In the second case, we estimate for positive constants $c$ and $C$
		\begin{align} 
			\int_{0}^{t\wedge T_{\sigma}} \int_{\cO} 
			\sum_{ k \in \Z} \lal^{2} \us^{n+\alpha-1} \glx^{2} \dx \ds 			 
			\le
			\sum_{k =1}^{\infty} \frac{\lal^{2}k^{2}8\pi^{2}}{L^{3}}
			\int_{0}^{t\wedge T_{\sigma}} \int_{\cO} 
			c\,\us^{\alpha+1} + C \dx \ds \, . 
		\end{align}
		Hence, in this case we may argue via Gronwall's lemma. 
		Collecting all terms and taking expectation, by the martingale property of the It\^o integral, we find for $C(\un ,\delta, T)>0$ independent of $\e$ for any $t' \in [0,T]$
		\begin{align} \notag
			\mathbb{E} \Bigg[
			&\mathcal{G}_{\alpha}[\us(t')] 
			+ S \int_{0}^{t' \wedge T_{\sigma}} \intort \us^{n+\alpha-3}  \usx^{2} \dx \ds
			\\ \notag 
			&+
			\int_{0}^{t' \wedge T_{\sigma}} \intort \us^{n+\alpha-1} \W''(\us) \usx^{2} \dx \ds
			+ 
			\frac{1}{\bar{\gamma}^{2}} \int_{0}^{t' \wedge T_{\sigma}} \intort \us^{n+1-2\bar{\gamma} +\alpha} (\us^{\bar{\gamma}})_{xx}^{2} \dx \ds 
			\\ \notag	
			&+(-(\bar{\gamma}-1)^{2} + \frac{1}{3} (2\bar{\gamma} -n - \alpha -1)(n + \alpha -2) )
			\int_{0}^{t' \wedge T_{\sigma}} \intort \us^{n + \alpha -3} \usx^{4} \dx \ds \, .
			\Bigg]
			\\ 
			\le
			&\mathbb{E} \Bigg[
			\mathcal{G}_{\alpha}[\un]
			+
			C
			\int_{0}^{t' \wedge T_{\sigma}} \int_{\cO} 
			\us^{\alpha+1} \dx \ds \Bigg] + C(\un ,\delta, T) \, . 
		\end{align} 
		Thus, the strict positivity of the initial data and a Gronwall argument as used in the proof of Lemma~\ref{lem:EnEnt} give 
		\begin{align} \notag
			\mathbb{E} \Bigg[
			&\mathcal{G}_{\alpha}[\us(t')] 
			+ S \int_{0}^{t' \wedge T_{\sigma}} \intort \us^{n+\alpha-3}  \usx^{2} \dx \ds
			\\ \notag 
			&+
			\int_{0}^{t' \wedge T_{\sigma}} \intort \us^{n+\alpha-1} \W''(\us) \usx^{2} \dx \ds
			+ 
			\frac{1}{\bar{\gamma}^{2}} \int_{0}^{t' \wedge T_{\sigma}} \intort \us^{n+1-2\bar{\gamma} +\alpha} (\us^{\bar{\gamma}})_{xx}^{2} \dx \ds 
			\\ \notag 	
			&+(-(\bar{\gamma}-1)^{2} + \frac{1}{3} (2\bar{\gamma} -n - \alpha -1)(n + \alpha -2) )
			\int_{0}^{t' \wedge T_{\sigma}} \intort \us^{n + \alpha -3} \usx^{4} \dx \ds \, .
			\Bigg]
			\\
			\le &C(\un, \delta, \alpha, T). 
		\end{align} 
		As $t' \in [0,T]$ was arbitrary, by Lemma~ \ref{lem:ConvergenceStoppingTimes} and the same arguments as in the proof of Lemma~\ref{lem:EnEnt}, we establish \eqref{eq:EnalphaEnt1}.
		\proofend  
		
		\begin{bem1}
			It is possible to improve estimate \eqref{eq:EnalphaEnt1} to also gain control of moments of any order $q >1$, using the Burkholder-Davis-Gundy inequality. 
		\end{bem1}
		
		Using Lemma~\ref{lem:EnalphaEnt1}, Fatou's lemma, and lower semi-continuity, we establish easily
		\begin{kor1}\label{kor:EnEntalpha}
			Solutions $\ut$ of \eqref{eq:stfe-Strat-equ1} derived by Theorem~\ref{theo:MainResult} satisfy for a positive constant $C(\unt, \delta, \alpha, T)$ and for every $t \in [0,T]$
			\begin{align}
				\mathbb{E} \Bigg[  \intort \ut^{\alpha+1}(t) \dx
				+ \intT \intort \ut^{n + \alpha -3} \uxt^{4} \dx \ds \Bigg]
				\le C(\unt,\delta ,\alpha, T) \, , 
			\end{align}
			where $\alpha \in (\frac{1}{2}-n, 2-n) \backslash \{-1\}$.
		\end{kor1}
		Having established Lemma~\ref{lem:EnalphaEnt1} and Corollary~\ref{kor:EnEntalpha}, we may prove the positivity properties of solutions $\ut$ as claimed in Theorem~\ref{theo:MainResult}. 	
		\begin{lem1}\label{lem:positivityProp}
			Let $\ut$ be a solution to \eqref{eq:stfe-Strat-equ1} as constructed in Theorem~\ref{theo:MainResult} (in particular, initial data are strictly positive). Then, for almost all paths,  $\ut(\omega, t, \cdot )$ is strictly positive for almost all $t \in [0,T]$.  
		\end{lem1}
		\proof
		The proof is along the lines of Theorem~4.1 in \cite{Beretta1995}. 
		For almost all $\omega \in \Omega$ and, we have for $\alpha < -1$ and $\alpha \in (\frac{1}{2}-n, 2-n)$
		\begin{align}\label{lem:positivityProp0}
			\intsp \ut^{\alpha+1}(t) \dx < \infty \, 
		\end{align}
		for every $t \in [0,T]$ due to Lemma~\ref{lem:EnalphaEnt1}. Let us fix $\omega \in \Omega$ where \eqref{lem:positivityProp0} holds and abbreviate $\ut(\omega, \cdot, \cdot)= \ut(\cdot, \cdot) $. We find for almost every $t \in [0,T]$ and points $(x_{0},t)$ with $\ut(x_{0},t) = 0$ a positive constant $C(t)$ such that for any $y \in \cO$
		\begin{align}\label{lem:positivityProp1}
			\ut(y,t) = \ut(y,t) - \ut(x_{0},t) \le C(t) |y - x_{0}|^{\frac{3}{n+\alpha+1}}
			\, , 
		\end{align}
		where we used $\ut^{\frac{n+\alpha+1}{4}}(t) \in W^{1,4}(\cO) $ for almost all $t \in [0,T]$, cf. Lemma~\ref{lem:EnEnt}, and Morrey's inequality. 
		For such $x_{0} \in \cO$ and $t \in [0,T]$, where \eqref{lem:positivityProp1} holds, we find
		\begin{align}\label{lem:positivityProp2}
			\intsp \ut^{\alpha+1}(t) \dx 
			\ge C(t)\intsp  |y - 	x_{0}|^{\frac{3(\alpha+1)}{n+\alpha+1}} \, dy\, .
		\end{align}		
		As $n >2$, we may choose $\alpha$ such that $\frac{3(\alpha +1)}{n+\alpha+1}  < -1$.
		Then, the right-hand side of \eqref{lem:positivityProp2} is infinite, a contradiction to \eqref{lem:positivityProp0}. Thus, roots $x_{0} \in \cO$ of $\ut$ may exist only for times $t \in [0,T]$, where \eqref{lem:positivityProp1} is not valid. The result follows.
		\proofend

		\noindent
		\textbf{Acknowledgment.}
		L.K. has been supported by the Graduiertenkolleg 2339 \textit{IntComSin} ''Interfaces, Complex Structures, and Singular Limits'' of the Deutsche Forschungsgemeinschaft (DFG, German Research Foundation) -- Project-ID 321821685. The support is gratefully acknowledged.
		
		\bibliographystyle{abbrv}
		\bibliography{bibliography-Grenz}

\begin{thebibliography}{10}

\bibitem{Beretta1995}
E.~Beretta, M.~Bertsch, and R.~{Dal Passo}.
\newblock Nonnegative solutions of a fourth-order nonlinear degenerate
  parabolic equation.
\newblock {\em Arch. Ration. Mech. Anal.}, 129(2):175--200, 1995.

\bibitem{BernisFinite}
F.~Bernis.
\newblock {Finite speed of propagation and continuity of the interface for thin
  viscous flows}.
\newblock {\em Adv. Differential Equations}, 1(3):337--368, 1996.

\bibitem{BernisFinite2}
F.~Bernis.
\newblock {Finite speed of propagation for thin viscous flows when $2\leq
  n<3$}.
\newblock {\em C. R. Math. Acad. Sci. Paris}, 322(12):1169--1174, 1996.

\bibitem{BernisInequ1996}
F.~Bernis.
\newblock Integral inequalities with applications to nonlinear degenerate
  parabolic equations.
\newblock In T.~S. Angell, L.~P. Cook, R.~E. Kleinmann, and W.~E. Olmstead,
  editors, {\em Nonlinear problems in Applied Mathematics}, volume in Honour of
  I. Stakgold, pages 57--65. SIAM, 1996.

\bibitem{Bertozzi1996}
A.~L. Bertozzi and M.~Pugh.
\newblock The lubrication approximation for thin viscous films: {R}egularity
  and long-time behavior of weak solutions.
\newblock {\em Comm. Pure Appl. Math.}, 49(2):85--123, 1996.

\bibitem{BertschDalPassoGarckeGruen}
M.~Bertsch, R.~{Dal Passo}, H.~Garcke, and G.~Grün.
\newblock {The thin viscous flow equation in higher space dimensions}.
\newblock {\em Adv. Differential Equations}, 3(3):417 -- 440, 1998.

\bibitem{BreitHofmanova2016}
D.~Breit, E.~Feireisl, and M.~Hofmanová.
\newblock Incompressible limit for compressible fluids with stochastic forcing.
\newblock {\em Arch. Ration. Mech. Anal.}, 222(2):895--926, 2016.

\bibitem{Brezniak2007}
Z.~Brzezniak and M.~Ondrej\'at.
\newblock Strong solutions to stochastic wave equations with values in
  {R}iemannian manifolds.
\newblock {\em J. Funct. Anal.}, 253:449--481, 2007.

\bibitem{DalPassoGiacoGruenWTP}
R.~{Dal Passo}, L.~Giacomelli, and G.~Grün.
\newblock A waiting time phenomenon for thin film equations.
\newblock {\em Ann. Sc. Norm. Super. Pisa Cl. Sci.}, (4) 30(2):437--463, 2001.

\bibitem{DGGG}
K.~Dareiotis, B.~Gess, M.~V. Gnann, and G.~Grün.
\newblock Non-negative martingale solutions to the stochastic thin-film
  equation with nonlinear gradient noise.
\newblock {\em Arch. Ration. Mech. Anal.}, 242:179--234, 2021.

\bibitem{Dareiotis2023solutions}
K.~Dareiotis, B.~Gess, M.~V. Gnann, and M.~Sauerbrey.
\newblock Solutions to the stochastic thin-film equation for initial values
  with non-full support.
\newblock {\em arXive e-prints, arXiv:2305.06017}, 2023.

\bibitem{DavidovitchStone2005}
B.~Davidovitch, E.~Moro, and H.~A. Stone.
\newblock Spreading of viscous fluid drops on a solid substrate assisted by
  thermal fluctuations.
\newblock {\em Phys. Rev. Lett.}, 95(24):244505, 2005.

\bibitem{DebusscheHofmanovaVovelle}
A.~Debussche, M.~Hofmanov\'a, and J.~Vovelle.
\newblock Degenerate parabolic stochastic partial differential equations:
  {Q}uasilinear case.
\newblock {\em Ann. Probab.}, 44(3):1916--1955, 2016.

\bibitem{Dur_n_Olivencia_2019}
M.~A. Durán-Olivencia, R.~S. Gvalani, S.~Kalliadasis, and G.~A. Pavliotis.
\newblock Instability, rupture and fluctuations in thin liquid films: {T}heory
  and computations.
\newblock {\em J. Stat. Phys.}, 174(3):579–604, 2019.

\bibitem{FischerGruenFSOP}
J.~Fischer and G.~Grün.
\newblock {Finite speed of propagation and waiting times for the stochastic
  porous medium equation: A unifying approach}.
\newblock {\em SIAM J. Math. Anal.}, 47:825--854, 2015.

\bibitem{FischerGruen}
J.~Fischer and G.~Grün.
\newblock Existence of positive solutions to stochastic thin-film equations.
\newblock {\em SIAM J. Math. Anal.}, 50(1):411--455, 2018.

\bibitem{GessGnann2020}
B.~Gess and M.~V. Gnann.
\newblock The stochastic thin-film equation: {E}xistence of nonnegative
  martingale solutions.
\newblock {\em Stoch. Process. Appl.}, 130(12):7260--7302, 2020.

\bibitem{GiacOtto2002}
L.~Giacomelli and F.~Otto.
\newblock Droplet spreading: {I}ntermediate scaling law by {PDE} methods.
\newblock {\em Commu. Pure and Appl. Math.}, 55(2):217--254, 2002.

\bibitem{GiacShish2005}
L.~Giacomelli and A.~Shishkov.
\newblock Propagation of support in one-dimensional convected thin-film flow.
\newblock {\em Indiana Univ. Math. J.}, 54(4):1181--1215, 2005.

\bibitem{GruenKlein24Exist}
G.~Gr\"{u}n and L.~Klein.
\newblock Existence of nonnegative energy-dissipating solutions to a class of
  stochastic thin-film equations under weak slippage: {P}art {I} -- positive
  solutions.
\newblock {\em Submitted for publication}.

\bibitem{GruenKlein24FSOP}
G.~Gr\"{u}n and L.~Klein.
\newblock Finite speed of propagation for a class of stochastic thin-film
  equations.
\newblock {\em Submitted for publication}.

\bibitem{GruenKlein22}
G.~Gr\"{u}n and L.~Klein.
\newblock Zero-contact angle solutions to stochastic thin-film equations.
\newblock {\em J. Evol. Equ.}, 22(3):Paper No. 64, 37, 2022.

\bibitem{GruenBernisInequ}
G.~Grün.
\newblock On {B}ernis’ interpolation inequalities in multiple space
  dimensions.
\newblock {\em Z. Anal. Anwend.}, 20(4):987--998, 2001.

\bibitem{GruenDropletWTP}
G.~Grün.
\newblock Droplet spreading under weak slippage: The waiting time phenomenon.
\newblock {\em Ann. I. H. Poincaré}, 21:255–269, 2004.

\bibitem{MeckeRauscher}
G.~Grün, K.~Mecke, and M.~Rauscher.
\newblock Thin-film flow influenced by thermal noise.
\newblock {\em J. Stat. Phys.}, 122(6):1261--1291, 2006.

\bibitem{HofmanovaWeak2016}
M.~Hofmanova, M.~R\"oger, and M.~{von Renesse}.
\newblock Weak solutions for a stochastic mean curvature flow of
  two-dimensional graphs.
\newblock {\em Probab. Theory Relat. Fields}, 168(1-2):373--408, 2017.

\bibitem{HofmanovaSeidler2012}
M.~Hofmanova and J.~Seidler.
\newblock On weak solutions of stochastic differential equations.
\newblock {\em Stoch. Anal. Appl.}, 30:100--121, 2012.

\bibitem{Jakubowski1994}
A.~Jakubowski.
\newblock The almost sure {S}korokhod representation for subsequences in
  nonmetric spaces.
\newblock {\em Theory Probab. Appl.}, 42(1):167--174, 1998.

\bibitem{Karatzas2014}
I.~Karatzas and S.~Shreve.
\newblock {\em Brownian Motion and Stochastic Calculus}.
\newblock Springer, Berlin, Heidelberg, 1998.

\bibitem{Krylov}
N.~V. Krylov.
\newblock A relatively short proof of {It\^o's} formula for {SPDEs} and its
  applications.
\newblock {\em {Stoch. Partial Differ. Equ. Anal. Comput.}}, 1(1):152--174,
  2013.

\bibitem{Kuo2006}
H.-H. Kuo.
\newblock {\em Introduction to Stochastic Integration}.
\newblock Universitext. Springer, New York, 2006.

\bibitem{MetzgerGruen2022}
S.~Metzger and G.~Gr\"{u}n.
\newblock Existence of nonnegative solutions to stochastic thin-film equations
  in two space dimensions.
\newblock {\em Interfaces Free Bound.}, 24(3):307--387, 2022.

\bibitem{Sauerbrey2021}
M.~Sauerbrey.
\newblock Martingale solutions to the stochastic thin-film equation in two
  dimensions.
\newblock {\em Ann. Inst. Henri Poincar\'{e} Probab. Stat.}, 60(1):373--412,
  2024.

\bibitem{Sauerbrey2024solutions}
M.~Sauerbrey.
\newblock Solutions to the stochastic thin-film equation for the range of
  mobility exponents $n\in (2,3)$.
\newblock {\em arXive e-prints, arXiv:2310.02765}, 2024.

\end{thebibliography}
			
              \end{document}